\documentclass[11pt]{article}

\makeatletter
\let\@fnsymbol\@arabic
\makeatother

\usepackage{appendix}
\usepackage{graphicx}
\usepackage{makecell}
\usepackage{amsmath,amsfonts,amssymb}
\usepackage{epsfig}
\usepackage{dsfont}
\usepackage{mathrsfs}
\usepackage{bbm}
\usepackage{multirow}
\usepackage[dvipsnames,usenames]{xcolor}
\usepackage{subcaption}
\usepackage{caption}
\usepackage{svg}
\usepackage{natbib}
\usepackage{mathtools}
\usepackage{enumitem}
\usepackage{float}
\usepackage{lscape}
\usepackage{rotating}
\usepackage{wrapfig}
\usepackage{colortbl}
\usepackage{booktabs}
\usepackage{diagbox}
\usepackage{threeparttable}
\usepackage{adjustbox}
\usepackage{boldline}
\usepackage[colorlinks=true]{hyperref}
\hypersetup{urlcolor=blue,citecolor=blue,linkcolor=blue}

\allowdisplaybreaks[4]
\setcitestyle{numbers,square}

\renewcommand{\Box}{\framebox{\rule{0.3em}{0.0em}}}

\newtheorem{theorem}{Theorem}[section]
\newtheorem{theorem*}{Theorem}[subsubsection]
\newtheorem{lemma}{Lemma}[section]
\newtheorem{proposition}{Proposition}[section]
\newtheorem{example}{Example}[section]

\newtheorem{definition}{Definition}[section]

\newtheorem{assumption}{Assumption}[section]

\newcommand{\setd}{{ d \kern -.15em l}}
\newcommand{\hatsetd}{ d \hat{\kern -.15em l }}
\newcommand{\dd}{\mathsf {d\kern -0.07em l}}

\newcommand{\vt}{{\vartheta}}

\newcommand{\bgeqn}{\begin{eqnarray}}
\newcommand{\edeqn}{\end{eqnarray}}
\newcommand{\bgeq}{\begin{eqnarray*}}
\newcommand{\edeq}{\end{eqnarray*}}
\newcommand{\bec}{\begin{center}}
\newcommand{\enc}{\end{center}}
\newcommand{\R}{{\rm I\!R}}

\newcommand{\inmat}[1]{\mbox{\rm {#1}}}

\newcommand{\F}{{\cal F}}

\newcommand{\B}{{\cal B}}

\newcommand{\be}{\begin{equation}}
\newcommand{\ee}{\end{equation}}

\def\bbe{{\Bbb{E}}}

\renewcommand{\Box}{\hfill \rule{2.3mm}{2.3mm}}

\newcommand{\keywords}[1]{\par\vspace{0.5em}\noindent\textbf{Keywords: }#1\par}

\numberwithin{equation}{section}

\graphicspath{{./figures/}}

\title{Risk-averse Decision Making with Contextual Information: Model, Sample Average Approximation, and Kernelization}

\author{Yuan Tao\thanks{Department of Systems Engineering \& Engineering Management, The Chinese University of Hong Kong, Shatin, N.T., Hong Kong. Email: yuantao@se.cuhk.edu.hk}
\and Erick Delage\thanks{Department of Decision Sciences, HEC Montr\'eal, Qu\'ebec, Canada. Email: erick.delage@hec.ca}
\and Huifu Xu\thanks{Department of Systems Engineering \& Engineering Management, The Chinese University of Hong Kong, Shatin, N.T., Hong Kong. Email: hfxu@se.cuhk.edu.hk}}

\date{}

\begin{document}

\maketitle

\begin{abstract}
We consider risk-averse contextual optimization problems where the decision maker (DM) faces two types of uncertainties: problem data uncertainty (PDU) and contextual uncertainty (CU) associated with PDU, the DM makes an optimal decision by minimizing the risk arising from PDU based on the present observation of CU and then assesses the risk of the optimal policy against the CU. A natural question arises as to whether the nested risk minimization/assessment process is equivalent to joint risk minimization/assessment against CU and PDU simultaneously. First, we demonstrate that the equivalence can be established by appropriate choices of the risk measures and give counterexamples where such equivalence may fail. One of the interesting findings is that the optimal policies are independent of the choice of the risk measure against the CU under certain conditions. Second, by using the equivalence, we propose a computational method for solving the risk-averse contextual optimization problem by solving a one-stage risk minimization problem. The latter is particularly helpful in data-driven environments. We consider a number of risk measures/metrics to characterize the DM’s risk preference for PDU and discuss the computational tractability for the resulting risk-averse contextual optimization problem. Third, when the risk-averse contextual optimization problem is defined in the reproducing kernel Hilbert space {(RKHS)}, we show the consistency of the optimal values obtained from solving sample average approximation problems. Some numerical tests, on the newsvendor problem and portfolio selection problem, are performed to validate the theoretical results.
\end{abstract}


\section{Introduction} \label{sec:intro}

Contextual optimization is concerned with
optimal decision-making where the decision maker (DM)
faces two types of uncertainties.
One is problem data uncertainty (PDU), which is inherent in the decision-making problem
and directly
affects
the objective and/or constraints
of the problem.
  The other is
  contextual uncertainty (CU), which represents random contextual information
associated with PDU.
The optimal decision is made before the realization of PDU based on the present observation of CU.
For example, in inventory management problems,
the level of demand (PDU) depends on many correlated features, such as seasonality, weather, location, and economic indicators (CU),
which are available before the order is placed and can be used to infer the distribution of random demand \citep{ban2019big,elmachtoub2022smart}.
In portfolio selection,
it is believed that the daily asset returns (PDU)  are closely related to some side information (CU) such as Volatility Index (VIX), 10-year Treasury Yield Index,  Crude Oil Index, S\&P 500 (GSPC), Dow Jones Index (DJI), which are available
when assembling the portfolio
and could be used to infer the distribution of future random returns \citep{chenreddy2022data,nguyen2024robustifying,wang2022robust}.

In order to incorporate the contextual information into the decision-making problem,
Ban and Rudin \cite{ban2019big} employ a parameterized mapping from contextual information to inventory decision and select parameters to achieve the best empirical performance based on the available data.
This method, termed `decision rule approach', is also applied to the portfolio selection problem by Bazier-Matte and Delage \cite{bazier2020generalization} and
extended to general stochastic programming problems using reproducing kernel Hilbert spaces {(RKHS)}
\cite{bertsimas2022data} and the non-convex piecewise affine decision rule
\cite{zhang2023data}.
Zhang et al. \cite{zhang2024optimal} explore a distributionally robust framework for the newsvendor problem and identify a class of nonparametric policies that
interpolate an optimal in-sample policy to unobserved feature values in a specific manner.
{
Chen et al. \cite{chen2024robust} focus on interpretable policy classes, such as tree-based policies, and integrate a robust satisficing model to handle distributional ambiguity.
}
Another stream of studies and approaches is named sequential learning and optimization (SLO),
which
first trains models to predict the conditional distribution of
problem data
based on observed contextual information,
and subsequently solves the associated conditional stochastic programming.
Bertsimas and Kallus \cite{bertsimas2020predictive} propose various methods for predicting the conditional distribution of problem data,
while
Zhu et al. \cite{zhu2022joint} focus on the point prediction
and take a robust optimization perspective to address the inaccuracy in prediction.
To obtain a distribution prediction from the point prediction model,
Kannan et al. \cite{kannan2024residuals} add the residuals of the prediction model to the point prediction, thereby constructing a conditional distribution prediction.
Distributionally robust optimization techniques are employed to enhance the predicted conditional distribution \cite{kannan2021heteroscedasticity,kannan2025data}.
To promote predictions that facilitate sequential decision-making,
the integrated learning and optimization (ILO) framework is proposed. This framework studies a bilevel optimization problem comprising a learning phase and
a decision-making phase, see
\cite{donti2017task,kallus2023stochastic,qi2025integrated}.
Elmachtoub and Grigas \cite{elmachtoub2022smart} investigate this framework from a regret perspective.
Further extensions and theoretical analyses of this method are provided in \cite{el2023generalization,elmachtoub2023estimate,ho2022risk,iyengar2023optimizer,liu2021riskbound}.
For a comprehensive overview on contextual optimization, see \cite{sadana2024survey}.

These studies mainly focus on risk-neutral decision-making problems.
To
address the issue of risk aversion in contextual optimization,
researchers have explored risk-averse stochastic optimization methods and robust optimization methods on this topic.
With a predicted conditional distribution of problem data,
Lin et al. \cite{lin2022data} introduce a value-at-risk (VaR) constraint in the newsvendor problem to ensure that the profit meets a specific target with a high probability in each context.
Rahimian and Pagnoncelli \cite{rahimian2023data,rahimiancontextual} investigate the risk-averse conditional stochastic programming problem using chance constraints and expected-value constraints.
Kallus and Mao \cite{kallus2023stochastic} consider variance-based and conditional value-at-risk (CVaR) portfolio optimization problems within the ILO framework.
Nguyen et al. \cite{nguyen2024robustifying} focus on the conditional portfolio optimization problem with mean-variance and mean-CVaR objective functions, and construct a distributionally robust ambiguity set for the conditional distribution of uncertain returns.
Wang et al. \cite{wang2022robust} investigate the portfolio optimization where the uncertain returns and the contextual information jointly follow a Gaussian Mixture (GM) distribution
and adopt a robust perspective on the parameters of returns.
The conditional robust optimization framework is introduced in \cite{chenreddy2024end,chenreddy2022data}, where the uncertainty set of problem data is predicted using contextual information.
The above literature mainly
addresses the risk-averse decision-making problem using predicted conditional distributions within SLO and ILO frameworks.
Regarding the decision-making approach,
Bazier-Matte and Delage \cite{bazier2020generalization} introduce expected utility theory to the contextual portfolio optimization problem. They characterize the investor's risk attitude using a strictly increasing concave utility function and evaluate the performance of the investment policy based on the expected utility over the joint distribution of contextual information and asset returns.
Liu and Zhu \cite{liu2024newsvendor} approach the contextual newsvendor problem from a risk management perspective.
They use CVaR to evaluate the risk of a parameterized policy over the joint distribution of covariates and problem data.
{
This method of assessing CVaR directly on the joint distribution is also adopted by
Chen et al. \cite{chen2025robust} to address the portfolio optimization problem.
}

In this paper, we follow the strand of research by
considering that
the DM makes an optimal decision
by minimizing the risk
arising from PDU
based on the
observation of CU and then assesses the risk of the optimal policy against the CU.
Specifically, since CU is observable at the time of decision making while the PDU remains unobserved, the
DM primarily
optimizes the decision (called policy) to deal with the risk arising from PDU,
which is
a conditional risk minimization
as investigated by \cite{lin2022data,nguyen2024robustifying,wang2022robust}.
Since the optimal policy is dependent on the contextual information, then it may be viewed as a function mapping from contextual information to the
space of decision vectors.
Once the optimal policy is selected,
the DM evaluates the risk
of the decision-making problem (based on the optimal policy)
across contextual uncertainty.
For instance, in the portfolio optimization, an investor
makes portfolio decisions based on the risk arising from asset returns (PDU), which can be partially predicted using economic indicators (CU),
and
then evaluates the risk of
the optimal policy across the market by considering all potential outcomes of these economic indicators (CU).
A natural question
arises as to whether
the nested risk minimization/assessment
process is equivalent to joint risk minimization/assessment against
CU and PDU. This paper aims to address the question and a number of
related issues.

The main contributions of the paper can be summarized as follows.

\begin{enumerate}
    \item[(a)]
    \textbf{Modeling.}
Unlike the existing risk-averse contextual
optimization models, which are
composed of a learning phase and a decision phase, we propose a model where an optimal decision/policy is
selected to minimize the risk of PDU and the risk of the policy is then evaluated against the CU.
We cast
the decision-and-evaluation process as
a risk-averse stochastic programming problem with nested risk measures for CU and PDU,
and then find an optimal policy
that minimizes the nested risk.
We call this ex post risk minimization approach. Since the resulting optimization problem is
difficult to solve particularly in a data-driven environment, we propose a one-stage ex ante risk minimization model where the optimal policy is selected to minimize the risk of the joint uncertainty of contextual information and the problem data (a combination of PDU and CU).
    To illustrate the versatility of the proposed models,
     we consider
     a number of risk measures/metrics characterizing the DM's risk preference.
Particular attention is given to
the case when the risk measures are
CVaR, optimized certainty equivalent and entropic risk measures.

    \item[(b)]
    \textbf{Analysis.}
    We discuss
    well-definedness of both models
    (Proposition~\ref{prop:random_variable_inner_risk} and  Proposition~\ref{prop:well-definedness-problem})
    and
    derive sufficient conditions for the equivalence of the two models.
    Under some moderate conditions, we
    use the interchangeability principle \cite{shapiro2017interchangeability} to
    demonstrate that the optimal     policy
    to the ex ante risk minimization problem is also optimal
to the conditional risk minimization problem
in each context
    (Theorem~\ref{thm: optimality-of-ent}).
    An interesting consequence is that
    the optimal policy is independent of
the choice of
the risk measure  for the contextual risk,
which allows us to
take a risk-neutral
approach
(with an expectation operator) for
the CU.
To construct ex ante risk mappings,
we first consider the case where the DM's risk attitude towards PDU is characterized by optimized certainty equivalent such as CVaR and entropic risk measure.
By taking expectation
over the risk arising from CU,
we derive an appropriate ex ante risk mapping
which ensures optimality of the conditional risk minimization problem (Example~\ref{example:expected-oce} and  Example~\ref{example:expected-ent}).
Moreover, based on
Kusuoka's representation theorem \cite{kusuoka2001law},
we extend
the model and analysis to a number of risk measures/metrics including
mean-upper-deviation, spectral risk measure, distortion risk metrics and quantile deviation (Examples \ref{example:coherent}-\ref{example:quantile-semide}).

    \item[(c)]
    \textbf{Consistency}.
To facilitate the numerical solution of
the ex ante contextual minimization problem, we
    propose to
use sample average approximation (SAA) to discretize
PDU and CU and restrict the policy to reproducing kernel Hilbert spaces.
    The RKHS approach is employed in risk-neutral contextual optimization by \cite{bertsimas2022data,wang2026data,wang2024generalization}.
    Here, we complement this line of research
  by deriving
  convergence/consistency of the SAA problem in terms of optimal values as the sample
    size goes to infinity
 (Theorem~\ref{theorem:saa_convergence_rkhs}).
  In the case where the RKHS is generated by the universal kernel,
  the optimal value
  has the potential to provide an arbitrarily accurate approximation to the optimal value in the space of all continuous policies.
    The result ensures that the optimal value obtained by solving
    the SAA
    problem in RKHS converges to
    its true counterpart (Theorem~\ref{theorem:optimal-rkhs}).
    To examine the performance of the proposed risk-averse contextual optimization framework and approximation schemes,
    we have conducted some numerical
    tests on the newsvendor problem and portfolio selection problem.
    Our results confirm the established theory.

\end{enumerate}

{\color{black}

The rest of the paper is organized as follows.
Section \ref{sec:preliminary} recalls
some basic
 definitions of the
 topology of weak convergence,
 risk measures
 and three popular paradigms in contextual optimization.
Section~\ref{sec:contextual-risk-measure}
formally establishes two risk-averse contextual optimization models and investigates their well-definedness.
Section~\ref{sec:contextual-ra-opt}
examines the properties of the proposed models.
Many examples are used to illustrate the models and their properties.
Section \ref{sec:DR-rkhs}
examines the sample average approximation of the risk-averse contextual optimization problem in the reproducing kernel Hilbert space and
investigates the
convergence of the optimal value
as the sample size increases.
Section \ref{sec:experiment} reports numerical test results
on the proposed
risk-averse contextual optimization models.
Section \ref{sec:extension} outlines directions for potential extensions.
Finally,
Section \ref{sec:conclusion} concludes with some remarks.
}

\section{Preliminaries} \label{sec:preliminary}

In this section, we
recall some basic notions and results regarding
risk measures and contextual stochastic optimization.
\subsection{Topology of weak convergence and risk measures}
\label{sec:risk-measure}

Consider a probability space $(\Omega, \mathcal{F}, \mathbb{P})$ and
a $\mathcal{F}$-measurable function $\xi: \Omega \to \R^{d}$.
Let $P_\xi=\mathbb{P}\circ \xi^{-1}$
denote
the probability measure/distribution\footnote{Throughout the paper,
we use terms probability measure and probability distribution interchangeably.
We also use $\xi$ for both random variable and its realization depending on the context.} on $\R^{d}$ induced by $\xi$,
and
$L_{d}^p$ denote the
linear space of $d$-dimensional random vectors defined over $(\Omega, \mathcal{F}, \mathbb{P})$
with $p$-th order finite moments
$$
\mathbb{E}_{\mathbb{P}}[\|\xi\|^p] = \int_\Omega \|\xi(\omega)\|^p \mathbb{P}(d\omega)= \int_{\R^{d}} \|t\|^p \mathbb{P}\circ \xi^{-1}(dt) = \int_{\R^{d}} \|t\|^p P_\xi(dt)<\infty,
$$
{
where $\|\cdot\|$
denotes a norm in $\R^d$.
Unless  specified otherwise,
it is treated as
$\|\cdot\|_2$.
}
Let $\mathscr{P}(\R^{d})$ denote the space of all probability measures over $(\R^{d},{\cal B}(\R^{d_\zeta}))$, {
where $\mathcal{B}(\R^{d_\zeta})$ denotes the Borel $\sigma$-algebra on $\R^{d_\zeta}$,
}
and
\begin{eqnarray*}
        {\cal M}^p_{d} =\left\{P\in \mathscr{P}(\R^{d}):
\int_{\R^{d}} \|t\|^p P(dt)<\infty\right\}.
\end{eqnarray*}
Define {${\cal C}_d^p$ the linear space of all continuous functions $h:\R^{d}\to \R$} for which there exists a positive constant $c_0$ such that
$
   |h(t)|\leq c_0(\|t\|^p+1), \forall t\in \R^{d}.
\label{eq:ht-phi-topo}
$
The $\|\cdot\|^p$-weak topology, denoted by $\tau_{\|\cdot\|^p}$, is the coarsest topology on
${\cal M}_{d}^p$ for which the mapping $g_h:{\cal M}_{d}^p\to \R$ defined by
$
g_h(P) :=\int_{\R^{d}} h(t) P(dt),\; h\in {\cal C}_{d}^p
$
 is continuous. A sequence  $\{P_k\} \subset {\cal M}_d^p$ is said to converge $\|\cdot\|^p$-weakly to $P\in {\cal M}_d^p$ written
 ${P_k} \xrightarrow[]{\|\cdot\|^p} P$ if it converges w.r.t. $\tau_{\|\cdot\|^p}$.
 In the case that $p=0$, {i.e., $\mathcal{C}_d^0$ is the linear space of bounded continuous functions,} $\tau_{\|\cdot\|^0}$ reduces to the ordinary topology of weak convergence.
 {
Since $\mathcal{C}_d^0 \subset \mathcal{C}_d^p$ for every $p \geq 0$,
the $\|\cdot\|^p$-weak topology on $\mathcal{M}_d^p$ is finer than the ordinary weak topology on $\mathcal{M}_d^0$.
In particular, $\|\cdot\|^p$-weak convergence implies ordinary weak convergence.
In the general definition, we do not require $p \geq 1$.
To ensure that our results apply to a broader range of settings, we impose the condition $p\geq 1$ in some results.
Specific assumptions will be made in the forthcoming discussions depending on the context.
}
 We refer readers to \cite{claus2016advancing} for comprehensive treatments of these concepts.
 The topology of weak convergence effectively poses restrictions
 on probability distributions of the associated random variables and their data structure.
 We need the concept
 for well-definedness of
the risk-averse optimization models to be introduced  in the forthcoming discussions.

 Recall that a function $\rho(\cdot):  L^p\rightarrow \R$ is
 called a \textit{monetary  risk measure} if it satisfies:
 (a) \textit{monotonicity}: for any $\xi,\zeta \in
 L^p
 $,
$\xi(\omega) \geqslant \zeta(\omega)$
almost surely
implies that $\rho(\xi) \geqslant \rho(\zeta)$, and
(b) \textit{cash invariance}: $\rho(\xi+c)=\rho(\xi)+c$ for any $\xi \in L^p$ and real number $c \in \R$. A monetary risk measure $\rho(\cdot)$ is said to be \textit{ convex}
if it also satisfies:
(c) \textit{convexity}: for any $\xi, \zeta \in L^p$ and $\lambda \in[0,1],$   $\rho(\lambda \xi+(1-\lambda)\zeta) \leqslant \lambda\rho(\xi)+(1-\lambda)\rho(\zeta)$.
A convex risk measure $\rho(\cdot)$ is said to be \textit{ coherent}
 if it further satisfies
(d) \textit{positive homogeneity}: $\text{for any }\xi \in L^p \text{ and any } \lambda \geqslant 0, \text{ it holds that } \rho(\lambda \xi) = \lambda\rho(\xi).$
A risk measure $\rho$ is (e) \textit{ law invariant}, if
$\rho(\xi)=\rho(\zeta)$
for any random variables $\xi,\zeta\in L^p$ with the same probability distribution.
We need a property slightly stronger than
monotonicity.
For any
two random variables $\xi,\zeta\in L^p$, we write $\xi \succ \zeta$ if $\xi(\omega)
\geq \zeta(\omega)$ almost surely
and $\mathbb{P}(\omega\in\Omega: \xi(\omega) > \zeta(\omega))>0$.
A risk measure $\rho:L^p\to \R$ is said to be
\textit{strictly monotone}, if
$\rho(\xi)>\rho(\zeta)$
for any $\xi,\zeta \in L^p$
with $\xi \succ \zeta$.
This property is discussed in Chapter 6.3.4 of \cite{shapiro2014lectures}.

It may be helpful to note that any law invariant
risk measure defined over a non-atomic $L^p$ space can be represented as a risk functional over the space of the probability distributions of the random variables.
Specifically, there exists a unique
functional $\varrho: {\cal M}_1^p\to
\R$ such that
\begin{eqnarray*}
    \rho(\xi) = \varrho(P_\xi) := \rho\left(F_{P_\xi}^{-1}(U)\right),
\end{eqnarray*}
where $U$ is a random variable uniformly distributed over $[0,1]$, $F_{P_\xi}$ denotes the cumulative distribution
function associated with $P_\xi$ and $F_{P_\xi}^{-1}(p):=\inf\{t\in\R:F_{P_\xi}(t)\geq p\}$ is a quantile function of $p\in[0,1]$.
{
In the forthcoming discussions, we
often need to indicate the dependence
of $\rho$ on $P_\xi$ by writing
$\rho_{P_\xi}$.}
Consider a random function $f:\R^{d_z}\times \R^{d_\zeta}$ with
$|f(z,\zeta)|\leq \eta(z)(\|\zeta\|^\iota+1), \forall z\in \R^{d_z},$ for some constant $\iota>0$,
and
a probability measure $Q$ over $(\R^{d_\zeta},\mathcal{B}(\R^{d_\zeta}))$.
Let $\theta(z,Q)$ denote the
image measure of $\delta_z\otimes Q$ under
$f(z,\zeta)$, i.e., $\theta(z,Q):=(\delta_z\otimes Q)\circ f^{-1}$,
where $\delta_z\otimes Q$ denotes the product probability measure of the Dirac probability measure $\delta_z$ at $z$ and $Q$.
Then $\theta(z,Q)$ is a probability measure on $(\R,\mathcal{B}(\R))$ induced by $f$, and
$\theta(z,Q)\in {\cal M}^p_1$
if $Q\in {\cal M}_{d_\zeta}^{\iota p}$
given the relationship
$
\int_{\R} |t|^p \theta(z,Q)(dt)=\int_{\R^{d_z}\times \R^{d_\zeta}} |f(z',\zeta)|^p (\delta_z
\otimes Q)(d(z',\zeta)) = \int_{\R^{d_\zeta}}|f(z,\zeta)|^pQ(d\zeta),\forall z\in \R^{d_z},
$
see Corollary 2.11 in \cite{claus2016advancing}. Consequently, we have
$\rho(f(z,\zeta))
=\varrho(\theta(z,Q))$.

In this paper, we will frequently use four well-known
classes of
risk measures:
optimized certainty equivalent (OCE),
conditional value-at-risk (CVaR, also termed expected shortfall),
entropic risk measure,
and shortfall risk measure.
The OCE
of a random loss $\xi$ is defined as
$
S_u(\xi) := \sup\limits_{s\in \R}\left\{s + \mathbb{E}[u(-\xi-s)]\right\},
$
where $u:\R\rightarrow\R$ is a proper closed concave and non-decreasing utility function.
By considering the negative of OCE,
we
obtain
a risk measure
\bgeqn
\label{eq:neg-OCE}
\rho_u(\xi) := -S_u(\xi) = \inf\limits_{s\in \R} \left\{-s - \mathbb{E}[u(-\xi-s)]\right\}
\edeqn
via particular choices of a utility function.
For instance,
by setting $u(x) = \frac{1}{1-\beta}x$ for $x \leq 0$ and $u = 0$ for $x > 0$,
we obtain \begin{eqnarray} \label{def:cvar}
    \inmat{CVaR}_{\beta}(\xi) := \inf\limits_{t\in\R}\left\{ t+\frac{1}{1-\beta}\mathbb{E} [(\xi-t)_+] \right\},
\end{eqnarray}
where $(x)_+ = \max\{x,0\}$, and { $\beta \in (0,1)$ is a specified probability level, typically selected close to $1$.
$\text{CVaR}_\beta$
represents the average
of the worst $(1-\beta)\%$ of the outcomes.}
On the other hand,
by setting
$u(x) = \frac{1}{\gamma}\left(1-e^{-\gamma x}\right)$,
we obtain
entropic risk measure
\begin{eqnarray}
\label{def:entropic-risk-measure}
    \rho_{\gamma}^{\text{ent}}(\xi) := \frac{1}{\gamma} \ln \mathbb{E}[e^{\gamma \xi}],
\end{eqnarray}
where $\gamma> 0$ is a positive parameter; see \cite{ahmadi2012entropic}. We refer readers to \cite{ben1986expected,ben2007old} for the original work on OCEs.

{
Let $l:\R\rightarrow \R$ be an increasing and non-constant loss function.
With a given acceptable loss $l_0\in\R$,
the loss function $l$ induces the shortfall risk measure defined by
\begin{eqnarray} \label{eqn:shortfall-rm}
    \rho_l(\xi):= \inf\left\{s\in\R|\mathbb{E}[l(\xi-s)]\leq l_0\right\}.
\end{eqnarray}
Many convex risk measures, including OCE, CVaR, and entropic risk measure, can be represented in the following parametric form:
\begin{eqnarray} \label{eqn:parametric-rm}
    \rho(\xi):=\inf_{s\in \mathcal{S}} \mathbb{E}[U(\xi,s)],
\end{eqnarray}
where $\mathcal{S}$ is a subset of a finite-dimensional vector space and $U:\R\times \mathcal{S}\rightarrow \R$ is a proper continuous function that is non-decreasing in the first attribute and jointly convex in both attributes;
see \cite{guigues2023risk} for a detailed discussion and more examples.
In the forthcoming discussions,
we will write $\rho_{u,P_\xi}$, $ \inmat{CVaR}_{\beta, P_\xi}$, $\rho_{\gamma,P_\xi}^{\text{ent}}$, and $\rho_{l,P_\xi}$
for $\rho_{u}$, $\inmat{CVaR}_{\beta}$, $\rho_{\gamma}^{\text{ent}}$, and $\rho_{l}$ to emphasize that the underlying expectation is taken w.r.t.~the probability distribution of $\xi$.
}

\subsection{Contextual stochastic optimization} \label{sec:contextual-sto-opt}

Consider a one-stage
stochastic
programming
problem:
\begin{eqnarray}
    \min_{z\in \mathcal{Z}} \mathbb{E}_{P_Y}[c(z, Y)],
    \label{eq:one-stg-SP}
\end{eqnarray}
where $z$ is a decision vector with $\mathcal{Z}$ being a nonempty subset of $\R^{d_z}$,
$Y:\Omega\rightarrow \mathcal{Y}\subseteq \R^{d_y}$
is a random vector defined on probability space $(\Omega,\mathcal{F},\mathbb{P})$,
$c: \R^{d_z}\times \R^{d_y} \rightarrow \R$ is a cost function,
and the mathematical expectation $\mathbb{E}$ is taken w.r.t. the probability distribution $P_Y=\mathbb{P}\circ Y^{-1}$.
The uncertain problem data $Y$ is unknown when making the decision,
but a covariate $X:\Omega\rightarrow \mathcal{X}\in \R^{d_x}$, which is correlated with the uncertain problem data $Y$, is observed before decision making.
The joint distribution of covariate $X$ and problem data $Y$ is denoted by $P:=P_{XY}$.
Specifically, given the
contextual information $x\in \mathcal{X}$ and
conditional distribution $P_{Y|X}$,
a risk-neutral
decision maker (DM) is interested in
finding an optimal
decision
$z^*\in \mathcal{Z}$
by solving the following conditional
stochastic optimization problem:
\begin{equation*}
    \min_{z\in \mathcal{Z}} \mathbb{E}_{P_{Y|X=x}}[c(z,Y)].
\end{equation*}
However, the relationship between covariate $X$ and problem data $Y$
is usually unknown in practice,
and it is difficult to figure out the conditional distribution $P_{Y|X}$
from the sample data $\{(x_i,y_i)\}_{i=1}^N$.
Therefore,
in the contextual optimization problem, the DM needs to find out the optimal
decision
while exploring the relationship between $X$ and $Y$.
This is a departure from
the
classical one-stage
stochastic programming model (\ref{eq:one-stg-SP}).
There are
three popular approaches
in the literature of  contextual optimization.
We refer readers to \cite{sadana2024data} for a detailed review on contextual optimization.

\textbf{
Decision rule approach}.
In this approach,
the relationship between covariate $X$ and the optimal
decision $z^*$ is characterized by a policy $g: \mathcal{X} \rightarrow \mathcal{Z}$,
and
the learning phase and decision-making phase
are incorporated into
a
one-stage
stochastic optimization framework.
By selecting a hypothesis space $\widetilde{\mathcal{G}}$ for the policies,
the contextual optimization problem can be formulated as
\begin{eqnarray*}
    \min_{g\in \widetilde{\mathcal{G}}} \mathbb{E}_{P_X}\left[\mathbb{E}_{ P_{Y|X}}[c(g(X),Y)]\right]
    =
    \mathbb{E}_{P}\left[c(g(X),Y)\right],
\end{eqnarray*}
where $P$ is the joint probability distribution of $(X,Y)$ and the equality
is due to the tower property of the expectation operator $\mathbb{E}$.
In the setting of contextual optimization, $P$ is unknown but
can be approximated by the empirical distribution $P^N(\cdot) = \frac{1}{N} \sum_{i=1}^N \mathds{1}_{(x_i,y_i)}(\cdot)$
based on sample data $\{(x_i,y_i)\}_{i=1}^N$.
Here and later on, we use $\mathds{1}_{S}(\cdot)$ to denote
indicator function of a set $S$ or a point.

\textbf{Sequential learning and optimization (SLO)}.
In the SLO framework,
the relationship between $X$ and $Y$
is learned from
sample data $\{(x_i,y_i)\}_{i=1}^N$
via a mapping $f$ which maps $x\in \mathcal{X}$ to a distribution supported on $\mathcal{Y}$.
After observing
the realization of the covariate
$X=x$, the conditional distribution of problem data $Y$ given $x$
is predicted by
$\hat{P}_{Y|X=x} = f(x)$
and the optimal decision is
then made
by solving
optimization problem:
\begin{eqnarray*}
   \inmat{(SLO)}\quad  \min_{z\in \mathcal{Z}} \mathbb{E}_{f(x)}[c(z,Y)].
\end{eqnarray*}
When the cost function $c(z,y)$ is linear in $y$, the
distribution estimation reduces to point estimation, i.e., $f:\mathcal{X}\rightarrow\mathcal{Y}$ and the optimal decision is made by solving:
$
    \min_{z\in \mathcal{Z}} c(z,f(x)).
$
Regression models in machine learning
can be used in the learning phase.

\textbf{Integrated learning and optimization (ILO)}.
ILO consists of a learning phase and a decision-making phase.
Compared with the SLO framework, ILO
evaluates the performance of the prediction model $f$ based on the quality of resulting decisions rather than on the precision of prediction models. It
combines the learning phase
at the upper level
and the decision-making phase
at the lower level
into a
bi-level stochastic program:
    \begin{eqnarray*}
     (\inmat{ILO})\qquad
        \min_{f\in\mathcal{F}}
        && \mathbb{E}_{P_X}\left[\mathbb{E}_{P_{Y|X}}[L(z^*(X,f),X,Y)]\right] = \mathbb{E}_{P} [L(z^*(X,f),X,Y)]\nonumber\\
        \inmat{s.t.} && z^*(x,f) \in \arg\min_{z\in \mathcal{Z}} \mathbb{E}_{f(x)}[c(z,Y)],\; \forall x \in \mathcal{X},
    \end{eqnarray*}
    where $\mathcal{F}$ is the hypothesis
     space
     of the prediction model,
     $z^*(x,f)$ is the optimal decision
     based on
     the predicted distribution $f(x)$,
     and
     $L(z^*(x,f),x,y)$
     represents the loss incurred
     from
     the decision $z^*(x,f)$ with contextual information $X=x$ and problem data $Y=y$.
     The form of $L$ depends on
     actual problems, for instance,
    Elmachtoub et al. \cite{elmachtoub2022smart}
    consider a
    regret-based loss function.

\section{Risk management in contextual optimization}
\label{sec:contextual-risk-measure}

In this paper, we investigate
contextual optimization from a risk management perspective.
We consider a DM who aims to minimize the risk of random losses/costs associated with
random problem
parameter/data
$Y$ based on the present observation of the contextual uncertainty $X=x$, and then adopts a risk measure
to evaluate the performance of the minimal risk (the optimal value) of the decision-making problem
across contexts.

\subsection{Conditional risk minimization and risk evaluation}
\label{sec:3.1}

Let $x$ be an
observation of $X$.
The DM's conditional risk minimization problem is defined as:
\begin{eqnarray} \label{eqn:conditional_risk_min}
   \psi(x) := \min_{z\in\mathcal{Z}}\rho_{P_{Y|X=x}}^{(2)}(c(z,Y)),
\end{eqnarray}
where $\rho_{P_{Y|X=x}}^{(2)}:L^p\rightarrow \R$ is a risk measure/function which captures the DM's risk preferences and
$P_{Y|X=x}$ is the conditional distribution of $Y$ on
the observation $X=x$,
$L^p$ denotes the space of
 random variables $\phi:{\cal Y}\to \R$
with $p$-th order moments, i.e.,
$\int |\phi(y)|^p P_{Y|X=x}(dy)<\infty$ for almost every $x\in {\cal X}$ according to $\mathbb{P}$.
Let
$\theta(z,P_{Y|X=x}) := \delta_z\otimes P_{Y|X=x}\circ c^{-1}$
denote the image distribution of $\delta_z\otimes P_{Y|X=x}$ under $c(z,y)$ conditional on $X=x$. Then
$
\rho_{P_{Y|X=x}}^{(2)}(c(z,Y))
=\varrho(\theta(z,P_{Y|X=x})),
$
for some distribution-based (law invariant) risk measure
as discussed in Section \ref{sec:risk-measure}.
{
Throughout the paper, we make a blanket assumption that the feasible set $\mathcal{Z}$ is nonempty and closed, and that the minimum
in \eqref{eqn:conditional_risk_min} is attainable,
which means that the set of optimal solutions is non-empty.
We also assume that $\mathcal{X},\mathcal{Y},\mathcal{Z}$ are nonempty Polish subspaces of $\R^{d_x},\R^{d_y},\R^{d_z}$ endowed with the relative topology inherited from the Euclidean spaces.
This assumption covers a range of sets that typically arise in applications, such as open sets, closed sets, and intervals/rectangles.
}

Let $z^*(x)$ denote an optimal solution of the problem.
Since it depends on $x$, we may regard it as a vector-valued function mapping from ${\cal X}$ to $\mathcal{Z}$
provided the optimal solution is unique for each fixed $x\in {\cal X}$.
Let $\psi(x)$ denote the optimal value of problem (\ref{eqn:conditional_risk_min}). Then $\psi(x)$ is the minimal risk when $X=x$ and is attained if the DM takes decision $z^*(x)$. Since
the decision is merely based on a single observation of
the contextual uncertainty, the DM may wish to evaluate
the risk of the decision by looking into
\begin{eqnarray*}
   \rho_{P_X}^{(1)}(\psi(X))
=
   \rho_{P_X}^{(1)} \left(\rho_{P_{Y|X}}^{(2)}(c(z^*(X),Y))\right),
\end{eqnarray*}
where
 $\rho_{P_X}^{(1)}:L^p\rightarrow \R$ is a
 law-invariant
 risk measure/function.
Since
\begin{eqnarray*}
\label{eqn:ex_post_opt-a}
    \rho_{P_X}^{(1)} \left(\rho_{P_{Y|X}}^{(2)}(c(z^*(X),Y))\right) =
    \rho_{P_X}^{(1)} \left(\min_{z\in\mathcal{Z}}\rho_{P_{Y|X}}^{(2)}(c(z,Y))\right),
\end{eqnarray*}
 and
under some appropriate conditions, we may apply
the well-known interchangeability principle (see e.g.,~\cite{shapiro2017interchangeability}) to establish
\begin{eqnarray} \label{eqn:ex-post-interchange}
\rho_{P_X}^{(1)} \left(\min_{z\in\mathcal{Z}}\rho_{P_{Y|X}}^{(2)}(c(z,Y))\right) =
\min_{g\in\widetilde{\mathcal{G}}}\rho_{P_X}^{(1)}\left(\rho_{P_{Y|X}}^{(2)}(c(g(X),Y))\right),
\end{eqnarray}
then evaluation of the risk of $\psi(X)$ is equivalent to
solving the minimization problem on the rhs of (\ref{eqn:ex-post-interchange}). Unfortunately this problem is not easy to solve
in that the conditional distribution $P_{Y|X}$ is challenging to approximate in a data-driven setting.
For instance, if the sample data $\{(x_i,y_i)\}_{i=1}^N$ are obtained from a continuous distribution $P$ of $(X,Y)$, then
with probability one, we will have $x_k \neq x_j, \forall k \neq j$, see \cite{zhang2024optimal}. Consequently,
the objective function in problem ({\ref{eqn:ex-post-interchange}}) cannot be effectively
approximated as the sample size increases.
This motivates us to consider an entirely different risk minimization problem:
\begin{eqnarray} \label{eqn:ex_ante}
    \min_{g\in \widetilde{\mathcal{G}}} \rho_P^{\text{EA}}\left(c(g(X),Y)\right),
\end{eqnarray}
where $\rho_P^{\text{EA}}:L^p\rightarrow \R$ is a law-invariant risk mapping evaluated using the joint distribution of $(X,Y)$ directly.
Problem (\ref{eqn:ex_ante}) is much easier
to solve (we will come back to this in Section \ref{sec:DR-rkhs}).
The remaining question to address  is under what conditions
problem (\ref{eqn:ex_ante}) is equivalent to the problem
on the rhs of (\ref{eqn:ex-post-interchange}).
The answer depends on the structure of
$\rho_{P_X}^{(1)}$ and $\rho_{P_{Y|X}}^{(2)}$. We will give some simple examples where $\rho_P^{\text{EA}}$ may not exist.
Likewise, the interchangeability
in (\ref{eqn:ex-post-interchange}) also depends on the structure of $\rho_{P_X}^{(1)}$.

We call
$\rho_{P_{X}}^{(1)}\circ
\rho_{P_{Y|X}}^{(2)}$
\textit{ex post risk mapping}
to
emphasize that
the
risk arising from
PDU is based on an observation of
the CU and
then the
conditional risk is assessed
over the CU.
Likewise,
we
call $\rho_P^{\text{EA}}$
\textit{ex ante risk
mapping}
to emphasize that
overall risk
is assessed based on the joint distribution of the
PDU and the CU simultaneously.
{
Throughout the paper, we call problem \eqref{eqn:conditional_risk_min} the conditional risk minimization problem, where we minimize the risk from problem data $Y$ conditional on observed data $X$.
We call both problem \eqref{eqn:ex-post-interchange}, based on ex post risk, and problem \eqref{eqn:ex_ante}, based on ex ante risk, the risk-averse contextual optimization problems, where we minimize risk before observing context $X$.
}

In the case that
$\rho_P^{\text{EA}} =\inmat{CVaR}_{\beta,P}$,
(\ref{eqn:ex_ante}) recovers
the risk-averse contextual optimization
model
considered by \cite{chen2025robust} and \cite{liu2024newsvendor}:
\begin{eqnarray} \label{eqn:ex-ante-cvar}
    \min_{g\in \widetilde{\mathcal{G}}} \inmat{CVaR}_{\beta,P}\left(c(g(X),Y)\right).
\end{eqnarray}
However, the
optimal policy obtained
from solving
problem (\ref{eqn:ex-ante-cvar}) is not necessarily
optimal for
the conditional risk minimization problem  (\ref{eqn:conditional_risk_min})
with
$\rho_{P_{Y|X=x}}^{(2)} = \inmat{CVaR}_{\beta,
P_{Y|X=x}}$ for $x\in\mathcal{X}$ almost surely, as discussed in Section 2.2 in \cite{marzban2023deep}.
We will come back to this in
Section \ref{sec:model-crao}.
The
next example illustrates
the difference and the relationship between (\ref{eqn:ex-post-interchange}) and (\ref{eqn:ex_ante}).

\begin{example} [Shortfall risk]
\label{example:shortfall-risk}
Consider
model (\ref{eqn:ex_ante})
using shortfall risk, defined in (\ref{eqn:shortfall-rm}), i.e.,
\begin{eqnarray*}
    \rho_P^{\text{EA}}(c(g(X),Y)) = \min\left\{s^{(0)}\in \R\big|\;\mathbb{E}_P [l^{(0)}(c(g(X),Y)-s^{(0)})]\leq l^{(0)}_0\right\},
\end{eqnarray*}
where $l^{(0)}$ is an increasing, non-constant loss function.
Under some
appropriate conditions,
$\rho_P^{\text{EA}}(c(g(X),Y))$ is well-defined, see Lemma 2.2 in \cite{delage2022shortfall}.
The shortfall risk
calculates the minimum capital $s^{(0)}$ to be injected such that the expected loss of the new position
$c(g(X),Y)-s^{(0)}$
does not exceed the specified level $l_0^{(0)}$.
However, this model does not differentiate the PDU and CU, and therefore it is unable to interpret its relationship to the conditional risk arising from PDU in specific contexts.
In contrast,
we consider model (\ref{eqn:ex-post-interchange}) with
\begin{eqnarray*}
    \rho_{P_{Y|X=x}}^{(2)}(c(g(x),Y)):= \min\left\{s^{(2)}\in \R\big|\; \mathbb{E}_{P_{Y|X=x}} [l^{(2)}(c(g(x),Y)-s^{(2)})]\leq l^{(2)}_0\right\},
\end{eqnarray*}
being a shortfall risk,
where $l^{(2)}$ is an increasing, non-constant and convex loss function.
In this formulation, the inner risk $\rho_{P_{Y|X=x}}^{(2)}(c(g(x),Y))$
is the minimum capital to be injected in
the context $x$ such that the expected utility loss
falls within the specified level $l^{(2)}_0$,
and $\rho_{P_X}^{(1)}$ integrates the required amount of cash across all contexts.
If we choose $\rho_{P_X}^{(1)} = \mathbb{E}_{P_X}$,
then (\ref{eqn:ex-post-interchange}) evaluates the expected amount of cash that needs to be injected prior to deploying the policy $g$ in order
to cover the cash requirements across all contexts.
If we choose $\rho_{P_X}^{(1)}$ as value-at-risk with $\beta = 95\%$,
then (\ref{eqn:ex-post-interchange})  evaluates the amount of cash that needs to be injected in order
to cover the cash requirements in $95\%$ of contexts.
Specifically, {if we choose $l^{(2)}(t)=e^{\gamma t}$ and $l^{(2)}_0:=l^{(2)}(0) = 1$,}
then $\rho_{P_{Y|X=x}}^{(2)}
$ is an entropic risk measure (see \cite{guo2019distributionally}).
Then we can reformulate model (\ref{eqn:ex-post-interchange}) into model (\ref{eqn:ex_ante}) by setting $\rho_{P_X}^{(1)}$ as the entropic risk measure.
The reformulation is derived from the definition of the entropic risk measure:
\begin{align*}
\rho_{P_X}^{(1)}\left(\rho_{P_{Y|X}}^{(2)}(c(g(X),Y))\right)
&= \frac{1}{\gamma}\ln\mathbb{E}_{P_X}\left[e^{\gamma (1/\gamma)\ln\mathbb{E}_{P_{Y|X}}\left[e^{\gamma c(g(X),Y)}\right]}\right]\\
&= \frac{1}{\gamma}\ln\mathbb{E}_{P_X}\left[\mathbb{E}_{P_{Y|X}}\left[e^{\gamma c(g(X),Y)}\right]\right]\\
&= \frac{1}{\gamma}\ln\mathbb{E}_{P}\left[e^{\gamma c(g(X),Y)}\right]=:\rho_P^{\text{EA}}\left(c(g(X),Y)\right).
\end{align*}

\end{example}
{
Note that the nature of
the nested risk measure here
differs from the Bayesian composite risk measures in \cite{shapiro2023bayesian,lin2022bayesian,ma2024bayesian,ma2026adaptive}.
In Bayesian DRO and Bayesian risk MDP/SOC,  the underlying PDU is represented by a parametric distribution with some unknown parameters, and the distribution of the parameters is learned via a Bayesian formula. The outer risk measure is used to capture the epistemic uncertainty represented by the Bayesian posterior distributions of the unknown parameters.
By contrast, in our risk-averse contextual optimization
framework, the outer risk arises from  aleatoric contextual uncertainty.

}

\subsection{Well-definedness of risk-averse contextual optimization problem}

Before proceeding to further discussions, we need to verify the well-definedness of problem (\ref{eqn:ex-post-interchange}).
To facilitate the discussion, we make the following assumptions.

\begin{assumption} \label{assumption:well-defined_rho}
    Consider models (\ref{eqn:ex-post-interchange}) and \eqref{eqn:ex_ante}.
    (a) The risk measures {$\rho^{(1)}$, $\rho^{(2)}$,} $\rho^{\text{EA}}$ are law-invariant convex risk measures.
    (b) $\widetilde{\mathcal{G}}$ is the  set of all continuous
    functions
    $g:\mathcal{X}\rightarrow \mathcal{Z}$.
    (c) The cost function $c$ is continuous in $(z,y)$. There
    exist a positive parameter $\iota \geq 0$ and     a locally bounded
    non-negative function
    $\eta:\mathcal{Z}\rightarrow\R$ such that
    $
    |c(z,y)|\leq \eta(z) (\|y\|^\iota + 1),
    \forall (z,y)\in \mathcal{Z}\times \mathcal{Y}.$
    (d) {For a fixed $p\geq 1$,}
    the conditional distribution $P_{Y|X=x}$ of $Y$ is continuous in $x\in \mathcal{X}$ with respect to the $\|\cdot\|^{\iota p}$-weak topology.
\end{assumption}

These conditions are often  used in the literature and cover several common cases.
Condition (a) allows us to consider a group of convex risk measures, e.g., CVaR and entropic risk measure.
Condition (b) restricts the discussion to continuous policies.
Conditions (c) and (d)
{ ensure that
$c(z,Y)\in L^p$ for all $P_{Y|X=x}$ and all $z$ because
\begin{eqnarray*}
    \int |c(z,y)|^{p} P_{Y|X=x}(dy)&\leq& \eta(z)^p \int (\|y\|^{\iota}+1)^p P_{Y|X=x}(dy)\\
    &\leq& \eta(z)^p \int 2^{p-1}(\|y\|^{\iota p}+1^p) P_{Y|X=x}(dy)\\
    &=&2^{p-1}\eta(z)^p\int \|y\|^{\iota p} P_{Y|X=x}(dy) + 2^{p-1}\eta(z)^p\\
    &< &\infty.
\end{eqnarray*}
They also guarantee}
the continuity of $\rho^{(2)}_{P_{Y|X=x}} ( c(g(x),Y)): \mathcal{X}\rightarrow \R$ in $x$ (see Proposition~\ref{prop:random_variable_inner_risk}).
Condition (c) {specifically stipulates
} the growth condition of the cost function in $y$ for each fixed $z$,
{which ensures the finiteness of the $p$-th moment of the cost function, see e.g.~\cite{claus2016advancing}.}
This condition is satisfied
when $c$ is H\"older continuous in $y$.
{
In the newsvendor problem
and the portfolio selection problem
to be studied in
Section~\ref{sec:experiment},
the function $c$ is Lipschitz continuous.
Condition (d) requires that the conditional distribution of $Y$ change smoothly
as $x$ varies.
To see this, let $r:=\iota p$ and denote by $\mathscr{P}_r(\mathbb R^{d_y})$ the set of probability measures with a finite $r$-th moment. On $\mathscr{P}_r(\mathbb R^{d_y})$, the $\|\cdot\|^{r}$-weak topology can be generated by the Wasserstein distance $W_r$ of order $r$ defined by
\begin{eqnarray*}
    W_r(P,Q):= \inf_{\pi\in \Pi(P,Q)} \left(\int \|\xi-\zeta\|^r \pi(d\xi,d\zeta) \right)^{1/r},
\end{eqnarray*}
where $\Pi(P,Q)$ denotes the set of all joint probability
distributions whose marginals are $P$ and $Q$, see \cite[Proposition~2.63]{claus2016advancing} and \cite[Theorem~6.9]{villani2009optimal} for details.
Assumption~\ref{assumption:well-defined_rho}~(d) means that $W_r\big(P_{Y\mid X=x'},P_{Y\mid X=x}\big)\to 0$ as
$x'\to x$ in $\mathcal X$.
Note that the $\|\cdot\|^r$-weak topology can
also be generated by a Fortet–Mourier metric, which implies that $P_{Y|X=x}$ is continuous on $\mathcal{X}$ in the sense of Fortet–Mourier metric; see \cite[Proposition 2.1]{zhang2024statistical} for details.
In particular, by the definition of $\|\cdot\|^{r}$-weak topology, all conditional $q$th moments with $q\in(0,r]$, i.e., $\mathbb{E}_{P_{Y|X=x}}[\|Y\|^q]$, depend continuously on $x$. That is to say, for any $q\in(0,r]$, the function $h(y)=\|y\|^{q}$ lies in $C_{d_y}^{r}$ since $\|y\|^{q}\leq 1+\|y\|^{r}$.
Assumption \ref{assumption:well-defined_rho} (d) ensures that
$
\mathbb{E}_{P_{Y|X=x'}}[\|Y\|^q] \to \mathbb{E}_{P_{Y|X=x}}[\|Y\|^q]
$ as $x'\to x$.
We use
two examples to illustrate
how Assumption~\ref{assumption:well-defined_rho}(d) may be satisfied
for some
popular conditional distributions.

\begin{example}
(i)    Suppose that the relationship between $X$ and $Y$ is given by
\begin{eqnarray} \label{eqn:struc-YX}
    Y = \mu(X) + \sigma(X)Z,
\end{eqnarray}
where
{$\mu:\R^{d_x}\rightarrow \R^{d_y}$ and $\sigma:\R^{d_x}\rightarrow \R^{d_y\times d_y}$ are continuous mappings}, and $Z$ is a random variable with distribution $P_Z$ satisfying $\mathbb{E}_{P_Z}[Z]=0$ and $\mathbb{E}_{P_Z}[\|Z\|^r]^{1/r}<\infty$.
For $x,x'\in \mathcal{X}$, $Y_x:= (Y|X=x) = \mu(x) + \sigma(x)Z$ and $Y_{x'}:= (Y|X=x') = \mu(x') + \sigma(x')Z'$.
We choose $\hat{\pi}\in \Pi(P_{Y|X=x},P_{Y|X=x'})$ such that $Z=\tilde{Z}$ and $\tilde{Z}$ has the same probability distribution as $Z'$.
Consequently, for $x'\rightarrow x$ in $\mathcal{X}$,
\begin{eqnarray*}
        W_{r}^r(P_{Y|X=x},P_{Y|X={x'}})&=& \inf_{\pi\in \Pi(P_{Y|X=x},P_{Y|X=x'})} \left(\int \|Y_x-Y_{x'}\|^r \pi(dY_x,dY_{x'}) \right) \\
        &\leq&  \int \|Y_x-Y_{x'}\|^r \hat{\pi}(dY_x,dY_{x'}) \\
        &\leq& \mathbb{E}_{P_Z}\left[\left\| \mu(x)-\mu(x') + \big(\sigma(x)-\sigma(x')\big)Z \right\|^r\right]\\
        &=& 2^{r-1} \left(\| \mu(x)-\mu(x')\|^r + \big\|\sigma(x)-\sigma(x')\big\|^r \mathbb{E}_{P_Z}[\|Z\|^r]\right) \rightarrow 0,
    \end{eqnarray*}
which shows that Assumption \ref{assumption:well-defined_rho}(d) is satisfied in this setting.

The relationship in \eqref{eqn:struc-YX} is commonly assumed in the contextual optimization literature, see e.g., \cite{kannan2024residuals,kannan2025data}, and the continuity result can be applied to several classical joint distributions, such as the normal distribution, as described below.

    If $(Y,X)$ follows a multivariate normal distribution $\mathcal{N}(\mu,\Sigma)$ with
    \begin{eqnarray*}
        \mu=
        \begin{bmatrix}
        \mu_Y\\
        \mu_X
        \end{bmatrix}, \quad
        \Sigma = \begin{bmatrix}
            \Sigma_{YY} & \Sigma_{YX}\\
            \Sigma_{XY} & \Sigma_{XX}
        \end{bmatrix}\succ 0,
    \end{eqnarray*}
    then the conditional distribution of $Y_x:=(Y|X=x)$ is a normal distribution $\mathcal{N}(\mu(x), \sigma(x))$ with $\mu(x) := \mu_Y+ \Sigma_{YX} \Sigma_{XX}^{-1}(x - \mu_X)$ and $\sigma(x) := \Sigma_{YY} - \Sigma_{YX} \Sigma_{XX}^{-1} \Sigma_{XY}$.
    Therefore, $Y_x$ can be written by $Y_x:= \mu(x)+ \sigma(x)Z$ where $Z$ follows a $d_y$-dimensional standard multivariable normal distribution,
    and satisfies Assumption \ref{assumption:well-defined_rho}(d).

\noindent
(ii)
Wang et al. \cite{wang2022robust} and Yoon et al. \cite{yoon2025data} consider the case where $(Y,X)$ follows a Gaussian mixture model,  $\inmat{GM}\left(\left\{p^k,\mu^k,\Sigma^k\right\}_{k=1}^K\right)$ with $p\in \Delta^K$,
    where $p^k$ are mixture weights and $(\mu^k,\Sigma^k)$ are the mean and covariance of the $k$-th Gaussian component, i.e.,
    \begin{eqnarray*}
        \inmat{GM}\left(\left\{p^k,\mu^k,\Sigma^k\right\}_{k=1}^K\right) = \sum_{k=1}^K p^k \mathcal{N}(\mu^k,\Sigma^k)
    \end{eqnarray*}
    with $\mathcal{N}\left(\mu^k,\Sigma^k\right)$ denoting the $k$-th multivariate normal distribution.
    By \cite[Lemma~1]{wang2022robust}, the conditional distribution $P_{Y|X=x}$ is again a Gaussian mixture as
    $$
    GM\left(\left\{p_{Y|X=x}^k,\mu_{Y|X=x}^k,\Sigma_{Y|X=x}^k\right\}_{k=1}^K\right)= \sum_{k=1}^K p_{Y|X=x}^k \mathcal{N}(\mu_{Y|X=x}^k,\Sigma_{Y|X=x}^k)$$
    where,
    by writing the parameters in block form as
    \begin{eqnarray*}
        \mu^k :=
\begin{bmatrix}
\mu_Y^k \\
\mu_X^k
\end{bmatrix}
\in \mathbb{R}^{d_X+d_Y},\quad
\Sigma^k :=
\begin{bmatrix}
\Sigma_{YY}^k & \Sigma_{YX}^k \\
\Sigma_{XY}^k & \Sigma_{XX}^k
\end{bmatrix}
\in \mathbb{S}_+^{d_X+d_Y},
    \end{eqnarray*}
the conditional mean, covariance and mixture weights are
    \begin{eqnarray*}
\mu_{Y|X=x}^k &=& \mu_Y^k + \Sigma_{YX}^k (\Sigma_{XX}^k)^{-1} (x - \mu_X^k), \\
\Sigma_{Y|X=x}^k &=& \Sigma_{YY}^k - \Sigma_{YX}^k (\Sigma_{XX}^k)^{-1} \Sigma_{XY}^k, \\
p_{Y|X=x}^k &=& \frac{p^k \mathcal{N}(x | \mu_X^k, \Sigma_{XX}^k)}{\sum_{j=1}^K p^j \mathcal{N}(x | \mu_X^j, \Sigma_{XX}^j)}.
    \end{eqnarray*}
Note that all parameters are continuous in $x$.
To facilitate the notation, we write $p_{x}^k:=p_{Y|X=x}^k$ and $\mathcal{N}_{x}^k:= \mathcal{N}(\mu_{Y|X=x}^k,\Sigma_{Y|X=x}^k)$.
For any $x,x'\in \mathcal{X}$,
\begin{eqnarray*}
    &&W_r(P_{Y|X=x},P_{Y|X={x'}})  = W_r \left(\sum_{k=1}^K p_{x}^k \mathcal{N}_{x}^k, \sum_{k=1}^K p_{x'}^k \mathcal{N}_{x'}^k \right)\\
    && \leq W_r \left(\sum_{k=1}^K p_{x}^k\mathcal{N}_{x}^k, \sum_{k=1}^K p_{x'}^k \mathcal{N}_{x}^k \right)
    +  W_r \left(\sum_{k=1}^K p_{x'}^k \mathcal{N}_{x}^k, \sum_{k=1}^K p_{x'}^k \mathcal{N}_{x'}^k \right)\\
    &&\leq \left(\frac{1}{2}\sum_{k=1}^K \left| p_{x}^k - p_{x'}^k \right|\right)^{1/r} \max_{1\leq i,j \leq K} W_r \left( \mathcal{N}_{x}^i, \mathcal{N}_{x}^j \right)
    +
    \sum_{k=1}^K p_{x'}^k W_r \left(\mathcal{N}_{x}^k, \mathcal{N}_{x'}^k \right),
\end{eqnarray*}
where the first inequality comes from the triangle inequality of the Wasserstein distance, and the second inequality comes from Lemma \ref{lemma:bound_wass-dist} and the convexity of the Wasserstein distance; see, e.g., \cite[Theorem 4.8]{villani2009optimal}.
By the continuity of $p_{x}^k$ in $x$ and the continuity of the conditional multivariate normal distribution $\mathcal{N}_x^k$ shown in part~(a), we have $W_r(P_{Y|X=x},P_{Y|X={x'}})\rightarrow 0$ as $x'\rightarrow x$, which shows that the Gaussian mixture model satisfies Assumption \ref{assumption:well-defined_rho}(d).

\noindent
(iii) Ban and Rudin \cite{ban2019big} use the kernel regression to approximate the conditional distribution $P_{Y|X}$:
Given i.i.d. samples $\{(x_i,y_i)\}_{i=1}^N$ and a kernel $k:\mathcal{X}\times \mathcal{X}\rightarrow \R_+$,
the conditional distribution of $Y$ given $X=x$ is predicted as
    \begin{eqnarray*}
        \hat{P}_{Y|X=x} = \sum_{i=1}^N\frac{ k(x,x_i) }{\sum_{j=1}^N k(x,x_j)}\delta_{y_i},
    \end{eqnarray*}
where $\delta_{y_i}$ is a Dirac measure concentrated at $y_i$.
Then for $x,x'\in \mathcal{X}$,
\begin{eqnarray*}
    &&W_r\left(\hat{P}_{Y|X=x}, \hat{P}_{Y|X=x'}\right)\\
    && =  W_r\left( \sum_{i=1}^N \frac{ k(x,x_i) }{\sum_{j=1}^N k(x,x_j)} \delta_{y_i}, \sum_{i=1}^N\frac{ k(x',x_i) }{\sum_{j=1}^N k(x',x_j)}\delta_{y_i} \right)\\
    && \leq  \left(\frac{1}{2}\sum_{k=1}^K \left| \frac{ k(x,x_k) }{\sum_{j=1}^N k(x,x_j)} - \frac{ k(x',x_k) }{\sum_{j=1}^N k(x',x_j)} \right|\right)^{1/r} \max_{1\leq i,j \leq K} W_r \left( \delta_i, \delta_j \right)\\
    && \leq  \left(\frac{1}{2}\sum_{k=1}^K \left| \frac{ k(x,x_k) }{\sum_{j=1}^N k(x,x_j)} - \frac{ k(x',x_k) }{\sum_{j=1}^N k(x',x_j)} \right|\right)^{1/r} \max_{1\leq i,j \leq K} \| y_i, y_j \|,
\end{eqnarray*}
where the first inequality comes from Lemma \ref{lemma:bound_wass-dist}.
If the kernel is selected to be nonnegative and continuous, then
the function $\frac{k(x,x_k)}{\sum_{j=1}^N k(x, x_j)}$ is also continuous in $x$. This continuity implies
that $W_r(\hat{P}_{Y|X=x},\hat{P}_{Y|X=x'}) \to 0$ as $x' \to x$, and therefore, the conditional distribution obtained via the kernel regression satisfies Assumption \ref{assumption:well-defined_rho}(d).

\end{example}

}

The next proposition addresses the well-definedness of $\rho_{P_X}^{(1)}\circ \rho_{P_{Y|X}}^{(2)}$ in (\ref{eqn:ex-post-interchange}).
It demonstrates that the inner part $\rho^{(2)}_{P(Y|X=\cdot)} ( c(g(\cdot),Y)): \mathcal{X}\rightarrow \R$
is continuous in $x$ over $\mathcal{X}$ and thus is an appropriate random function of $X$.

\begin{proposition} \label{prop:random_variable_inner_risk}
    Consider problem (\ref{eqn:ex-post-interchange}).
    Under
    Assumption \ref{assumption:well-defined_rho},
        $\rho^{(2)}_{P_{Y|X=x}} ( c(g(x),Y)): \mathcal{X}\rightarrow \R$         is continuous in  $x$ over $\mathcal{X}$.
\end{proposition}

\noindent
\textbf{Proof.}
See Section \ref{proof:Prop1}.
\hfill$\Box$

Next, we need to ensure the problem is well-defined over the hypothesis set $\widetilde{\mathcal{G}}$.

\begin{assumption} \label{assumption:well-defined_obj}
    Consider problems
    (\ref{eqn:ex-post-interchange}) and (\ref{eqn:ex_ante}).
        (a) The support sets $\mathcal{X},\mathcal{Y}$
        of random variables $X$ and $Y$
        are bounded.
        (b) The cost function $c(z,y)$ is H\"older continuous in $z\in\mathcal{Z}$, i.e., there exists a continuous function $r:\mathcal{Y}\rightarrow \R_+$ with
        supremum norm $\|r\|_{\infty,\mathcal{Y}}:=\sup\{|r(y)|:y\in\mathcal{Y}\}$
        being
        finite and a positive constant $\nu>0$ such that
        $
        |c(z_1,y) - c(z_2,y)| \leq r(y)\|z_1-z_2\|^\nu, \forall z_1,z_2 \in\mathcal{Z}
        $
        uniformly for all $y\in\mathcal{Y}$.
\end{assumption}

With the continuity of the cost function $c$ and the policy $g$, condition (a) ensures that $c(g(X),Y)$ is a bounded random function of $(X,Y)$, which will facilitate the forthcoming discussions.
        In practice, the hypothesis
        space $\widetilde{\mathcal{G}}$ is usually restricted
        to a
        class of continuous functions with finite parameters.
        In that case, the boundedness of $\mathcal{X},\mathcal{Y}$ can be relaxed. We will discuss this case in Section \ref{sec:DR-rkhs}.
The H\"older continuity of the cost function in condition (b) is widely
used in the literature and can be verified in many practical applications.
        Assumption (b) also ensures the growth condition in Assumption \ref{assumption:well-defined_rho} (c).
{The next proposition establishes
the continuity of objective functions in problems (\ref{eqn:ex_ante}) and (\ref{eqn:ex-post-interchange}) in $g$.}

\begin{proposition} \label{prop:well-definedness-problem}
   Under Assumptions \ref{assumption:well-defined_rho} and \ref{assumption:well-defined_obj}, the objective function $\rho_P^{\text{EA}}\big(c(g(X),Y)\big)$ and $\rho_{P_X}^{(1)}\left(\rho_{P_{Y|X}}^{(2)}\big(c(g(X),Y)\big)\right)$ are continuous in $g\in\widetilde{\mathcal{G}}$.
\end{proposition}

\noindent
\textbf{Proof.}
See Section \ref{proof:Prop2}.
\hfill$\Box$

\subsection{Contextual consistency}

{
A natural requirement on the choice of the ex ante risk is that, for a conditional risk measure $\rho_{P_{Y|X}}^{(2)}$ and
any two policies $g_1,g_2\in\widetilde{\mathcal{G}}$,
if $g_1$ has better performance than $g_2$ in all contexts $x\in \mathcal{X}$ in the sense of $\rho_{P_{Y|X=x}}^{(2)}(c(g(x),Y))$ for $x\in \mathcal{X}$ almost surely, then $g_1$ should incur a lower ex ante risk than $g_2$.}
A property called contextual consistency
is introduced to fulfill this requirement.

\begin{definition} [Contextual consistency] \label{def:contextual-consist}
    An ex ante risk  $\rho_{P}^{\text{EA}}$ is \textit{contextually consistent}
    with
    risk measure $\rho_{P_{Y|X}}^{(2)}$
    if for any joint probability distribution of $(X,Y)$, any continuous cost function $c(z,y)$ and any continuous policies $g_1,g_2: \mathcal{X} \rightarrow \mathcal{Z}$ that satisfy
    \bgeqn
    \label{eq:context-consist-1}
    \rho_{P_{Y|X=x}}^{(2)}(c(g_1(x),Y)) \leq \rho_{P_{Y|X=x}}^{(2)}(c(g_2(x),Y)), \forall x\in \mathcal{X}\; \inmat{almost surely},
    \edeqn
    the following inequality holds:
    \bgeqn
    \label{eq:context-consist-2}
        \rho_{P}^{\text{EA}}(c(g_1(X),Y)) \leq \rho_{P}^{\text{EA}}(c(g_2(X),Y)).
    \edeqn
    An ex ante risk $\rho_{P}^{\text{EA}}$ is
    said to be
    \textit{strictly contextually consistent} if additionally
   strict inequality in (\ref{eq:context-consist-2}) holds
   when
   strict inequality in (\ref{eq:context-consist-1})
    holds with strictly positive probability.

\end{definition}

{
This definition does not require $\rho_P^{\text{EA}}=\rho_{P_X}^{(1)}\circ \rho_{P_{Y|X}}^{(2)}$.
Nevertheless, if
an ex ante risk can be decomposed
into
a nested form, i.e., there exists a monotone risk measure $\rho_{P_X}^{(1)}$
such that $\rho_{P}^{\text{EA}} = \rho_{P_X}^{(1)}\circ \rho_{P_{Y|X}}^{(2)}$, then the ex ante risk is contextually consistent.
If, in addition, $\rho_{P_X}^{(1)}$ is strictly monotone,
then $\rho_{P}^{\text{EA}}$ is \textit{strictly contextually consistent}.
This observation allows us to
connect
the decomposability of ex ante risk
with
the contextual consistency of ex post risk.}
This
property
coincides with
time consistency in dynamic risk measure
as stated in Section 11.2 in \cite{follmer2011stochastic}.
However, the purpose
of discussing the contextual consistency here
differs from
that of time consistency in
dynamic risk management.
The latter is used to derive
dynamic Bellman equations for
finding a
holistic policy
in a multistage
risk minimization problem
whereas
the former is on the opposite, that is,
we solve
(\ref{eqn:ex_ante}) to obtain
an optimal solution of
(\ref{eqn:ex-post-interchange}) under the
contextual consistency.
Moreover, the dynamic consistency emphasizes
{consistency} of
the  decision-making at each stage,
while the contextual consistency
is more focused on the relationship
between the conditional risk arising from PDU and the overall integrated risk arising from CU.

\section{Risk-averse contextual optimization and examples}
\label{sec:contextual-ra-opt}

In this section, we
establish the equivalence between the conditional risk minimization problem (\ref{eqn:conditional_risk_min}) and {risk-averse contextual optimization problems (\ref{eqn:ex-post-interchange}) and (\ref{eqn:ex_ante})}.
We examine the application of some classical risk measures and their properties in the setting of contextual optimization.

\subsection{Optimality of risk-averse contextual optimization} \label{sec:model-crao}

We now move on to
discuss properties of
problems (\ref{eqn:ex-post-interchange}) and (\ref{eqn:ex_ante}),
and demonstrate the conditional optimality of the optimal policy derived by risk-averse contextual optimization problems.
To facilitate the discussion, we make
an assumption which
guarantees well-definedness
of the problem.

\begin{assumption} \label{assumption:convex_finite-minimum}
{
(a) For any $\hat{x}\in \mathcal{X}$, there exists a constant $\kappa \in \R$ and a compact set $\hat{\mathcal{Z}}\subset \mathcal{Z}$ such that for every $x$ in a neighborhood of $\hat{x}$, the lower level set $L_x(\kappa):=\left\{ z\in \mathcal{Z}: \rho_{ P_{Y|X=x}}^{(2)}(c(z,Y))\leq \kappa \right\}$ is nonempty and contained in $\hat{\mathcal{Z}}$.
}
(b) The optimal solution to $\min\limits_{z\in \mathcal{Z}} \rho_{ P_{Y|X=x}}^{(2)}(c(z,Y))$ is unique for $x\in\mathcal{X}$.

\end{assumption}

{
Assumption \ref{assumption:convex_finite-minimum} (a)
is known as the inf-compactness condition in the literature of optimization, see,
e.g., \cite[Section 4]{bonnans2013perturbation}.
The uniform inf-compactness assumption holds automatically when the feasible set $\mathcal{Z}$ is compact.
}
Assumption \ref{assumption:convex_finite-minimum} (b) is satisfied in some important  applications.
For example, when $\rho_{ P_{Y|X=x}}^{(2)}$ is chosen as an entropic risk measure and $c(z,y)$ is convex in $z$,
$\rho_{ P_{Y|X=x}}^{(2)}(c(z,Y))$
is strictly convex in $z$ as proved in \cite[Proposition 2.6]{cherny2007divergence} and thus the optimal solution to $\min\limits_{z\in \mathcal{Z}} \rho_{ P_{Y|X=x}}^{(2)}(c(z,Y))$ is unique.
When {$\rho_{ P_{Y|X=x}}^{(2)}=\inmat{CVaR}_{\beta,P_{Y|X=x}}$} is applied to the newsvendor problem with
$
    c(z,y) = h(z-y)_+ + b(y-z)_+
$
where
$h,b>0$ and $(a)_+ := \max\{a,0\}$,
the optimal solution can be represented
via the quantile function $F^{-1}$ of $P_{Y|X=x}$, i.e.,
$
    z^{*} = \frac{h}{h+b}F^{-1}\left(\frac{b(1-\beta)}{h+b}\right) + \frac{b}{b+h}F^{-1}\left(\frac{h\beta+b}{b+h}\right),
$
see
\cite{gotoh2007newsvendor}
for details.
In this case, $z^*$ is unique provided that $P_{Y|X=x}$ has a continuous distribution.
{
We verify these assumptions in the newsvendor problem and the portfolio selection problem in Section \ref{sec:vei-assumption-application}.
}

The next theorem states that, under some additional moderate
conditions, the policies $g^*$ obtained
from solving problem (\ref{eqn:ex-post-interchange})
and
(\ref{eqn:ex_ante})
are
optimal in conditional risk minimization problem (\ref{eqn:conditional_risk_min}) under $\rho_{P_{Y|X}}^{(2)}$ almost surely, i.e.,
$
g^*(x) \in \arg\inf\limits_{z\in \mathcal{Z}} \rho_{P_{Y|X=x}}^{(2)} \big(c(z,Y)\big), \forall x\in \mathcal{X}
$
almost surely.
Therefore, by employing the strictly contextually consistent ex ante risk mapping,
the risk-averse contextual optimization problem
is not only
computationally tractable but also exhibits good theoretical properties.

\begin{theorem} \label{thm: optimality-of-ent}
    Under Assumptions~\ref{assumption:well-defined_rho} and
    \ref{assumption:convex_finite-minimum},
    the following assertions hold.
    \begin{itemize}
        \item[(i)] The         optimal solution mapping
        $z
        ^*(x) := \arg\min\limits_{z\in \mathcal{Z}} \rho_{ P_{Y|X=x}}^{(2)}(c(z,Y))$ is continuous in~$x$.

        \item[(ii)]
    For problem (\ref{eqn:ex_ante}),
    {if $\rho_{P}^{\text{EA}}$ is
    contextually consistent
    with $\rho_{P_{Y|X}}^{(2)}$,
    then the optimal solution to problem (\ref{eqn:ex_ante}) exists;}
    moreover, if $\rho_{P}^{\text{EA}}$ is
    strictly contextually consistent     with $\rho_{P_{Y|X}}^{(2)}$,  then the optimal policy $g^*$ obtained
    from solving problem (\ref{eqn:ex_ante}) is also  optimal in problem (\ref{eqn:conditional_risk_min}) in all contexts almost surely.

        \item[(iii)]
        For problem (\ref{eqn:ex-post-interchange}),
        {the optimal policy $g^*$ to problem (\ref{eqn:ex-post-interchange}) exists;}
        moreover
        if $\rho_{P_X}^{(1)}$ is strictly monotone,
        then
        $g^*$
        is also an optimal policy
        to
        problem (\ref{eqn:conditional_risk_min})         almost surely.
    \end{itemize}
\end{theorem}

\noindent
\textbf{Proof.}
See Section \ref{proof:Thm1}.
\hfill$\Box$

Theorem \ref{thm: optimality-of-ent} (iii) is analogous to Proposition 3.1 in \cite{shapiro2017interchangeability}
(without uniqueness).
Here we impose uniqueness
to ensure the
continuity of the optimal policy.
Given the particular interest of this theorem in contextual optimization, we provide the proof in Section \ref{proof:Thm1}.
 The result implies that the optimal policy is independent of the choice of $\rho_{P_X}^{(1)}$ if it is strictly monotone,
{
which is consistent with the intuition behind the modeling of the risk-averse contextual optimization, as explained below.
The conditional risk minimization problem \eqref{eqn:conditional_risk_min} minimizes $\rho_{P_{Y|X=x}}^{(2)}(c(g(x),Y))$ for each context $X=x$, and hence the resulting conditionally optimal policy does not involve the contextual risk measure $\rho_{P_X}^{(1)}$.
Problem \eqref{eqn:ex-post-interchange} is introduced as an equivalent counterpart of problem \eqref{eqn:conditional_risk_min} under the assumptions of Theorem \ref{thm: optimality-of-ent}, and
thus its optimal policy should also be independent of the choice of $\rho_{P_X}^{(1)}$.

From a mathematical perspective, this property follows from the strict monotonicity of $\rho_{P_X}^{(1)}$.
For any fixed policy $g$, the conditional risk satisfies
\begin{eqnarray*}
\rho_{P_{Y|X=x}}^{(2)}(c(g(x),Y))\geq \psi(x),\; \forall x\in \mathcal{X} \text{ almost surely},
\end{eqnarray*}
and the monotonicity of $\rho_{P_X}^{(1)}$ yields
$$\rho_{P_{X}}^{(1)}\left( \rho_{P_{Y|X=x}}^{(2)}(c(g(x),Y))\right)\geq \rho_{P_{X}}^{(1)}(\psi(x)).$$
Thus, any conditionally optimal policy is also optimal for \eqref{eqn:ex-post-interchange}, regardless of $\rho_{P_{X}}^{(1)}$.
Nevertheless, the monotonicity of $\rho_{P_X}^{(1)}$ alone does not ensure that an optimal policy for problem \eqref{eqn:ex-post-interchange} is conditionally optimal for problem  \eqref{eqn:conditional_risk_min}, as illustrated in Example~\ref{example:nested-cvar}.
For this issue, we assume strict monotonicity of $\rho_{P_{X}}^{(1)}$ to ensure the conditional optimality of the policy for problem \eqref{eqn:ex-post-interchange}.
We will illustrate this observation in Example~\ref{example:expected-ent}.

Note that if
$\rho_P^{\text{EA}}\big(c(g(X),Y)\big)
=\rho_{P_X}^{(1)}\left(\rho_{P_{Y|X}}^{(2)}\big(c(g(X),Y)\big)\right)$,
then the strict monotonicity of $\rho_{P_X}^{(1)}$
is a sufficient
condition for the strict contextual consistency of $\rho_P^{\text{EA}}$.
}
The next two examples show that this theorem cannot be extended to the general risk measures/mappings including ex ante CVaR, i.e.,
$\inmat{CVaR}_{\beta, P}(c(g(X),Y))$,
and
nested CVaR, i.e.,
$\inmat{CVaR}_{\beta, P_X}\left(\inmat{CVaR}_{\beta, P_{Y|X}}(c(g(X),Y))\right)$.

\begin{example} [CVaR]
\label{example:ex-ante-cvar}
The policy derived from ex ante CVaR may not be optimal for the conditional risk minimization problem in every context.
To see this,
{
consider a risk-averse contextual optimization problem \eqref{eqn:ex_ante} with $\rho_P^{\text{EA}}= \text{CVaR}_{\beta,P}$,
}
\begin{eqnarray} \label{eqn:cvar-ex-ante}
    \min_{g\in\widetilde{\mathcal{G}}} \inmat{CVaR}_{\beta,P}(c(g(X),Y)) =  \min_{t\in\R,g\in\widetilde{\mathcal{G}}} \mathbb{E}_{P}\left[t+\frac{1}{1-\beta}\left[c(g(X),Y) - t\right]_+\right].
\end{eqnarray}
Likewise, for $x\in \mathcal{X}$ almost surely,
consider { the conditional risk minimization problem \eqref{eqn:conditional_risk_min} with $\rho_{P_{Y|X}}^{(2)}(\cdot)=\inmat{CVaR}_{\beta,P_{Y|X}}(\cdot)$, }
\begin{eqnarray} \label{eqn:cvar-each-contexts}
    &&\min_{z\in\R^{d_z}} \inmat{CVaR}_{\beta,P_{Y|X=x}}(c(z,Y))\nonumber\\
    &=&  \min_{t\in\R,z\in\R^{d_z}} \mathbb{E}_{P_{Y|X=x}}\left[t+\frac{1}{1-\beta}\left[c(z,Y) - t\right]_+\right].
\end{eqnarray}
Since
$t$ in problem (\ref{eqn:cvar-ex-ante}) is
independent of the realization of covariate $X$,
whereas $t$ in conditional risk minimization problem (\ref{eqn:cvar-each-contexts}) depends on the realization of $X$,
there always exists some subset of  $\mathcal{X}$ with non-zero measure such that
the optimal solution
$(g^*,t^*)$ to
problem (\ref{eqn:cvar-ex-ante})
is not an optimal solution of
problem (\ref{eqn:cvar-each-contexts})
provided that
$t^*$ in (\ref{eqn:cvar-each-contexts})
is not constant for
almost every
$x\in \mathcal{X}$.
{
This example illustrates that
if $\rho_P^{\text{EA}}$ is not contextually consistent for the inner risk measure $\rho_{P_{Y|X}}^{(2)}$, then Theorem~\ref{thm: optimality-of-ent}(ii) does not apply here and a policy obtained from problem \eqref{eqn:ex_ante} may fail to be conditionally optimal.
}
\end{example}

\begin{example} [Nested-CVaR]
\label{example:nested-cvar}
The policy derived from nested-CVaR may not be conditionally optimal as well.
Consider
{
a risk-averse contextual optimization problem \eqref{eqn:ex-post-interchange} with nested-CVaR by setting $\rho_{P_X}^{(1)}(\cdot)=\inmat{CVaR}_{\beta,P_X}(\cdot)$ and $\rho_{P_{Y|X}}^{(2)}(\cdot)=\inmat{CVaR}_{\beta,P_{Y|X}}(\cdot)$,}
\begin{eqnarray} \label{eqn:nested-cvar-obj}
    \min_{g\in\widetilde{\mathcal{G}}} \inmat{CVaR}_{\beta,P_X}\left(\inmat{CVaR}_{\beta,P_{Y|X}}(c(g(X),Y))\right),
\end{eqnarray}
which is not equivalent to ex ante CVaR model (\ref{eqn:cvar-ex-ante}), see e.g., \cite{iancu2015tight}.
By Theorem~\ref{thm: optimality-of-ent}(i), $z^*(x)$
obtained from solving problem (\ref{eqn:cvar-each-contexts}) for $x\in \mathcal{X}$ almost surely
is an optimal solution to problem
(\ref{eqn:nested-cvar-obj}).
However, since CVaR is not strictly monotone
\cite[Proposition 6.37]{shapiro2014lectures})
and only focuses on the tail behavior, it is possible for an optimal policy $g^*$ of problem (\ref{eqn:nested-cvar-obj}) to satisfy the inequality
$$\inmat{CVaR}_{\beta,P_{Y|X=x}} \big(c(g^*(x),Y)\big) > \min_{z\in\mathcal{Z}}\inmat{CVaR}_{\beta,P_{Y|X=x}} \big(c(z,Y)\big)$$
on a subset of $\mathcal{X}$ with non-zero measure.
{
This example shows that if $\rho_{P_X}^{(1)}(\cdot)$ is not strictly monotone, then the assumption of Theorem \ref{thm: optimality-of-ent}(iii) is violated and a policy obtained from problem \eqref{eqn:ex-post-interchange} may still fail to be conditionally optimal.
Hence, one cannot conclude that every optimizer of \eqref{eqn:ex-post-interchange} is conditionally optimal.
}
\end{example}

\subsection{Examples of ex ante risk mappings}
\label{sec:example-oce-ent}

To further illustrate the relationship between problem (\ref{eqn:ex-post-interchange}) and problem (\ref{eqn:ex_ante})
and show how to construct contextually consistent ex ante risk mappings using Theorem~\ref{thm: optimality-of-ent},
we provide the following examples.
In the first example, we deal with the case where DM's risk attitude is characterized as
OCE (including CVaR and entropic risk measure)
in conditional risk minimization, {and
in the second example, we consider a more general class of risk measures, represented as the worst expected value over a family of parametric functions, studied in \cite{guigues2023risk} and \cite{huang2020stochastic}.}

\begin{example}[Expected optimized certainty equivalent]
\label{example:expected-oce}
    Consider the case that
    $ \rho_{P_X}^{(1)} = \mathbb{E}_{P_X}$ and $\rho_{P_{Y|X}}^{(2)}$ is OCE
  defined
  as in (\ref{eq:neg-OCE}).
  Then Theorem \ref{thm: optimality-of-ent}(iii) is applicable to ensure the conditional optimality in problem (\ref{eqn:conditional_risk_min}).
   Moreover, observe that model (\ref{eqn:ex-post-interchange}) is equivalent to
    \begin{eqnarray} \label{eqn:expected-oce}
        &&\min_{g\in\widetilde{\mathcal{G}}} \mathbb{E}_{P_X}\left[\rho_{P_{Y|X}}^{(2)}(c(g(X),Y))\right]
        = \min_{g\in\widetilde{\mathcal{G}}}\mathbb{E}_{P_X} \left[ - \max_{s \in \R} \left\{\mathbb{E}_{P_{Y|X}}\big[s + u(-c(g(X),Y)-s)\big] \right\}  \right]\nonumber\\
        && \qquad = \min_{g\in\widetilde{\mathcal{G}},s\in \widetilde{\mathcal{S}}} \bigg\{ - \mathbb{E}_{P} \big[s(X) + u(-c(g(X),Y)-s(X))\big]\bigg\},
    \end{eqnarray}
    where
     $\widetilde{\mathcal{S}}$ denotes the set of all measurable functions from $\mathcal{X}\rightarrow \R$ and
    the second
    equality is due to
    the interchangeability principle \cite[Theorem 2.2]{hiai1977integrals}.
    In
    this
    example,
    the objective function is written in the form of $\rho_{P}^{\text{EA}}$ and takes the expectation over the joint distribution $P$, which
    will greatly facilitate the data-driven optimization.
    Specifically, when $u(x) = \frac{1}{1-\beta}x$ for  $x\leq 0$ and $u = 0$ for $x>0$, we obtain the Expected CVaR mapping
    \begin{eqnarray}
    \label{eq:CVaR-opt-X}
        \min_{g\in\widetilde{\mathcal{G}},t\in \widetilde{\mathcal{T}}} \mathbb{E}_{P} \left[t(X) + \frac{1}{1-\beta}(c(g(X),Y)-t(X))_+\right],
    \end{eqnarray}
    where $\widetilde{\mathcal{T}}$
    is
    the set of all
    measurable functions $t:\mathcal{X}\rightarrow\R$.
    In this way, we can reformulate model (\ref{eqn:ex-post-interchange}) into
    an augmented version of
    model (\ref{eqn:ex_ante}).

\end{example}

Example \ref{example:expected-oce} shows that
the decision maker could consider expected CVaR \eqref{eq:CVaR-opt-X} whose properties are better than ex ante CVaR (Example \ref{example:ex-ante-cvar}) and Nested-CVaR (Example \ref{example:nested-cvar}) for a risk-averse contextual optimization problem { in terms of conditional optimality
ensured by strict monotonicity of $ \rho_{P_X}^{(1)} = \mathbb{E}_{P_X}$ and
computational tractability, since it only involves an expectation over the joint distribution as in \eqref{eq:CVaR-opt-X}.
This kind of reformulation can also be obtained
for a broader class of parametric risk measures using Legendre–Fenchel conjugate.
}

{
\begin{example} [Nested parametric risk measure]
    Consider the case that $\rho_{P_X}^{(1)}$ and $\rho_{P_{Y|X}}^{(2)}$ are both parametric risk measures defined in \eqref{eqn:parametric-rm}, i.e., $\rho_{P_X}^{(1)}(\cdot)=\inf_{s_1\in \mathcal{S}_1} \mathbb{E}[U_1(\cdot,s_1)]$ and $\rho_{P_{Y|X}}^{(2)}(\cdot)=\inf_{s_2\in \mathcal{S}_2} \mathbb{E}[U_2(\cdot,s_2)]$, where $U_1:\R\times \mathcal{S}_1\rightarrow \R$ and $U_2:\R\times \mathcal{S}_2\rightarrow \R$ are proper, continuous functions that are jointly convex on $\R\times \mathcal{S}_1$ and $\R\times \mathcal{S}_2$, respectively, and are non-decreasing in their first argument.
    Then the problem at the lhs of equality \eqref{eqn:ex-post-interchange} can be written as
    \begin{eqnarray} \label{eqn:nest-parametric-rm1}
       &&\min_{s_1\in \mathcal{S}_1}  \mathbb{E}_{P_{X}}\left[ U_1\left( \min_{z\in \mathcal{Z}}\min_{s_2 \in \mathcal{S}_2} \mathbb{E}_{P_{Y|X}}\Big[ U_2\big(c(z,Y),s_2\big) \Big],s_1  \right) \right] \nonumber\\
       &=& \min_{\substack{g(\cdot)\in \widetilde{\mathcal{G}} \\
       s_1\in \mathcal{S}_1,
       s_2(\cdot)\in \widetilde{\mathcal{S}}_2}
       }  \mathbb{E}_{P_{X}}\left[ U_1\left(  \mathbb{E}_{P_{Y|X}}\Big[ U_2\big(c(g(X),Y),s_2(X)\big) \Big],s_1  \right) \right],
    \end{eqnarray}
    where
    $\widetilde{\mathcal{S}}_2$ is the set of all
    measurable functions $s_2:\mathcal{X}\rightarrow \mathcal{S}_2$,
    and the equality follows from the monotonicity of $U_1(w,s_1)$ in $w$ and the interchangeability principle.
    Since $U_1(w,s_1)$ is in general nonlinear in $w$,  problem \eqref{eqn:nest-parametric-rm1} is nonlinear in $\mathbb{E}_{P_{Y|X}}$.
However, we can follow Dai et al.~\cite{dai2017learning} to represent $U_1(\cdot,s_1)$ in terms of its
 Legendre–Fenchel conjugate and consequently obtain a linear representation of the term.
Recall that
the Legendre–Fenchel conjugate of $U_1(\cdot,s_1)$ is a function $U_1^*:\R\times \mathcal{S}_1 \rightarrow \R$ defined by $U_1^*(v,s_1):= \sup_{w\in \R} \left\{ wv - U_1(w,s_1) \right\}$.
Since $U_1(w,s_1)$ is proper, continuous, and convex in $w$, $U_1^*(v,s_1)$ is also proper, convex, and lower semicontinuous in $v$, and
$
    U_1(w,s_1) := \sup_{v\in \R} \left\{ wv - (U_1^*(v,s_1)) \right\}.
$
Substituting this expression into the above problem and invoking interchangeability again,
we arrive at
\begin{eqnarray}\label{eqn:nested-parametric-rm}
       &&\min_{g(\cdot), s_1, s_2(\cdot)}    \mathbb{E}_{P_{X}}\left[
       \sup_{v\in \R} \left\{ \mathbb{E}_{P_{Y|X}}\Big[ U_2\big(c(g(X),Y),s_2(X)\big) \Big] v-U_1^*(v,s_1) \right\}\right]\nonumber\\
       &=&\min_{g(\cdot), s_1, s_2(\cdot)} \sup_{v(\cdot)\in \widetilde{\mathcal{V}}}   \mathbb{E}_{P}\left[
          U_2\big(c(g(X),Y),s_2(X)\big) v(X)-U_1^*(v(X),s_1) \right],
\end{eqnarray}
where $\widetilde{\mathcal{V}}$ is the set of all measurable functions $v:\mathcal{X}\rightarrow\R$.

The objective function
on the rhs of \eqref{eqn:nested-parametric-rm}
in linear in $\bbe_P$.
Since $U_1(w,s_1)$ is jointly convex on $\R\times \mathcal{S}_1$,
then
$U_1^*$ is concave in $s_1$.
Moreover, since $U_1(w,s_1)$ is non-decreasing in $w$, we have $U_1^*(v,s_1)=+\infty$ for any $v< 0$ and $s_1\in\mathcal{S}_1$, which implies that the supremum in $v(\cdot)$ in problem \eqref{eqn:nested-parametric-rm} must be attained at some $v(\cdot)$ such that $v(x)\geq 0$ for all $x\in \mathcal{X}$.
Furthermore, since $U_2(w,s_2)$ is jointly convex on $\R\times \mathcal{S}_2$,
 the above arguments imply that the objective function in \eqref{eqn:nested-parametric-rm} is convex in $(s_1,s_2(\cdot))$ and concave in $v(\cdot)$.
Consequently,
we may  interpret
$$
\rho(c(g(X),Y)):=\min_{ s_1\in\mathcal{S}_1, s_2(\cdot)\in \widetilde{\mathcal{S}}_2} \sup_{v(\cdot)\in \widetilde{\mathcal{V}}}   \mathbb{E}_{P}\left[
          U_2\big(c(g(X),Y),s_2(X)\big) v(X)-U_1^*(v(X),s_1) \right]$$
as a
saddle-point risk measure discussed in \cite{huang2020stochastic}
despite they focus on parameters in finite-dimensional spaces.
The minimax optimization problem on the rhs of \eqref{eqn:nested-parametric-rm}
may be solved by
existing algorithms
for stochastic saddle-point optimization (see, e.g., \cite{du1995minimax,nemirovski2009robust,nemirovski2002efficient}) after some appropriate pre-treatment such as  discretization and use of decision rules to convert the infinite-dimensional problem into an approximate finite-dimensional problem.

Consider the special case when $\rho_{P_X}^{(1)}$ and $\rho_{P_{Y|X}}^{(2)}$ are both OCE defined in \eqref{eq:neg-OCE}, i.e.,
$U_1(w,s_1)= - s_1-u_1(-w-s_1)$, $U_2(w,s_2)= - s_2-u_2(-w-s_2)$, and $\mathcal{S}_1=\mathcal{S}_2=\R$,
where $u_1,u_2:\R \rightarrow\R$ are proper, closed, concave, and non-decreasing utility functions.
Since $-u_1$ is proper, convex, and continuous,
we have
\begin{eqnarray*}
    U_1^*(v,s_1) &=& \sup_{w\in \R} wv-(- s_1-u_1(-w-s_1)) = \sup_{w'\in \R} (-w'-s_1)v-(- s_1-u_1(w'))\\
    &=& s_1-s_1v+ \sup_{w'\in \R} w'(-v)-(-u_1(w'))  = s_1-s_1v+ (-u_1)^*(-v),
\end{eqnarray*}
and thus \eqref{eqn:nested-parametric-rm} can be written as follows:
\begin{eqnarray}\label{eqn:nested-oce}
    &&
    \min_{g,s_1,s_2(\cdot)} \sup_{v(\cdot)}  \mathbb{E}_{P}\Bigg[\Big(s_2(X) +u_2\big(-c(g(X),Y)- s_2(X)\big) \Big)v(X)
    \nonumber\\
    &&\qquad \qquad\qquad\quad
    - \left(s_1+s_1v(X)+(-u_1)^*(v(X))\right)\Bigg]
\end{eqnarray}
In particular, we consider the two special cases.
\begin{itemize}
    \item[(i)]
$\rho_{P_X}^{(1)}= \mathbb{E}_{P_X}$.
The corresponding utility function is the identical function, i.e., $u_1(w)=w$, and thus $(-u_1)^*(v)=0$ if $v=-1$; and $(-u_1)^*(v)=+\infty$ if $v\neq -1$.
In this case, the optimal $v(\cdot)$ in problem \eqref{eqn:nested-oce} is $v(x)=-1$ for all $x\in \mathcal{X}$, and problem \eqref{eqn:nested-oce} reduces to
\begin{eqnarray*}
    \min_{g,s_2(\cdot)}  \mathbb{E}_{P}\Bigg[ -\Big( s_2(X) +u_2\big(-c(g(X),Y)- s_2(X)\big)    \Big)  \Bigg],
\end{eqnarray*}
which is equivalent to problem \eqref{eqn:expected-oce} in Example \ref{example:expected-oce}.

\item[(ii)] $\rho_{P_X}^{(1)}= \text{CVaR}_{\beta,P_X}$, $\rho_{P_{Y|X}}^{(2)}= \text{CVaR}_{\beta,P_{Y|X}}$, see  Example~\ref{example:nested-cvar} on nested-CVaR.
The corresponding functions are $u_1(w) = u_2(w) = \frac{1}{1-\beta}x$ for $x \leq 0$ and $u_1(w) = u_2(w) = 0$ for $x > 0$.
Then $(-u_1)^*(v) = 0$ if $v\in [-\frac{1}{1-\beta},0]$; $(-u_1)^*(v) = +\infty$ otherwise.
In this case, problem \eqref{eqn:nested-oce} reduces to
{\small
\begin{eqnarray} \label{eqn:nested-cvar-robust1}
    &&\min_{g,s_1,s_2(\cdot)} \sup_{v(X)\in [-\frac{1}{1-\beta},0]}  \mathbb{E}_{P}\Bigg[\Big( s_2(X) -\frac{1}{1-\beta}\big(c(g(X),Y)+ s_2(X)\big)_+ \Big)v(X) - (s_1 + s_1 v(X)) \Bigg]\nonumber\\
    &=&\min_{g,s_1,s_2(\cdot)} \sup_{v(X)\in [0, \frac{1}{1-\beta}]}  \mathbb{E}_{P}\Bigg[ \Big(- s_2(X) +\frac{1}{1-\beta}\big(c(g(X),Y)+ s_2(X)\big)_+ \Big)v(X) - (s_1 - s_1 v(X)) \Bigg]\nonumber\\
    &=& \min_{g,t_2(\cdot)} \sup_{\substack{v(X)\in [0, \frac{1}{1-\beta}]\\\mathbb{E}_{P_X}[v(X)]=1} } \mathbb{E}_{P}\Bigg[ \Big( t_2(X) +\frac{1}{1-\beta}\big(c(g(X),Y)-t_2(X)\big)_+ \Big)v(X) \Bigg].
\end{eqnarray}}
In this way, we can evaluate the objective function of problem \eqref{eqn:nested-cvar-obj} in Example~\ref{example:nested-cvar} based on the joint distribution $P$.
 Since the objective function of problem \eqref{eqn:nested-cvar-robust1} is linear in $v(\cdot)$ and convex in $t_2(\cdot)$, and the feasible set of $v(\cdot)$ is convex,  we can use a minimax theorem (\cite{fan1953minimax,sion1958general})
 to interchange the order of the minimization over $t_2(\cdot)$ with the supremum over $v(\cdot)$:
\begin{eqnarray*}
    &&\min_{g} \sup_{\substack{v(X)\in [0, \frac{1}{1-\beta}]\\\mathbb{E}_{P_X}[v(X)]=1} } \min_{t_2(\cdot)}\mathbb{E}_{P}\Bigg[ \Big( t_2(X) +\frac{1}{1-\beta}\big(c(g(X),Y)-t_2(X)\big)_+ \Big)v(X) \Bigg]\\
    &=&\min_{g} \sup_{\substack{v(X)\in [0, \frac{1}{1-\beta}]\\\mathbb{E}_{P_X}[v(X)]=1} } \mathbb{E}_{P_X}\Bigg[ \min_{t_2}\mathbb{E}_{P_{Y|X}}\Big[ t_2 +\frac{1}{1-\beta}\big(c(g(X),Y)-t_2\big)_+ \Big] v(X) \Bigg]\\
    &=& \min_{g} \sup_{\substack{v(X)\in [0, \frac{1}{1-\beta}]\\\mathbb{E}_{P_X}[v(X)]=1} } \mathbb{E}_{P_X}\Bigg[ \text{CVaR}_{\beta,P_{Y|X}}\big(c(g(X),Y)\big) v(X) \Bigg],
\end{eqnarray*}
where the first equality follows from the interchangeability principle to swap the positions of $\mathbb{E}_{P_X}$ and $\min_{t(\cdot)}$, and the second equality comes from the definition of CVaR in \eqref{def:cvar}.
By setting $v(\cdot)$ equal to the Radon–Nikodym derivative $dQ_X/dP_X$, we obtain
\begin{eqnarray} \label{eqn:nested-cvar-robust}
    \min_{g\in\widetilde{\mathcal{G}}} \sup_{Q_X \in \mathcal{Q}(P_X) } \mathbb{E}_{Q_X}\Bigg[ \text{CVaR}_{\beta,P_{Y|X}}\big(c(g(X),Y)\big) \Bigg],
\end{eqnarray}
where $\mathcal{Q}(P_X):= \{Q_X\ll P_X:0\leq \frac{dQ_X}{dP_X}\leq \frac{1}{1-\beta}\}$.
This formula also follows from a robust/dual representation of CVaR, i.e., $\text{CVaR}_{\beta,P_\xi}(\xi) = \sup_{Q\in \mathcal{Q}(P_\xi)} \mathbb{E}_Q[\xi]$; see e.g., \cite{follmer2010convex}.
Let $(g^*,Q_X^*)$ be an optimal solution to problem \eqref{eqn:nested-cvar-robust}.
In the case that there is
$\hat{\mathcal{X}}\subset \mathcal{X}$ such that
$Q_X^*(\hat{\mathcal{X}})=0$ and
$P_X(\hat{\mathcal{X}})>0$,
$g^*(x)$ does not have an impact on the objective function for $x\in \hat{\mathcal{X}}$.
Such $g^*(x)$ is not necessarily an optimal solution to problem \eqref{eqn:conditional_risk_min}.
This observation is consistent with our observation in Example~\ref{example:nested-cvar}.
\end{itemize}
\end{example}
}

\begin{example} [Entropic risk measure and expected entropic risk measure] \label{example:expected-ent}
    In Example~\ref{example:shortfall-risk}, we
    consider the case that
    $
    \rho_{P_X}^{(1)}=
    \rho_{P_{Y|X}}^{(2)}
=\rho_\gamma^{\text{ent}}
$
with $\gamma>0$.
Since the entropic risk measure is strictly monotone,
Theorem \ref{thm: optimality-of-ent} can be
applied
and
the optimal policy $g^*$ obtained by solving
    \begin{eqnarray} \label{eqn:example4.4ent}
        \min\limits_{g\in\widetilde{\mathcal{G}}} \rho_{\gamma,P_X}^{\text{ent}}\left( \rho_{\gamma,P_{Y|X}}^{\text{ent}} (c(g(X),Y)) \right) = \min\limits_{g\in\widetilde{\mathcal{G}}} \frac{1}{\gamma}\ln\left(\mathbb{E}_P\left[e^{\gamma\left(c(g(X),Y)\right)}\right]\right)
    \end{eqnarray}
    is optimal in conditional risk minimization problem (\ref{eqn:conditional_risk_min}) in all contexts almost surely.
    On the other hand, since the entropic risk measure can be represented in the form of optimized certainty equivalent by setting $u(x) = \frac{1}{\gamma}\left(1-e^{-\gamma x}\right)$, then by Example \ref{example:expected-oce} the optimal policy $\hat{g}^*$ obtained by solving
    \begin{eqnarray} \label{eqn:example4.4eent}
        &&\min\limits_{g\in\widetilde{\mathcal{G}}} \mathbb{E}_{P_X}\left[ \rho_{\gamma,P_{Y|X}}^{\text{ent}} (c(g(X),Y)) \right] \nonumber\\
        &=& \min\limits_{g\in\widetilde{\mathcal{G}},t \in \widetilde{\mathcal{T}}} \mathbb{E}_P\left[t(X)- \frac{1}{\gamma}\left(1-e^{\gamma\left(c(g(X),Y)-t(X)\right)}\right)\right]
    \end{eqnarray}
    is also optimal in problem (\ref{eqn:conditional_risk_min}) almost surely.
    Thus,
    the above two problems
    result in the same optimal policies, i.e., $g^*(x) = \hat{g}^*(x)$ for $x\in \mathcal{X}$ almost surely.
    Indeed, this equivalence can be alternatively
    established by observing that $\widetilde{\mathcal{T}}$ in problem (\ref{eqn:example4.4eent}) can be set to any non-empty set.
By definition of expected entropic risk measure,
\begin{eqnarray*}
&&\min\limits_{g\in\widetilde{\mathcal{G}},t \in \widetilde{\mathcal{T}}} \mathbb{E}_P\left[t(X)+ \frac{1}{\gamma}e^{\gamma\left(c(g(X),Y)-t(X)\right)}\right]
    \nonumber\\
    &=& \min\limits_{t \in \widetilde{\mathcal{T}}}\min\limits_{g\in\widetilde{\mathcal{G}}} \mathbb{E}_{P_X}\left[t(X)+  \frac{1}{\gamma e^{t(X)}}\mathbb{E}_{P_{Y|X}}\left[e^{\gamma\left(c(g(X),Y)\right)}\right]\right]\\
    &=& \min\limits_{t \in \widetilde{\mathcal{T}}} \mathbb{E}_{P_X}\left[t(X)+  \frac{1}{\gamma e^{t(X)}}\min_{z\in\mathcal{Z}}\mathbb{E}_{P_{Y|X}}\left[e^{\gamma\left(c(z,Y)\right)}\right]\right],
\end{eqnarray*}
where the last
equality is due to the interchangeability principle.
Therefore, $g\in \widetilde{\mathcal{G}}$ will always minimize the inner conditional risk minimization problem, irrespective of the choice of $\widetilde{\mathcal{T}}$.
By setting $\widetilde{\mathcal{T}}= \{0\}$, problem (\ref{eqn:example4.4eent}) is equivalent to
$
    \min\limits_{g\in\widetilde{\mathcal{G}}} \frac{1}{\gamma}\mathbb{E}_P\left[e^{\gamma\left(c(g(X),Y)\right)}\right].
$
It is easy to see that the optimal solution to this problem is equivalent to that of the problem (\ref{eqn:example4.4ent}).
\end{example}

{
In Examples \ref{example:ex-ante-cvar}-\ref{example:expected-ent}, the negatives of the OCEs defined via the concave utility function $u$, including CVaR and the entropic risk measure, are law-invariant convex risk measures.
When the cost function $c$ satisfies certain continuity and growth conditions and the conditional distribution $P_{Y|X}$ satisfies the continuity conditions as required in Assumption \ref{assumption:well-defined_rho} and Assumption \ref{assumption:convex_finite-minimum},
we can use
Theorem \ref{thm: optimality-of-ent}(i)
to establish continuity of
the conditionally optimal policy $g^*$ in $x$.
In Section \ref{sec:experiment}, we will conduct numerical tests
on the proposed models
in some practical problems.
To ensure these models are
well posed, we will verify the assumptions
in Section \ref{sec:vei-assumption-application}.
In such settings, if the hypothesis space $\widetilde{\mathcal{G}}$ is rich enough to contain $g^*$ and $\rho^{(1)}_{P_X}$ is strictly monotone,
then the optimal policy in problem (\ref{eqn:ex-post-interchange})
 coincides with the conditionally optimal policy $g^*$.
}
This observation allows us to choose $\rho_{P_X}^{(1)}$ as the expectation operator $\mathbb{E}_{P_X}$ for the contextual uncertainty,
and provides guidance on the choice of hypothesis space
which should be dense in the space of continuous functions, such as RKHS generated by universal kernels.
We will come back in Section~\ref{sec:optimal-saa-rkhs}.

\subsection{Extensions to general risk measures/metrics}
We extend this construction method to general risk measures/metrics
and
summarize our examples in
Table~\ref{tab:example}.
We begin
by recalling
the well-known Kusuoka's representation theorem
\cite[Theorem 4]{kusuoka2001law}.

\begin{lemma} [Kusuoka's representation, Theorem 4 in \cite{kusuoka2001law}]
\label{lemma:kusuoka}
    A coherent risk measure $\rho$ is continuous from above and law-invariant if and only if
    \begin{eqnarray*}
        \rho(\xi) = \sup_{\mu\in \mathcal{M}} \int_{[0,1)} \inmat{CVaR}_{\beta}(\xi)\mu(d\beta),
    \end{eqnarray*}
    for some subset $\mathcal{M}$ of the set of all probability measures/distributions on $[0,1)$.
\end{lemma}

\begin{table}
    \centering
    \caption{Risk metrics and their properties. The expected-OCE includes expected-CVaR and expected entropic risk measure. $\circ$ means the problem is solvable if the risk measures/metrics in Example \ref{example:coherent}, \ref{example:spectral}, \ref{example:distortion} satisfy certain conditions.}
    \begin{tabular}{|c|cccc|}
    \hline
        Risk & \makecell[c]{Represented
        in\\ the form of\\ ex ante risk} & \makecell[c]{Represented
        in \\ the form of\\ ex post risk} & \makecell[c]{Optimal\\
        in each context} & Solvability  \\
        \hline
        Ex ante CVaR & $\checkmark$ & $\times$ & $\times$ & $\checkmark$\\
        Nested-CVaR & $\times$ & $\checkmark$ & $\times$ & {\color{purple} $\checkmark$} \\
        Entropic risk measure & $\checkmark$ & $\checkmark$ & $\checkmark$ & $\checkmark$\\
        Expected-OCE & $\checkmark$ & $\checkmark$ & $\checkmark$ &  $\checkmark$ \\
        \makecell[c]{Expected\\coherent risk measure} & $\checkmark$ & $\checkmark$ & $\checkmark$ & $\circ$ \\
        \makecell[c]{Expected\\mean-upper-semideviation} & $\checkmark$ & $\checkmark$ & $\checkmark$ & $\checkmark$ \\
       \makecell[c]{Expected\\spectral risk measure} & $\checkmark$ & $\checkmark$ & $\checkmark$ & $\circ$ \\
        \makecell[c]{Expected\\distortion risk metric} & $\checkmark$ & $\checkmark$ & $\checkmark$ & $\circ$\\
       \makecell[c]{Expected\\quantile deviation}
        & $\checkmark$ & $\checkmark$ & $\checkmark$ & $\checkmark$\\
        \hline
    \end{tabular}
    \label{tab:example}
\end{table}

\begin{example} [Coherent risk measure]
\label{example:coherent}
    Consider a composition of expectation $\mathbb{E}_{P_X}$ and a coherent risk measure $\rho_{P_{Y|X}}$.
    By exploiting Kusuoka's representation theorem (Lemma \ref{lemma:kusuoka}) and the interchangeability principle for Polish spaces (Lemma \ref{lemma:interchange-polish-space}),
we have
\begin{eqnarray*}
   &&\min_{g\in\mathcal{\widetilde{\mathcal{G}}}} \mathbb{E}_{P_X} \left[\rho_{P_{Y|X}}\left(c(g(X),Y)\right)\right]
   \nonumber\\
   &=&
   \min_{g\in\mathcal{\widetilde{\mathcal{G}}}} \mathbb{E}_{P_X}\left[ \sup_{\mu\in \mathcal{M}} \int_{[0,1)} \inmat{CVaR}_{\beta,P_{Y|X}}(c(g(X),Y))\mu(d\beta) \right]\nonumber\\
   &=& \min_{g\in\mathcal{\widetilde{\mathcal{G}}}} \sup_{\tilde{\mu}(\cdot)\in \widetilde{\mathcal{M}}} \mathbb{E}_{P_X}\left[  \int_{[0,1)} \inmat{CVaR}_{\beta,P_{Y|X}}(c(g(X),Y))\tilde{\mu}(d\beta;X) \right] \nonumber\\
    &=& \min_{g\in\mathcal{\widetilde{\mathcal{G}}}}  \sup_{\tilde{\mu}(\cdot)\in \widetilde{\mathcal{M}}} \inf_{\bar{t}\in\bar{\mathcal{T}}}  \int_{[0,1)} \mathbb{E}_P\left[ \bar{t}(X,\beta) + \frac{1}{1-\beta}[c(g(X),Y)-\bar{t}(X,\beta)]_+ \right] \tilde{\mu}(d\beta;X),
\end{eqnarray*}
where
$\mathcal{M}\subseteq \mathscr{P}([0,1))$,
$\widetilde{\mathcal{M}}$
denotes
 the set
 of
 measurable mappings
 $\tilde{\mu}: \mathcal{X}\rightarrow \mathcal{M}$
 such that
 $\widetilde{\mu}(x)\in{\cal M}$, $\forall x\in {\cal X}$,
 and $\bar{{\cal T}}$
denotes
 a
 set of  measurable functions $\bar{t}:\mathcal{X}\times [0,1)\rightarrow\R$.
 \end{example}

The objective function in Example \ref{example:coherent} can be evaluated by the joint distribution $P$.
However, it is difficult to solve this problem since there is a supremum problem over $\widetilde{\mathcal{M}}$.
Two exceptions are the mean-upper-semideviation with order $p=1$ and spectral risk measures (SRM),
which
subsumes a number of important risk measures including Wang’s proportional Hazard transform and Gini’s risk measure.

\begin{example} [Mean-upper-semideviation, MSD] \label{example:mean-upper-semideviation}
The mean-upper-semideviation with order $p$ is defined as
$
    \rho^{\textrm{MSD}}(\xi) := \mathbb{E}[\xi] + \eta \left(
\mathbb{E}\left[ (\xi - \mathbb{E}[\xi])_+^p \right] \right)^{1/p}
$.
$\rho^{\textrm{MSD}}$ is a coherent risk measure if $\eta\in[0,1]$.
By Kusuoka representation (\cite{kusuoka2001law}), the mean-upper-semideviation risk measure with $\eta\in [0,1]$ and order $p=1$ can be represented as
$$
    \rho^{\textrm{MSD}}(\xi)
    = \inf_{t\in \R}\sup_{s\in [0,1]} \mathbb{E}[\xi + \eta s (t-\xi) + \eta(\xi-t)_+].
$$
See Examples 6.23 and 6.41 in \cite{shapiro2014lectures} for details.
By the interchangeability principle, $\min\limits_{g\in\mathcal{\widetilde{\mathcal{G}}}}\mathbb{E}_{P_X}\left[ \rho_{P_{Y|X}}^{\textrm{MSD}}(c(g(X),Y))\right]$ with $\eta\in [0,1]$ and order $p=1$ is equivalent to
    \begin{eqnarray*}
        &&\min\limits_{g\in\mathcal{\widetilde{\mathcal{G}}}}\mathbb{E}_{P_X}\left[ \rho_{P_{Y|X}}^{\textrm{MSD}}(c(g(X),Y))\right]\\
        &=&\min\limits_{g\in\mathcal{\widetilde{\mathcal{G}}}} \mathbb{E}_{P_X}\left[ \inf_{t\in \R}\sup_{s\in [0,1]} \mathbb{E}_{P_{Y|X}}  \left[ c(g(X),Y) + \eta s(t- c(g(X),Y)) + \eta(c(g(X),Y)-t)_+\right]\right]\\
        &=&
        \inf_{g\in\mathcal{\widetilde{\mathcal{G}}},t\in \widetilde{\mathcal{T}}}\sup_{s\in\widetilde{\mathcal{S}}} \mathbb{E}_{P}[c(g(X),Y) + \eta s(X)(t(X)- c(g(X),Y)) + \eta(c(g(X),Y)-t(X))_+],
    \end{eqnarray*}
    where $\widetilde{\mathcal{T}}$ is
    the set of all measurable functions $t:\mathcal{X}\rightarrow \R$,
    and $\widetilde{\mathcal{S}}$ is the set of all measurable functions $s:\mathcal{X}\rightarrow [0,1]$.
\end{example}

Another important exception of the coherent risk measure is the spectral risk measure.
It is defined as the weighted average of the quantile function of a random variable over $[0,1]$, and subsumes a number of important risk measures including Wang’s proportional Hazard transform and Gini’s risk measure.
Spectral risk measure is comonotonic additive \cite{acerbi2002spectral} and is widely used in portfolio problem, e.g., \cite{acerbi2002portfolio,adam2008spectral}.
Since the spectral risk measures can be represented via CVaR by Kusuoka's representation as discussed in Remark 4.4 in \cite{acerbi2002spectral} and Theorem 7 in \cite{kusuoka2001law}, then the corresponding ex ante risk mapping can be constructed by the interchangeability principle in a similar way.

\begin{lemma} [Kusuoka's representation for spectral risk measure]
\label{lemma:kusuoka-spectral}
    A spectral risk measure $\rho$ is continuous from above and law-invariant if and only if
    \begin{eqnarray*}
        \rho(\xi) = \int_{[0,1)} \inmat{CVaR}_{\beta}(\xi)\mu(d\beta),
    \end{eqnarray*}
    for some $\mu\in \mathscr{P}([0,1))$.
\end{lemma}

\begin{example} [Spectral risk measures] \label{example:spectral}
Let $F_\xi(t):=\mathbb{P}(\xi\leq t)$ be the cumulative distribution function of $\xi$,
$F_\xi^{-1}(p):=\inf\{t\in\R:F_\xi(t)\geq p\}$ be a quantile function of $p\in[0,1]$,
and $\sigma :[0,1]\rightarrow \R$ be a non-negative, decreasing function with the normalization property $\int_0^1 \sigma(p)dp=1$. The SRM of $\xi$ is defined as
$
\rho^{\text{SRM}}(\xi):=-\int_0^1 F_\xi^{-1}(p)\sigma(p) dp
$
\cite{acerbi2002spectral}.
    Consider a composition of  $\mathbb{E}_{P_X}$ and a spectral risk measure $\rho_{P_{X|Y}}$. By Kusuoka's representation for SRM (Lemma \ref{lemma:kusuoka-spectral}), we have
\begin{eqnarray*}
&&\min_{g\in\mathcal{\widetilde{\mathcal{G}}}} \mathbb{E}_{P_X}\left[\rho_{P_{Y|X}}^{\text{SRM}}\left(c(g(X),Y)\right)\right]\nonumber\\
    &=& \min_{g\in\mathcal{\widetilde{\mathcal{G}}}} \mathbb{E}_{P_X}\left[ \int_{[0,1)} \inmat{CVaR}_{\beta,P_{Y|X}}(c(g(X),Y))\mu(d\beta) \right]\nonumber\\
    &=& \min_{g\in\mathcal{\widetilde{\mathcal{G}}},\bar{t}\in \bar{\mathcal{T}}} \int_{[0,1)}
     \mathbb{E}_P\left[ \bar{t}(X,\beta) + \frac{1}{1-\beta}[c(g(X),Y)-\bar{t}(X,\beta)]_+ \right] \mu(d\beta),
\end{eqnarray*}
where $\bar{\mathcal{T}}$ is the set of measurable functions $\bar{t}:\mathcal{X}\times [0,1)\rightarrow\R$.

\end{example}

In addition, our framework can also be extended to the distortion risk metrics, which are not necessarily monotone and thus are not included in the classic notion of monetary risk measures.
A distortion risk metric is defined by the signed Choquet integrals
$
    \rho(\xi)=\int_{-\infty}^0(h(\mathbb{P}(\xi \geq p))-h(1)) \mathrm{d} p+\int_0^{\infty} h(\mathbb{P}(\xi \geq p)) \mathrm{d} p,
$
where $h:[0,1]\rightarrow \R$ is a function of bounded variation with
$h(0)=0$.
In the case where $h$ is a concave function, the distortion risk metrics are convex.
The notion of distortion risk metric encompasses many monetary risk measures and deviation measures, such as mean-median deviation, Gini's deviation and inter-CVaR range. We refer interested readers to Table 1 in \cite{wang2020distortion} for more examples.

\begin{lemma}[Theorem 5 in \cite{wang2020distortion}]
\label{lemma:representation-distortion}
    A functional $\rho:L\rightarrow \R$ where $L\subseteq L^1$ is {a convex distortion risk metric} if and only if there exist finite Borel measures $\mu,\nu$ on $[0,1]$ such that
    \begin{eqnarray}
        \rho(\xi) = \int_0^1 \inmat{CVaR}_\beta(\xi)\mathrm{d}\mu(\beta) + \int_0^1 \inmat{CVaR}_\beta(-\xi)\mathrm{d}\nu(\beta).
    \end{eqnarray}
    Moreover, if $\rho$ is increasing, then we can take $\nu=0$.
\end{lemma}

\begin{example} [Distortion risk metrics] \label{example:distortion}
Consider a composition of  $\mathbb{E}_{P_{X}}$ and a convex distortion risk metric $\rho_{P_{Y|X}}$.
By the CVaR-based representation (Lemma \ref{lemma:representation-distortion}), we have
\begin{eqnarray*}
    &&\min_{g\in\mathcal{\widetilde{\mathcal{G}}}} \mathbb{E}_{P_X}\left[\rho_{P_{Y|X}}\left(c(g(X),Y)\right)\right]\\
    &=& \min_{g\in\mathcal{\widetilde{\mathcal{G}}}} \mathbb{E}_{P_X}\left[ \int_{0}^1 \inmat{CVaR}_{\beta,P_{Y|X}}(c(g(X),Y))\mu(d\beta) + \int_{0}^1 \inmat{CVaR}_{\beta,P_{Y|X}}(-c(g(X),Y))\nu(d\beta) \right]\\
    &=& \min_{g\in\mathcal{\bar{\mathcal{G}}},\bar{t}_1,\bar{t}_2\in\mathcal{\tilde{T}}} \int_{0}^1
     \mathbb{E}_P\left[ \bar{t}_1(X,\beta) + \frac{1}{1-\beta}[c(g(X),Y)-\bar{t}_1(X,\beta)]_+ \right] \mu(d\beta)\\
    && \qquad \qquad\quad + \int_{0}^1
     \mathbb{E}_P\left[ \bar{t}_2(X,\beta) + \frac{1}{1-\beta}[-c(g(X),Y)-\bar{t}_2(X,\beta)]_+ \right] \nu(d\beta),
\end{eqnarray*}
where $\bar{\mathcal{T}}$ is the set of measurable functions $\bar{t}:\mathcal{X}\times [0,1]\rightarrow\R$ and $\mu,\nu$ are finite Borel measures on $[0,1]$.
\end{example}

Particular attention is paid to quantile semideviation \cite{ogryczak2002dual}, which is also known as CVaR-deviation \cite{rockafellar2006generalized}.

\begin{example}[Quantile deviation, QDEV]
\label{example:quantile-semide}
    QDEV is defined as $\rho(\xi) = \min\limits_{t\in \R}\mathbb{E}[\max\{\epsilon_1(t-\xi),\epsilon_2(\xi-t)\}]$ with $\epsilon_1,\epsilon_2>0$.
    Consider a composition of  $\mathbb{E}_{P_{X}}$ and QDEV $\rho_{P_{Y|X}}$, we have
    \begin{eqnarray}
        &&\min_{g\in\widetilde{\mathcal{G}}} \mathbb{E}_{P_X}\left[\rho_{P_{Y|X}}\left(c(g(X),Y)\right)\right] \nonumber\\
        &=& \min_{g\in\widetilde{\mathcal{G}},t\in \widetilde{\mathcal{T}}}\mathbb{E}_P[\max\{\epsilon_1(t(X)-c(g(X),Y)),\epsilon_2(c(g(X),Y)-t(X))\}].
        \label{eqn:quantilesemideviation1}
    \end{eqnarray}
    Ruszczy\'nski and Shapiro \cite{ruszczynski2006optimization} show that the QDEV is equivalent to $-\epsilon_1 \mathbb{E}[\xi] + \epsilon_1 \text{CVaR}_\beta(\xi)$ with $\beta = \frac{\epsilon_2}{\epsilon_1 + \epsilon_2}$.
    Then problem (\ref{eqn:quantilesemideviation1}) is equivalent to
    \begin{eqnarray}
    \label{eqn:quantile-semid}
        \min_{g\in\widetilde{\mathcal{G}},t\in \widetilde{\mathcal{T}}} - \epsilon_1\mathbb{E}_P[c(g(X),Y)]+\epsilon_1\mathbb{E}_{P} \left[t(X) + \frac{1}{1-\beta}(c(g(X),Y)-t(X))_+\right].
    \end{eqnarray}
\end{example}

Deviation metrics quantify the variability of random cost.
To minimize the level and variability of the random cost simultaneously, the DM may consider a convex combination of the expectation and the deviation metric in Examples \ref{example:distortion}-\ref{example:quantile-semide}, i.e., $\lambda\mathbb{E}_{P_{Y|X}} + (1-\lambda)\rho_{P_{Y|X}}$ with $\lambda\in (0,1)$.
In general, the objective function for problem (\ref{eqn:quantile-semid}) may not be convex due to  $- \epsilon_1\mathbb{E}_P[c(g(X),Y)]$. Alternatively, we can ensure  convexity of the problem by considering $\lambda\mathbb{E}_{P_{Y|X}} + (1-\lambda)\rho_{P_{Y|X}}$ and setting $\lambda \geq \frac{\epsilon_1}{1+\epsilon_1}$.

\section{Risk-averse contextual optimization in
reproducing kernel Hilbert space
}
\label{sec:DR-rkhs}

{\color{black}
An important
ingredient of
the contextual optimization is to
choose an appropriate
hypothesis
space $\widetilde{\mathcal{G}}$
for the decision policy.
In Section \ref{sec:SAA}, we discuss the consistency of the SAA method for general hypothesis spaces.
In this section, we move on to consider
the case that the hypothesis space is
a reproducing kernel Hilbert space (RKHS),
and this will enable us to
arbitrarily approximate any real-valued function as the sample
size increases.
}
Let's start with a brief introduction of RKHS.

\subsection{Reproducing kernel Hilbert space}

Let $\mathcal{H}(\mathcal{X})$ be a Hilbert space of functions mapping from $\mathcal{X}$ to $\R$ with inner product $\langle \cdot,\cdot \rangle$.
Let $k:\mathcal{X}\times \mathcal{X} \rightarrow \R$ be a kernel over the support set $\mathcal{X}$ of covariate, that is, there is a feature map $\Phi:\mathcal{X}\rightarrow \mathcal{H}(\mathcal{X})$ such that $k(x,x) = \langle \Phi(x),\Phi(x)\rangle$.

\begin{definition} \label{def: RKHS}
    $\mathcal{H}(\mathcal{X})$ is said to be a reproducing kernel Hilbert space (RKHS) if there is a kernel $k$ such that: (a) $k(\cdot,x)\in \mathcal{H}(\mathcal{X})$ for all $x\in \mathcal{X}$ and (b) $f(x)=\langle f,k(\cdot, x)\rangle$ for all $f\in \mathcal{H}(\mathcal{X})$ and $x\in \mathcal{X}$. The corresponding norm is denoted by {$\|\cdot\|_\mathcal{H}:= \sqrt{\langle \cdot,\cdot \rangle}$}.
    The corresponding RKHS $\mathcal{H}(\mathcal{X})$ is said to be generated by the kernel $k$.
\end{definition}

A kernel is symmetric if $k(x',x'')= k(x'',x')$ for $x',x''\in \mathcal{X}$. A kernel is positive semi-definite if for any finite set $\{x_1,\dots,x_m\}\subset \mathcal{X}$, the $m\times m$ matrix $k[x]$ whose {$(i,j)$-th entry $k(x_i,x_j)$ is positive semi-definite}.
Let $k:\mathcal{X}\times \mathcal{X} \rightarrow \R$ be a Mercer kernel, i.e., a continuous, symmetric, semi-definite kernel,
{and
$$\mathcal{H}' = \left\{ \sum_{i=1}^N \alpha_i k(x_i,\cdot): N=1,2,\dots, \;
\alpha_i\in\R,\; x_i\in \mathcal{X} , \forall i=1,\dots,N \right\}$$
with the inner product $$\left\langle \sum_{i=1}^N \alpha_i k(x_i,\cdot), \sum_{j=1}^N \beta_j k(x_j,\cdot) \right\rangle = \sum_{i,j=1}^N \alpha_i \beta_j k(x_i,x_j)$$ for $N=1,2,\dots$,
then $\mathcal{H}'$ can be completed into the RKHS, i.e., the closure of $\mathcal{H}'$ is an RKHS. See Section 4.1.2 in \citep{boucheron2005theory} for details.
}

For the minimization problem in RKHS, the well-known representer theorem states that a minimizer of a regularized objective function over RKHS can be represented as a finite linear combination of $k(\cdot,x_i), i=1,\dots,N$.
Here we use a multidimensional version of the representer theorem with an inseparable objective function, which is an extension of the version established in \cite{bertsimas2022data,scholkopf2001generalized,shafieezadeh2019regularization}.
{We denote the Cartesian product of $d$ copies of $\mathcal{H}$ by $\mathcal{H}^{d}$, i.e., $\mathcal{H}^{d}=\mathcal{H}\times \cdots \times \mathcal{H}$, and for $f\in\mathcal{H}^{d}$, we denote the $t$-th component of $f(x)$ by $(f(x))^{(t)}$, for $t=1,\dots,d$.}
{
Likewise, for $\alpha\in \R^{N\times d}$,
we denote
the $t$-th column of $\alpha$ by $\alpha^{(t)}\in \R^N$ for $t=1,\dots,d$, and the $i$-th element of $\alpha^{(t)}$ corresponding to sample $x_i$ by $\alpha_i^{(t)}\in \R$ for $i=1,\dots,N$.
}

\begin{lemma}
\label{lemma:multidimensional-representer}
    Fix a sample
    $\{(x_i,y_i)\}_{i=1}^N$ with sample size $N$, $x_i\in\mathcal{X},y_i\in \mathcal{Y}, i=1,\dots,N$.
    Let $\mathcal{H}$ be the RKHS generated by a  kernel $k:\mathcal{X} \times \mathcal{X} \rightarrow {\R}$,
    $W:\left(\mathcal{X} \times \mathcal{Y} \times \R^{d} \right)^N \rightarrow \R$ be an arbitrary cost/loss function,
    {$\mathcal{H}^{d}$ be a hypothesis
    space of continuous
    mappings $f:\mathcal{X}\rightarrow \R^{d}$ with $\|f\|_{\mathcal{H}^{d}} = \sqrt{\sum_{t=1}^d \|f^{(t)}\|_{\mathcal{H}}^2}$,} and
    {$R:\R_+ \rightarrow \R_+$ be a strictly monotonically increasing real-valued  function.}
    Then {any $f\in \mathcal{H}^d$ minimizing the regularized problem}
    \begin{eqnarray}
    \label{eqn:proof-multidim-representer-ec}
       \min_{h\in \mathcal{H}^{d}} W((x_1,y_1,f(x_1)), \dots, (x_N,y_N,f(x_N))) + R(\|f\|_{{\mathcal{H}^{d}}})
    \end{eqnarray}
    {admits a representation of the form}
    \begin{eqnarray} \label{eqn:optimal-solution-form-ec}
   (f^*(x))^{(t)} = \sum_{i=1}^N k(x_i,x)\alpha_i^{(t)},\quad t = 1,\dots,d.
    \end{eqnarray}
    If, in addition, $R(t)=\lambda t^2$,
    with $\lambda> 0$, then $\alpha$ solves
    \begin{eqnarray} \label{eqn:optimal-alpha_form-ec}
        \min_{\alpha^{(1)},\dots,\alpha^{(d)}\in \R^N}&& W((x_1,y_1,(\hat{K} \alpha^{(1)})_1, \dots, (\hat{K} \alpha^{(d)})_1), \dots, (x_N,y_N,(\hat{K} \alpha^{(1)})_N, \dots, (\hat{K} \alpha^{(d)})_N))\nonumber\\
        &&+ \lambda \sum_{t=1}^{d} {{\alpha^{(t)}}^T \hat{K} {\alpha^{(t)}}},
    \end{eqnarray}
    where $\hat{K}\in \R^{N\times N}$ is the matrix with component $\hat{K}_{ij} = k(x_i,x_j)$.
\end{lemma}

\noindent
\textbf{Proof.}
See Section \ref{proof:Lemma4}.
\hfill$\Box$

{
Lemma \ref{lemma:multidimensional-representer}
characterizes the form of the optimal solutions to problem \eqref{eqn:proof-multidim-representer-ec} when the optimal solutions exist.
Moreover, the existence of optimal solutions can be guaranteed under mild additional conditions on the objective function $W$ and the regularizer $R$.
For instance, it is sufficient that $W$ and $R$ are continuous and coercive; see, e.g., \cite[Section 6]{vito2004some} for a detailed discussion.
}
{
After applying the representer theorem in Lemma~\ref{lemma:multidimensional-representer},
problem \eqref{eqn:proof-multidim-representer-ec} reduces to a finite-dimensional regularized SAA problem, which can be solved by standard optimization solvers in numerical test Section~\ref{sec:experiment}.
Gradient descent methods can also be employed for large-scale implementations.
Since this is not the primary focus of this paper, we refer interested readers to e.g., \cite{dai2017learning} for further details.
}

\subsection{Consistency of
SAA under
RKHS}
\label{sec:consist-saa-rkhs}

In data-driven problems, the true distribution of $(X,Y)$ is unknown
whereas the sample data are
available.
This raises a question as to whether
one can use the available sample data to obtain an approximate optimal value and optimal solution of problem
(\ref{eqn:ex_ante}).
Let $\{(x_i,y_i)\}_{i=1}^N$ be
i.i.d. samples of $(X,Y)$ and
$
    P^N(\cdot) = \frac{1}{N} \sum_{i=1}^N \mathds{1}_{(x_i,y_i)}(\cdot),
$
where $\mathds{1}_{(x_i,y_i)}(\cdot)$ denotes Dirac probability measure at $(x_i,y_i)$.
To facilitate the analysis, we concentrate on the case that (\ref{eqn:ex_ante}) takes a particular form
\begin{eqnarray} \label{eqn:saa-form-rkhs-primal}
    \vt(P) &:=& \min_{(z(\cdot),s(\cdot))\in \mathcal{H}^{d_z}\times\mathcal{H}} V(\mathbb{E}_{P}[h(z(X),s(X),Y)]),
\end{eqnarray}
{where $V:\R \rightarrow \R$ is a non-decreasing continuous function, and $h:\R^{d_z}\times \R \times \R^{d_y} \rightarrow\R$ is a continuous function where the first argument is the decision variable $z\in \mathcal{Z}$, the second is the auxiliary variable $s\in \R$ arising from the risk measure, and the third is the problem data $y\in \mathcal{Y}$. For some risk measures, the auxiliary argument is unnecessary, e.g., entropic risk measure. We will soon provide examples to illustrate this model.}
We approximate problem \eqref{eqn:saa-form-rkhs-primal}
by its regularized sample average counterpart:
\begin{eqnarray}
\label{eqn:saa-form-rkhs}
    \vt_N(P_N,\lambda_N) := \min_{(z(\cdot),s(\cdot))\in \mathcal{H}^{d_z}\times\mathcal{H}} V(\mathbb{E}_{P^N}[h(z(X),s(X),Y)])+ \lambda_N R(z,s),
\end{eqnarray}
where
{
$R(z,s) := \|z\|_{\mathcal{H}^{d_z}}^2 + \|s\|_\mathcal{H}^2 = \sum_{t=1}^{d_z} \|z^{(t)}\|_\mathcal{H}^2 + \|s\|_\mathcal{H}^2$ is
 the Tikhonov regularization
 }
 adopted to tackle potential ill-conditionedness of the problem
and $\lambda_N>0$ is a regularization parameter.
By Lemma \ref{lemma:multidimensional-representer}, the optimal solution to problem (\ref{eqn:saa-form-rkhs})
takes the form
(\ref{eqn:optimal-solution-form-ec}).

Especially,
(\ref{eqn:saa-form-rkhs}) subsumes
the two most important risk-averse contextual optimization problems:
risk-averse contextual optimization with an entropic risk measure,
and risk-averse contextual optimization with an expected OCE.
In particular,
if we let  $V(\mathbb{E}_{P^N}[h(z(X),s(X),Y)]) = \frac{1}{\gamma}\ln(\mathbb{E}_{P^N}[h(z(X),s(X),Y)])$ and $h(z,s,y) = e^{\gamma c(z,y)}$,
then
the problem (\ref{eqn:saa-form-rkhs}) reduces to
risk-averse contextual optimization with an entropic risk measure:
\begin{eqnarray*}
    \min_{z^{(1)}(\cdot),\dots,z^{(d_z)}(\cdot)\in \mathcal{H}} \frac{1}{\gamma} \ln\left(\frac{1}{N} \sum_{i=1}^N e^{\gamma c(z^{(1)}(x_i),\dots,z^{(d_z)}(x_i),y_i)} \right) + \lambda_N \sum_{t=1}^{d_z} \|z^{(t)}\|_\mathcal{H}^2.
\end{eqnarray*}
If  $V(\mathbb{E}_{P^N}[h(z(X),s(X),Y)]) = \mathbb{E}_{P^N}[h(z(X),s(X),Y)]$ and $h(z,s,y) = - s - u(-c(z,y)-s)$, then problem (\ref{eqn:saa-form-rkhs}) reduces to risk-averse contextual optimization with expected OCE:
\begin{eqnarray*}
    &\min\limits_{z^{(1)}(\cdot),\dots,z^{(d_z)}(\cdot)\in \mathcal{H},s(\cdot)\in \mathcal{H}}& \frac{1}{N}\sum_{i=1}^N \bigg[-s(x_i) - u\big(-c(z^{(1)}(x_i),\dots,z^{(d_z)}(x_i),y_i)-s(x_i)\big)\bigg] \nonumber\\
    &&+ \lambda_N \left(\sum_{t=1}^{d_z} \|z^{(t)}\|_\mathcal{H}^2 + \|s\|_\mathcal{H}^2\right).
\end{eqnarray*}
Moreover, if $V(\mathbb{E}_{P^N}[h(z(X),s(X),Y)]) = \mathbb{E}_{P^N}[h(z(X),s(X),Y)]$, $h(z,s,y) = c(z,y)$ and $\widetilde{\mathcal{G}}$ is a reproducing kernel Hilbert space, then problem (\ref{eqn:saa-form-rkhs}) reduces to the risk-neutral problem considered in Section 5 in \cite{bertsimas2022data}.

In order to
establish
consistency of problem (\ref{eqn:saa-form-rkhs}),
we need the following assumptions.
Unless specified otherwise,
the support set $\mathcal{X}\times \mathcal{Y}\subseteq \R^{d_x}\times  \R^{d_y}$ may be unbounded,
{
the RKHS $\mathcal{H}$ is induced by a continuous kernel $k$ and therefore consists of continuous functions,
}
and the feasible sets for $z(\cdot)$ and $s(\cdot)$ are bounded in the sense of $\|\cdot\|_{\mathcal{H}^{d_z}}$ and $\|\cdot\|_{\mathcal{H}}$, respectively, for the rest of the paper.
To ease the notation, we denote the bounded feasible sets of $z(\cdot),s(\cdot)$ by $\bar{\mathcal{H}}_z^{d_z},\bar{\mathcal{H}}_s$.

\begin{assumption} \label{assump:saa-convergence-unbounded}
\begin{enumerate}
    \item[(a)]
    There exists
    a continuous function $\phi(x,y): \mathcal{X}\times \mathcal{Y}\rightarrow \R_+$ such that for all $z\in \bar{\mathcal{H}}_z^{d_z}$ and $s\in\bar{\mathcal{H}}_s$,
    $
        |h(z(x),s(x),y)| \leq \phi(x,y), \; \forall (x,y)\in \mathcal{X}\times \mathcal{Y},
    $
    where
     $   \int_{\mathcal{X}\times \mathcal{Y}} \phi(x,y) P(dxdy) < \infty$,
  and $\phi(x,y) \rightarrow \infty$ as $\|(x,y)\|\rightarrow \infty$.

    \item[(b)]
    The feasible set $\bar{\mathcal{H}}_z^{d_z}$
    of $z$
    and
    $\bar{\mathcal{H}}_s$ of $s$ are bounded with {$\|z\|_{\mathcal{H}^{d_z}}=\sqrt{\sum_{t=1}^{d_z} \|z^{(t)}\|_{\mathcal{H}}^2} \leq \beta_1$} for all $z\in \bar{\mathcal{H}}_z^{d_z}$ and $\|s\|_{\mathcal{H}}\leq \beta_2$ for all $s\in \bar{\mathcal{H}}_s$.

    \item[(c)]
    For any compact subset $\mathcal{X}_0\subseteq \mathcal{X}$,
    the set of functions
    $\{k(\cdot,x):x\in \mathcal{X}_0\}$ is
    continuous
    over $\mathcal{X}_0$ uniformly for $x\in {\cal X}$, i.e.,
    for any $\epsilon>0$, there exists a constant $\eta_\epsilon>0$ such that
    $
        \| k(\cdot, x') - k(\cdot, x) \|_\mathcal{H} \leq \epsilon, \forall x,x' \in \mathcal{X}_0 \text{ with } \|x'-x\|\leq \eta_\epsilon.
    $

{
    \item[(d)] For any compact set $\mathcal{Z}_0 \times \mathcal{S}_0\subseteq \mathcal{Z} \times \mathcal{S}$,
    there exist a continuous  function $r:\mathcal{Y}\rightarrow \R_+$
        and a positive constant $\nu\in [0,1)$  such that
        $$
            |h(z_1,s_1,y)-h(z_2,s_2,y)| \leq r(y) \left(\|z_1-z_2\| + \|s_1-s_2\|\right)^\nu,
        $$
        for all $z_1,z_2\in \mathcal{Z}_0, s_1,s_2 \in{\mathcal{S}}_0$, $y\in \mathcal{Y}$.
 }

    \item[(e)]
    The moment generating function $m(t,z,s) := \mathbb{E}_{P}\big[e^{t[h(z(X),s(X),Y) -\mathbb{E}_{P}[h(z(X),s(X),Y)]]}\big]$
    is finite valued for all $z\in\bar{\mathcal{H}}_z^{d_z}$, $s\in \bar{\mathcal{H}}_s$, and $t$ in a neighborhood of zero.

\end{enumerate}
\end{assumption}

Assumption \ref{assump:saa-convergence-unbounded} (a) is known as a growth condition which controls the growth of $h(z(\cdot),s(\cdot),y)$ as $\|(x,y)\|$
increases.
This assumption is trivially satisfied when $\mathcal{X}\times \mathcal{Y}$ is compact.
When $\mathcal{X}\times \mathcal{Y}$ is unbounded, $\phi$ depends on the concrete structure of $h$.
Assumption \ref{assump:saa-convergence-unbounded} does not directly require the equi-continuity, which can instead be established by
more verifiable Assumption \ref{assump:saa-convergence-unbounded} (b) and (c).
The former is widely
used in the literature
of
convergence in RKHS, see e.g.,~\cite{guo2023statistical,norkin2010convergence,zhang2024statistical},
while the latter is satisfied
by the linear kernel, Laplace kernel, Gaussian kernel, polynomial kernel and so on, see e.g.,~\cite{guo2023statistical,shafieezadeh2019regularization}.
{
Assumption (d) requires $h(\cdot,\cdot,y)$ to be H\"older continuous over
any compact subset of $\mathcal{Z}\times \mathcal{S}$.
This assumption can be verified in many specific applications.
For example, in the newsvendor problem with expected CVaR,
all terms in $h(z,s,y)$ are piecewise linear in $(z,s)$ and thus Assumption (d) is satisfied automatically;
in the newsvendor problem with an entropic risk measure, Assumption (d) is satisfied by selecting $r(y)$ in accordance with the feasible set $\mathcal{Z}_0 \times \mathcal{S}_0$.
}
Assumptions (d)-(e) are common assumptions in SAA literature, see e.g.,~\cite{guo2023statistical}.
{We verify Assumptions (a), (d), and (e) for the risk-averse contextual newsvendor and portfolio selection problems in Section \ref{sec:experiment}.}

Now we are ready to state the consistency of SAA method in RKHS.

\begin{theorem} \label{theorem:saa_convergence_rkhs}
    Let Assumption
    \ref{assump:saa-convergence-unbounded} hold. Then the following assertions hold.
    \begin{enumerate}
        \item[(i)] \label{theorem:saa_convergence_rkhs-inner} For any $\delta > 0$, there exist positive constants $\epsilon < \frac{\delta}{6}$, $C(\epsilon,\delta)$ and $\beta(\epsilon,\delta)$, independent of $N$ and a positive number $N_0$ such that
    \begin{eqnarray}
    \label{eqn:saa_unbounded-rkhs-uniform-boundedness-rkhs}
        &&\inmat{Prob}\left(\sup_{(z(\cdot),s(\cdot))\in \bar{\mathcal{H}}_z^{d_z}\times\bar{\mathcal{H}}_s} | \mathbb{E}_{P^N}[h(z(X),s(X),Y)] + \lambda_N R(z,s) - \mathbb{E}_{P}[h(z(X),s(X),Y)] |\geq \delta\right)\nonumber\\
        &\leq& C(\epsilon,\delta)e^{-N{\beta(\epsilon,\delta)}},
    \end{eqnarray}
   for all  $N>N_0$ and $\lambda_N \leq \epsilon/(\beta_1^2+\beta_2^2)$.

    \item[(ii)] If
    $V(t)$ is differentiable and $|V'(t)|$ is upper bounded by $L_V>0$,
    then for any $\delta > 0$, there exist positive constants $\epsilon < \frac{\delta}{6L_V}$, $C(\epsilon,\delta/L_V)$ and $\beta(\epsilon,\delta/L_V)$, independent of $N$ and a positive number $N_0$ such that
$\inmat{Prob}\left(|\vt_N(P^N,\lambda_N)-\vt(P)|\geq \delta\right) \leq C(\epsilon,\delta/L_V)e^{-N{\beta(\epsilon,\delta/L_V)}},
 $
 for all $N>N_0$ and $\lambda_N \leq \epsilon/(\beta_1^2+\beta_2^2)$.

    \item[(iii)] Let $H^*:= \min\limits_{(z(\cdot),s(\cdot))\in \bar{\mathcal{H}}_z^{d_z}\times\bar{\mathcal{H}}_s} \mathbb{E}_{P}[h(z(X), s(X),Y)]$.
    If
    the function $V:\R\rightarrow \R$ takes a special form with $V(t) = \ln t$, $\mathbb{E}_{P^N}[h(z(X),s(X),Y)]>0$, and \linebreak$\mathbb{E}_{P}[h(z(X),s(X),Y)]>0$, then for any $\delta > 0$, there exist positive constants $\epsilon < \frac{\delta H^*}{6(1+\delta)}$,
    $C(\epsilon,\frac{\delta }{1+ \delta}H^*)$ and $\beta(\epsilon,\frac{\delta }{1+ \delta}H^*)$, independent of $N$, and a positive number $N_0$ such that
    \bgeq
        \inmat{Prob}\left(|\vt_N(P^N,\lambda_N)-\vt(P)|\geq \delta\right) \leq C\left(\epsilon,\frac{\delta }{2+ \delta}H^*\right)e^{-N{\beta(\epsilon,\frac{\delta }{2+ \delta}H^*)}},
    \edeq
for all $N>N_0$ and $\lambda_N \leq \epsilon/(\beta_1^2+\beta_2^2)$.
    \end{enumerate}
\end{theorem}

\noindent
\textbf{Proof.}
See Section \ref{proof:Thm2}.
\hfill$\Box$

{
Theorem \ref{theorem:saa_convergence_rkhs} (i) establishes a uniform bound on the difference between the regularized SAA objective $\mathbb{E}_{P^N}[h(z(X),s(X),Y)]+ \lambda_N R(z,s)$ and the true objective $\mathbb{E}_{P}[h(z(X),s(X),Y)]$ over the feasible set $\mathcal{H}^{d_z}\times\mathcal{H}$. This serves as an interim result for parts (ii)–(iii).
Parts (ii)–(iii) provide a bound on the gap between the optimal value of the regularized SAA problem, $\vt_N(P_N,\lambda_N)$, and that of the original problem, $\vt(P)$, but for different risk measures:
part (ii) recovers the expected OCE in Example \ref{example:expected-oce}, and part (iii) recovers the entropic risk measure in Example \ref{example:expected-ent}.
}

It might be helpful to
discuss the difference between the established
convergence results and those in
the literature, such as \cite{bazier2020generalization,bertsimas2022data}.
{
Note that Bertsimas and Koduri \cite{bertsimas2022data} study the risk-neutral objective function $\mathbb{E}_P[c(z(X),Y)]$, which is obtained in the problem (\ref{eqn:saa-form-rkhs}) by choosing $V(\mathbb{E}_{P^N}[h(z(X),s(X),Y)]) = \mathbb{E}_{P^N}[h(z(X),s(X),Y)]$ with $h(z(x),s(x),y) = c(z(x),y)$.
Bazier-Matte and Delage \cite{bazier2020generalization} consider instead an objective based on the certainty equivalent of cost $c(z(x),y)$ using a utility function $u:\R \rightarrow \R$, namely $u^{-1}(\mathbb{E}_P[u(c(z(X),Y))])$, which corresponds to specifying $V = u^{-1}$ and $h(z(x),s(x),y) = u(c(z(x),y))$ in the problem (\ref{eqn:saa-form-rkhs}).
In both \cite{bazier2020generalization,bertsimas2022data}, the auxiliary variable $s$ is superfluous.
To keep our notations consistent and our statements precise, we express the objective in \cite{bazier2020generalization,bertsimas2022data} using a function $\hat{h}:\R^{d_z}\times \R^{d_y}\rightarrow \R$ with $\hat{h}(z(x),y)=c(z(x),y)$ in \cite{bertsimas2022data} and $\hat{h}(z(x),y)=u\big(c(z(x),y)\big)$ in~\cite{bazier2020generalization}.
}

In \cite{bertsimas2022data}, the authors evaluate
the quality of the SAA solution
{$z_{P^N,\lambda_N}^*$} by looking into the difference
$
\mathbb{E}_{P}[{\hat{h}}({z_{P^N,\lambda_N}^*}(X),Y)] - \mathbb{E}_{P^N}[{\hat{h}}({z_{P^N,\lambda_N}^*}(X),Y)]
\label{eq:BK22-BMD20}
$
and they give a probabilistic statement
for the upper bound of the difference, which is known as out-of-performance or
finite sample error.
Here we evaluate the performance of the SAA solution to problem \eqref{eqn:saa-form-rkhs} for the decision variable {$z_{P^N,\lambda_N}^*$} and the auxiliary variable {$s_{P^N,\lambda_N}^*$}
by quantifying the differences
\begin{eqnarray} \label{eqn:diff-saa-true}
    |V(\mathbb{E}_{P^N}[h({z_{P^N,\lambda_N}^*}(X), {s_{P^N,\lambda_N}^*(X)},Y)]) - V(\mathbb{E}_{P}[h({z_{P}^*}(X), {s^*_P(X)},Y)])|,
\end{eqnarray}
{ where $z_P^*$ and $s_P^*$ are the solutions to the problem \eqref{eqn:saa-form-rkhs-primal} under the true distribution~$P$.}

In Proposition 5 in \cite{bertsimas2022data}, the authors  show the convergence of
$\mathbb{E}_{P^N}[{\hat{h}}({z_{P^N,\lambda_N}^*}(X),Y)]$ to $\mathbb{E}_{P}[{\hat{h}}({z_P^*}(X),Y)]$
in probability,
but do not give a rate of convergence.
In Theorem 2 in \cite{bazier2020generalization}, the authors
derive the so-called suboptimality performance bound
between $V\left(\mathbb{E}_{P}[{\hat{h}}({z_{P^N,\lambda_N}^*}(X),Y)]\right)$ and $V\left(\mathbb{E}_{P}[{\hat{h}}({z_P^*}(X),Y)]\right)$, which is
similar to~(\ref{eqn:diff-saa-true}).
Nevertheless, this upper bound has a term
$V'(\mathbb{E}_{P}[{\hat{h}}({z_{P^N,\lambda_N}^*}(x),y)] + \varepsilon)$
which is dependent on the sample and hence is undesirable.
Here, by adopting the
specific form of the risk measure,
we have successfully
avoided such dependence.

Note also that our convergence results are fundamentally based on the Cram\'{e}r’s large deviation theorem, whereas the results in \cite{bazier2020generalization,bertsimas2022data} are derived by virtue of McDiarmid's inequality.
The latter requires the boundedness of the cost function or the boundedness of the support set of the underlying random variables.
For instance,
condition 4 in Proposition 4 in \cite{bertsimas2022data} stipulates  that there exists a constant $C_0\geq 0$ such that
$c(0,y)\leq C_0/2$ and $c(z,y)\geq -C_0/2$,
for all $y\in \mathcal{Y}$ and $z$.
In Assumption 1,2 in \cite{bazier2020generalization},
the support sets $X$ and $Y$ are explicitly assumed to be bounded.
In our convergence results,
we have replaced these
conditions with some weaker or
more verifiable conditions such as Assumption $\ref{assump:saa-convergence-unbounded}$ (a)-(b).

\subsection{Optimality of SAA solution in RKHS with universal kernel}
\label{sec:optimal-saa-rkhs}

{
Section \ref{sec:model-crao} shows that the conditionally optimal policy is continuous in $x$ and can be obtained by solving a risk-averse contextual optimization problem over continuous functions.
In Section \ref{sec:consist-saa-rkhs},
we restrict the feasible set to a bounded RKHS}
and demonstrate convergence of the optimal value
$\vt_N(P_N,\lambda_N)
\to \vt(P)$  as $N\to \infty$.
In the case that $V(t)=t$,
\begin{eqnarray}
\label{eqn:saa-form-rkhs-a}
    \vt_N(P_N,\lambda_N) := \min_{(z(\cdot),s(\cdot))\in \bar{\mathcal{H}}^{d_z}_{z} \times \bar{\mathcal{H}}_s}\mathbb{E}_{P^N}[h(z(X),s(X),Y)]+ \lambda_N R(z,s),
\end{eqnarray}
where the hypothesis space is the Cartesian product of {bounded RKHSs.}
{To extend the result back to the continuous functions, we consider}
the contextual optimization problem
\begin{eqnarray} \label{eqn:optimal-continuous-policy}
    \min\limits_{(z(\cdot),s(\cdot))\in {\mathcal{C}(\mathcal{X})}^{d_z}\times \mathcal{C}(\mathcal{X})} \mathbb{E}_P[h(z(X),s(X),Y)],
\end{eqnarray}
where the hypothesis set is
the space of continuous functions. This raises a question as to whether we can use
problem (\ref{eqn:saa-form-rkhs-a})
to approximate
problem (\ref{eqn:optimal-continuous-policy}). In this subsection, we address this question by restricting
the hypothesis set in (\ref{eqn:saa-form-rkhs-a}) to
RKHSs generated by the universal kernel.
We begin by recalling the definition of universal kernels introduced by  \cite{micchelli2006universal}.

\begin{definition} [Universal kernel]\label{def:universal-kernel}
    A kernel $k: \mathcal{X}\times \mathcal{X}\rightarrow \R$ is said to be universal if for any compact subset $\hat{\mathcal{X}}\subseteq \mathcal{X}$, any continuous function $f\in \mathcal{C}(\hat{\mathcal{X}})$ and any positive constant $\epsilon>0$,
    there exists a function $\hat{f}$ in the RKHS generated by the kernel $k(\cdot,\cdot)$
    such that $\sup_{x\in \hat{\mathcal{X}}}|\hat{f}(x)-f(x)|\leq \epsilon$. Here, $\mathcal{C}(\hat{\mathcal{X}})$ denotes the set of all continuous functions defined on $\hat{\mathcal{X}}$.
\end{definition}

By definition, for any compact subset $\hat{\mathcal{X}}$, the RKHS of a universal kernel is dense in the set of all continuous functions. Detailed analysis and discussions on the universal kernel can be found in \cite{micchelli2006universal,sriperumbudur2011universality}.
In \cite{micchelli2006universal} and Corollary 4.58 in \cite{steinwart2008support}, they list some classes of universal kernels,
including Gaussian kernel, exponential kernel, and binomial kernel.
In particular, Gaussian kernel satisfies Assumption~\ref{assump:saa-convergence-unbounded} as discussed and thus is suitable for risk-averse contextual optimization.
{
Note that
we can use a vector-valued function in the Cartesian product RKHS space, $\mathcal{H}^d:= \mathcal{H}\times \cdots\times \mathcal{H}$, to approximate a vector-valued continuous function by approximating it componentwise,
where each component $\mathcal{H}$ is induced by a universal kernel in Definition \ref{def:universal-kernel}.
Specifically, for any compact subset $\hat{\mathcal{X}}\subseteq \mathcal{X}$, any positive constant $\epsilon>0$ and any $f\in \mathcal{C}^{d}(\mathcal{X})$, there exists a function $\hat{f}\in\mathcal{H}^d$ with the
$t$-th component $\hat{f}^{(t)}$ belonging to the RKHS generated by a universal kernel, such that $\sup_{x\in \hat{\mathcal{X}}} |(\hat{f}(x))^{(t)}-(f(x))^{(t)}|\leq \frac{\epsilon}{\sqrt{d}}$ for $t=1,\dots,d$, and thus $\sup_{x\in \hat{\mathcal{X}}} \|\hat{f}(x)-f(x)\|\leq \epsilon$.
}

The next theorem states that when the RKHS is generated by a universal kernel,
the optimal value of problem \eqref{eqn:saa-form-rkhs-a}
can
approximate that of \eqref{eqn:optimal-continuous-policy} arbitrarily closely as the sample size goes to infinity.

\begin{theorem} \label{theorem:optimal-rkhs}
    Let $\mathcal{X}\times \mathcal{Y}$ be a compact subset of $\R^{d_x}\times \R^{d_y}$.
    Consider an RKHS where $\mathcal{H}$ is generated by a universal kernel $k$ on $\mathcal{X}$, {and $\mathcal{H}^{d_z}$ is a Cartesian product of $d_z$ copies of $\mathcal{H}$.}
    Let $(z^*(\cdot),s^*(\cdot))$ be an optimal solution to (\ref{eqn:optimal-continuous-policy}). Under
    Assumption~\ref{assump:saa-convergence-unbounded}(d) {with $\bar{r}:=\mathbb{E}_P[r(Y)]<\infty$},     the following assertions hold.
    \begin{itemize}
        \item[(i)] For any $\delta_1>0$, there exists $(\hat{z}(\cdot),\hat{s}(\cdot))\in {\mathcal{H}^{d_z}\times\mathcal{H}}$ such that
        \begin{eqnarray*}
        \mathbb{E}_P[h(\hat{z}(X),\hat{s}(X),Y)] -
            \mathbb{E}_P[h(z^*(X),s^*(X),Y)] \leq \delta_1.
        \end{eqnarray*}

        \item[(ii)] If, in addition,
        Assumption \ref{assump:saa-convergence-unbounded} holds with $\beta_1: = \|\hat{z}\|_{\mathcal{H}^{d_z}}$ and $\beta_2: = \|\hat{s}\|_{\mathcal{H}}$ in Assumption \ref{assump:saa-convergence-unbounded} (b), then for any $\delta_2 > 0$, there exist positive constants $N_0$, $\epsilon < \frac{\delta_2}{6}$, $C(\epsilon,\delta_2)$ and $\beta(\epsilon,\delta_2)$, independent of $N$ such that for $N>N_0$ and $\lambda_N \leq \epsilon/(\beta_1^2+\beta_2^2)$,
    \begin{eqnarray*}
        \inmat{Prob}\bigg( \bigg| \vartheta_N(P_N,\lambda_N) - \mathbb{E}_{P}[h(z^*(X),s^*(X),Y)] \bigg|\geq \delta_1 + \delta_2 \bigg)
        \leq \hat{C}(\epsilon,\delta_2,\delta_1)e^{-N{\hat{\beta}(\epsilon,\delta_2,\delta_1)}},
    \end{eqnarray*}
    { }
    \end{itemize}
\end{theorem}

\noindent
\textbf{Proof.}
See Section \ref{proof:Thm3}.
\hfill$\Box$

{
Theorem \ref{theorem:optimal-rkhs} establishes a link between the continuous policy from problem \eqref{eqn:optimal-continuous-policy} and the RKHS policy from problem \eqref{eqn:saa-form-rkhs-a} via a decomposition of errors.
Part (i) leverages the universality of the kernel to control the approximation error, i.e., for any continuous optimal policy $(z^*(\cdot),s^*(\cdot))$ and any tolerance $\delta_1>0$, there is a policy $(\hat{z}(\cdot),\hat{s}(\cdot))$ in the RKHS such that the difference between their optimal values is smaller than $\delta_1$.
Part (ii) then bounds the error by combining the consistency of the SAA method in Theorem \ref{theorem:saa_convergence_rkhs} with the approximation error in part (i).
Consequently, restricting the hypothesis class to an RKHS generated by a universal kernel is asymptotically optimal in value for the original problem \eqref{eqn:optimal-continuous-policy} in the space of continuous functions, which justifies
the RKHS method adopted in Section~\ref{sec:DR-rkhs}.
}

\section{Numerical
tests}
\label{sec:experiment}

To examine the performance of
the proposed risk-averse contextual optimization and the RKHS approach,
we conduct some numerical tests.
{
The numerical tests are conducted using the standard optimization solvers Gurobi and Mosek on an Apple MacBook Air with an M4 chip. The implementation code is available.\footnote{\url{https://github.com/T4OYU4N/Risk-averse-Decision-Making-with-Contextual-Information}}
}

\subsection{Newsvendor problem}
\label{sec:newsvender}

\subsubsection{Models}
\label{sec:newsvendor-model}
Consider a newsvendor problem where the decision maker needs to decide the ordering quantity $z\geq 0$ before a random demand $Y$ is observed.
Let $h = 0.2,b=1$ represent the unit holding cost and the unit backorder cost, respectively.
For the ordering quantity $z$ and the demand $y$, the total cost is computed as
$$
    c(z,y) = h(z-y)_+ + b(y-z)_+,
$$
where $(a)_+ = \max\{a,0\}$.
Let us consider the conditional risk minimization problem \eqref{eqn:conditional_risk_min} with CVaR
\begin{eqnarray}
    && \min_{z\in\R}\inmat{CVaR}_{\beta,P_{Y|X=x}} \left( h(z-Y)_+ + b(Y- z)_+ \right)\nonumber\\
    & = & \min_{t,z\in\R} \mathbb{E}_{P_{Y|X=x}}\bigg[t+\frac{1}{1-\beta}\left[h(z-Y)_+ + b(Y- z)_+ - t\right]_+\bigg], \label{eqn:newsvendor_cond_cvar}
\end{eqnarray}
and the conditional risk minimization problem \eqref{eqn:conditional_risk_min} with entropic risk measure
\begin{eqnarray}
    && \min_{z\in\R} \rho_{\gamma,P_{Y|X=x}}^{\text{ent}} \left( h(z-Y)_+ + b(Y- z)_+ \right)\nonumber\\
    & = & \min_{z\in\R} -\frac{1}{\gamma} \ln  \mathbb{E}_{P_{Y|X=x}} \left[e^{\gamma (h(z-Y)_+ + b(Y- z)_+)}\right].\label{eqn:newsvendor_cond_ent}
\end{eqnarray}
{
To facilitate
discussions, we
call the optimal policies obtained from solving problems \eqref{eqn:newsvendor_cond_cvar} and \eqref{eqn:newsvendor_cond_ent}
the conditional-CVaR policy and conditional-entropic-RM policy, respectively.}
Both policies require the knowledge of the conditional distribution $P_{Y|X}$, which is typically unavailable in data-driven settings.
We present these policies solely as benchmarks for comparison.
Then we examine the risk-averse contextual optimization {problems \eqref{eqn:ex-post-interchange} and \eqref{eqn:ex_ante} under the following choices of risk mappings}:

\begin{itemize}
    \item {Ex-ante-CVaR policy:}
    {problem \eqref{eqn:ex_ante} with ex ante CVaR and $\beta = 0.9$} (Example~\ref{example:ex-ante-cvar}):
\begin{eqnarray} \label{eqn:newsvendor-ex-ante-cvar}
&&\min_{z(\cdot)\in\tilde{\mathcal{G}}} \inmat{CVaR}_{\beta,P} \left( h(z(X)-Y)_+ + b(Y- z(X))_+ \right)\\
&=&\min_{z(\cdot)\in\tilde{\mathcal{G}},t\in \R} \mathbb{E}_P\bigg[t + \frac{1}{1-\beta}\big[h(z(X)-Y)_+ + b(Y- z(X))_+ - t\big]_+\bigg]\nonumber
\end{eqnarray}
where the equation comes from the Rockafellar/Uryasev reformulation in
\cite{rockafellar2000optimization};

\item {Expected-CVaR policy:}
{problem \eqref{eqn:ex-post-interchange} with expected CVaR and $\beta = 0.9$} (Example \ref{example:expected-oce}), i.e., $\mathbb{E}_{P_X}\circ\inmat{CVaR}_{\beta,P_{Y|X}}$:
\begin{eqnarray} \label{eqn:newsvendor-expected-cvar}
&& \min_{z(\cdot)\in\tilde{\mathcal{G}}} \mathbb{E}_{P_X} (\inmat{CVaR}_{\beta,P_Y|X}(h(z(X)-Y)_+ + b(Y- z(X))_+)) \\
    &=& \min_{z(\cdot)\in\tilde{\mathcal{G}}, t(\cdot)\in \tilde{\mathcal{T}}} \mathbb{E}_{P} [t(X)]+ \mathbb{E}_{P}\left[\bigg(h(z(X)-Y)_+ + b(Y- z(X))_+-(t(X))\bigg)_+\right];\nonumber
\end{eqnarray}

\item {Ex-ante-entropic-RM policy:} {problem \eqref{eqn:ex_ante} with entropic risk measure and $\gamma = 0.5$} (Entropic-RM, Example \ref{example:expected-oce}):
\begin{eqnarray}\label{eqn:newsvendor-ex-ante-ent}
&& \min_{z(\cdot)\in\tilde{\mathcal{G}}} \rho_{\gamma,P}^{\text{ent}} \left( h(z(X)-Y)_+ + b(Y- z(X))_+ \right)\\
    & = & \min_{z(\cdot)\in\tilde{\mathcal{G}}} -\frac{1}{\gamma} \ln  \mathbb{E}_P \left[e^{\gamma (h(z(X)-Y)_+ + b(Y- z(X))_+)}\right];\nonumber
\end{eqnarray}

\item {Expected-entropic-RM policy:} {problem \eqref{eqn:ex-post-interchange} with expected entropic risk measure and $\gamma = 0.5$} (Expected Entropic-RM, Example~\ref{example:expected-oce}):
\begin{eqnarray}\label{eqn:newsvendor-expected-ent}
&& \min_{z(\cdot)\in\tilde{\mathcal{G}}} \mathbb{E}_{P_X} (\rho^\text{ent}_{\gamma,P_{Y|X}}(h(z(X)-Y)_+ + b(Y- z(X))_+)) \\
    &=& \min_{z(\cdot)\in\tilde{\mathcal{G}}, t(\cdot)\in \tilde{\mathcal{T}}} \mathbb{E}_{P} [t(X)]-\frac{1}{\gamma}\mathbb{E}_{P}\left[1-e^{\gamma\left(h(z(X)-Y)_+ + b(Y- z(X))_+-(t(X))\right)}\right]. \nonumber
\end{eqnarray}

\end{itemize}

{
These problems satisfy
Assumption~\ref{assumption:well-defined_rho},
Assumption~\ref{assumption:convex_finite-minimum},
and Assumption~\ref{assump:saa-convergence-unbounded}
required for the theoretical developments across Sections \ref{sec:contextual-risk-measure}-\ref{sec:DR-rkhs}.
We verify some of them in
Section \ref{sec:veri-assum-newsvendor}.
}

\subsubsection{Obtaining a conditionally optimal policy in a linear setting.}
\label{sec:experiment-newsvend-linear}
In order to showcase the characteristics of the optimal policies derived
from different risk measures/metrics,
we first consider a simple case where the demand function
is linearly dependent on the contextual covariates $X$.
The covariate $X$
denotes the newsworthiness and follows a log-normal distribution, where its underlying normal component $N(1,0.25)$ is truncated within $\pm 2$ standard deviations from the mean.
The quantity of demand is
$
    Y = \max \{5X + \epsilon+100,0\},
$
where $\epsilon$ is a truncated normal distribution
noise where the normal component
$N(0, (20-x)^2)$ is truncated to $\pm 2$ standard deviations from the mean.
This comes from the fact that when the newsworthiness is lower, the sales are more likely to be influenced by other factors and thus are more random.
To facilitate the presentation of the result, we first restrict the policies to linear decision rules, i.e.,
\begin{eqnarray*}
    \tilde{\mathcal{G}}:= \left\{z(\cdot): z(x)= z_0+z_1x,\; z_0,z_1\in\R \right\}, \quad \tilde{\mathcal{T}}:= \left\{t(\cdot): t(x)= t_0+t_1x,\; t_0,t_1\in\R \right\}.
\end{eqnarray*}

We compare the policies by solving the problem in Section \ref{sec:newsvendor-model} with the risk-neutral optimal policy (Risk-neutral policy) derived by solving the conditional risk minimization problem (\ref{eqn:conditional_risk_min}) with $\rho_{P_{Y|X}}^{(2)}=\mathbb{E}_{P_{Y|X}}$:
$$\inf\limits_{z\in\R} \mathbb{E}_{P_{Y|X=x}}\bigg[h(z-Y)_+ + b(Y- z)_+\bigg],\;\forall x.$$
The closed-form of the risk-neutral policy is
$
    z^*(x) = F^{-1}_{P_{Y|X=x}}\left(\frac{b}{h+b}\right),
$ where $F^{-1}_{P_{Y|X=x}}$ is the quantile function of $P_{Y|X=x}$.
For the truncated normal distribution $N(20 x +100, \left(20-x\right)^2)$ limited to $\pm 2$ standard deviations, we have
$F^{-1}_{P_{Y|X=x}}\left(p\right)= \mu_{Y|X=x} + \sigma_{Y|X=x} z_{p'} = (20x + 100) + (20-x) z_{p'}$, where $z_{p'}$ is the $p'$-quantile of the standard normal distribution, $p' = (2\Phi(2)-1)p + \Phi(-2)$ and $\Phi(\cdot)$ is the cumulative distribution function of the standard normal distribution.
Therefore, the risk-neutral policy $z^*(x)$ is linear in $x$.

We first compare the {ex-ante-CVaR policy derived from problem \eqref{eqn:newsvendor-ex-ante-cvar} and expected-CVaR policy from problem \eqref{eqn:newsvendor-expected-cvar} with}
{conditional-CVaR policy}
from problem \eqref{eqn:newsvendor_cond_cvar}  with $\beta = 0.9$ for almost every $x\in \mathcal{X}$.
The closed-form of conditional-CVaR policy
is
\begin{eqnarray*}
    z^{*}(x) &=& \frac{h}{h+b}F_{Y|X=x}^{-1}\left(\frac{b(1-\beta)}{h+b}\right) + \frac{b}{b+h}F_{Y|X=x}^{-1}\left(\frac{h\beta+b}{b+h}\right),\\
    t^*(x) &=& \frac{hb}{h+b} F_{Y|X=x}^{-1}\left(\frac{h \beta+b}{h+b}\right)-\frac{bh}{h+b} F_{Y|X=x}^{-1}\left(\frac{b(1-\beta)}{h+b}\right),
\end{eqnarray*}
see \cite{gotoh2007newsvendor}.
Therefore, as stated above,
both $z^*(x)$ and $t^*(x)$ are linear in $x$.

Similarly, we compare the {ex-ante-entropic-RM policy from problem \eqref{eqn:newsvendor-ex-ante-ent} and expected-entropic-RM policy from problem \eqref{eqn:newsvendor-expected-ent} with} the
{conditional-entropic-RM policy} from problem \eqref{eqn:newsvendor_cond_ent} with $\gamma = 0.5$.
This policy does not have a closed form, but simulation indicates that it is nearly linear in $x$ in our setting.
(The R-squared value of the linear regression between covariate $X$ and the conditional-entropic-RM policy $z^*(X)$ is $1.00$ for the simulated data, indicating a nearly linear relationship.)

\begin{figure}[htb]
    \centering
    \includegraphics[width=0.9\linewidth]{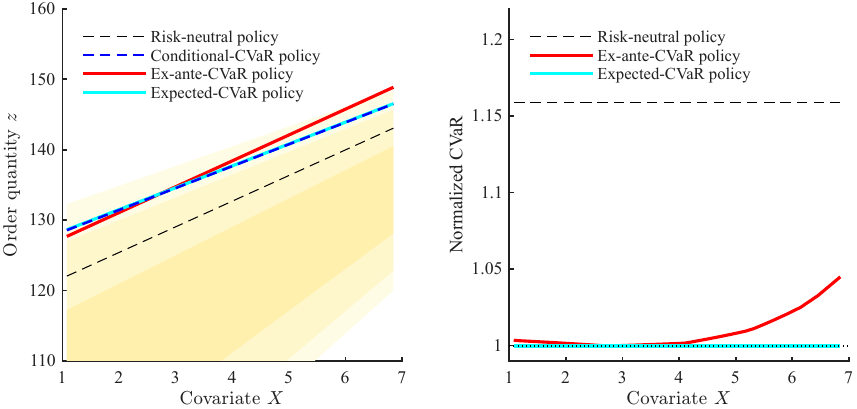}
    \caption{Ex-ante-CVaR policy and expected-CVaR policy, together with the CVaR of their losses in each context normalized by that of the conditional-CVaR policy.
    The yellow area and shades indicate the probability density of the demand $Y$ conditional on $X$. }
    \label{fig:experiment1-linear-cvar}
\end{figure}

\begin{figure} [H]
    \centering
    \includegraphics[width=0.9\linewidth]{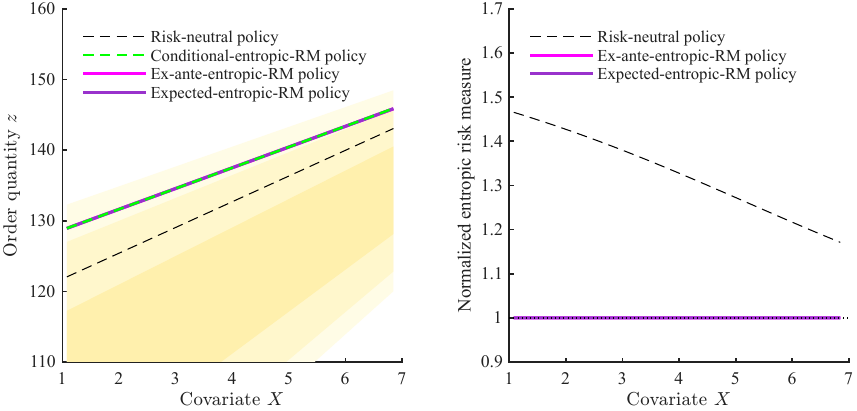}
    \caption{Ex-ante-entropic-RM policy and expected-entropic-RM policy, together with the entropic risk measure of their losses in each context normalized by that of the conditional-entropic-RM policy.
    The yellow area and shades indicate the probability density of the demand $Y$ conditional on $X$.
    Note that the ex-ante-entropic-RM policy and the expected-entropic-RM policy almost completely overlap in both figures.}
    \label{fig:experiment1-linear-ent}
\end{figure}

{
Fig. \ref{fig:experiment1-linear-cvar}-\ref{fig:experiment1-linear-ent} depict these policies along with their respective distances to the corresponding conditionally optimal policies.
}
We have the following observations:
\begin{enumerate}
{
    \item[(i)]
    Since the unit holding cost $h=0.2$ is lower than the unit backorder cost $b=1$, a more conservative policy tends to choose a larger order quantity.
    This is reflected in the left panels of Fig.~\ref{fig:experiment1-linear-cvar}-\ref{fig:experiment1-linear-ent}, where all risk-averse policies lie above the risk-neutral policy (black dashed line).

    \item[(ii)] In Fig.~\ref{fig:experiment1-linear-cvar},
    the ex-ante-CVaR policy differs substantially from the conditional-CVaR policy (left panel)
    and incurs a large normalized CVaR of loss in most contexts (right panel), which indicates that the ex-ante-CVaR policy is not conditionally optimal in most contexts.

    \item[(iii)] In contrast,
    the conditional-CVaR policy (blue dashed) and the expected-CVaR policy (blue solid) {nearly overlap} in the left panel of Fig.~\ref{fig:experiment1-linear-cvar}, and the normalized CVaR of loss is nearly $1$ for almost every $x\in \mathcal{X}$ in the right panel.
    Similarly, in Fig.~\ref{fig:experiment1-linear-ent},
    the conditional-entropic-RM policy (green dashed), ex-ante-entropic-RM policy (purple solid), and expected-entropic-RM policy (dark purple solid) {nearly overlap} in the left panel; and the normalized entropic risk measure of loss is nearly~$1$ in the right panel.
    This observation is consistent with our analysis in Section \ref{sec:model-crao},
     suggesting that the expected-CVaR policy, ex-ante-entropic-RM policy,
    and expected-entropic-RM policy
    can achieve the conditional optimality.

\item[(iv)]
    The yellow region
(which represents the spread of $Y$) narrows down as $x$ increases, indicating that the variance of $Y$ conditional on $X=x$
decreases with respect to~increase of $x$.
 However, when $x$ is large,
 the order quantity under the ex-ante-CVaR policy is larger than that under the conditional-CVaR policy.
 Given that the unit cost of shortage (with $b=1$) is  higher than
 the unit cost
 of unit overage holding costs
 $h=0.2$,
this indicates a relatively more conservative strategy against shortage.
    By contrast, when $x$ is small, the variance of $Y$ is high, yet the ex-ante-CVaR policy chooses a smaller order quantity than the conditional-CVaR policy,
    which indicates a more aggressive strategy against oversupply (holding cost).
    This outcome conflicts with the intuition that higher variance should call for a more conservative policy.}
\end{enumerate}

    {
    The outcome observed in $(iv)$ can be explained as follows.
Consider the conditional risk minimization problem \eqref{eqn:newsvendor_cond_cvar} with $x\in \mathcal{X}$, $h= 0.2, b= 1$ and fixed $t\in \left[0, \frac{2}{3}\sigma_{Y|X=x}\right]$.
We have
\begin{eqnarray*}
    &&\min_{z\in\R} \mathbb{E}_{P_{Y|X=x}}\bigg[t+\frac{1}{1-\beta}\left[h(z-Y)_+ + b(Y- z)_+ - t\right]_+\bigg]\\
    &=& \min_{z\in\R}  t + \frac{1}{1-\beta}\int_{-\infty}^{z-5t} \left[0.2(z-y)-t \right] dF_{Y|X=x}(y) + \int_{z+t}^\infty \left[ (y-z) - t \right] dF_{Y|X=x}(y)
\end{eqnarray*}
where $F_{Y|X=x}(y)$ denotes the cumulative distribution function of $Y|X=x$, and the equality follows from the fact that
$\left[ h(z-y)_+ + b(y-z)_+ - t \right]_+ $
    is positive
    only when $y< z- t/h = z-5t$ or $y > z+ t/b = z+ t$.
Since $P_{Y|X=x}$ is a continuous distribution,
the first-order optimality condition is
\begin{eqnarray} \label{eqn:newsvendor-opt}
    0.2 F_{Y|X=x}(z(x,t)-5t) - 1+ F(z(x,t)+t) = 0,
\end{eqnarray}
where $z(x,t)$ emphasizes that the optimal solution depends on both $x$ and $t$.
Differentiating both sides of the first-order optimality condition w.r.t. $t$, we obtain
\begin{eqnarray*}
    0.2f_{Y|X=x}(z(x,t)-5t) \left( z_t'(x,t) -5 \right) + f_{Y|X=x} (z(x,t)+t) (z_t'(x,t)+1) = 0.
\end{eqnarray*}
Therefore,
\begin{eqnarray*}
    z_t'(x,t) = \frac{f_{Y|X=x}(z(x,t)-5t) - f_{Y|X=x} (z(x,t)+t)}{0.2f_{Y|X=x}(z(x,t)-5t) + f_{Y|X=x} (z(x,t)+t)}.
\end{eqnarray*}
Thus, if $f_{Y|X=x}(z(x,t)-5t) > f_{Y|X=x} (z(x,t)+t)$, the optimal solution $z(x,t)$ increases in $t$.
Note that the conditional distribution $P_{Y|X=x}$ is a truncated normal distribution, which is symmetric and unimodal.
By the optimality condition \eqref{eqn:newsvendor-opt}, namely,
\begin{eqnarray*}
    0.2 \mathbb{P}_{Y|X=x} \left( Y< z(x,t)-5t \right) = \mathbb{P}_{Y|X=x} \left( Y \geq z(x,t)+t \right),
\end{eqnarray*}
we have
$\left| \big( z(x,t)-5t \big) - \mu_{Y|X=x} \right| < \left| \big( z(x,t) + t \big) - \mu_{Y|X=x} \right|$.
Since the density of a symmetric unimodal distribution decreases as the distance from its expectation increases, it follows that
\begin{eqnarray*}
    f_{Y|X=x}(z(x,t)-5t) > f_{Y|X=x} (z(x,t)+t).
\end{eqnarray*}
Consequently, for fixed $x$ and any $t\in \left[0, \frac{2}{3}\sigma_{Y|X=x}\right]$, the optimal solution $z(x,t)$ always increases in $t$.

On the other hand, Rockafellar and Uryasev \cite{rockafellar2000optimization} show that {the optimal threshold $t$ in CVaR reformulation corresponds to VaR}.
This implies that, when \(P_{Y|X=x}\) has larger variance, the conditional risk minimization problem with CVaR tends to yield a larger optimal threshold \(t\).
For example, when $x=1$, the variance of the underlying normal distribution is $19^2$, and the numerical experiment shows that the optimal threshold is $t=9.41$.
Conversely, when \(P_{Y|X=x}\) has smaller variance, it tends to yield a smaller optimal threshold \(t\).
This corresponds to the case when $x=7$, the variance is $13^2$, and thus the optimal threshold in experiments is $t=4.83$.

In contrast, the risk-averse contextual optimization problem with ex-ante CVaR accounts for all contexts and employs a moderate threshold \(t=7.98\) by the numerical test.
Therefore, when the actual conditional risk is high, e.g., when $x=1$, the ex-ante CVaR tends to underestimate risk, i.e. $t=7.98<9.41$, and thus selects a policy with a smaller order quantity.
Conversely, when the actual conditional risk is low, e.g., when $x=7$, it tends to overestimate risk, i.e., $t=7.98>4.83$, and consequently selects a policy with a larger order quantity.
These results are consistent with the observations in Fig. \ref{fig:experiment1-linear-cvar}.

}

\subsubsection{Performance of the SAA method in RKHS for the nonlinear case with two-dimensional contextual features.}

Next, we study the performance of
the SAA method in RKHS.
In this part,
we assess the performance of RKHS generated by Gaussian kernel and consider the case where demanding quantity is not linearly dependent on covariates:
$$
    Y = \max\{5\inmat{exp}(1+0.5X_1) + 20 \cos (2X_2) + (X_1+X_2+12)\epsilon + 100, 0\} ,
$$
where $X_1,X_2,\epsilon$ are independent random variables, each following a truncated standard normal distribution $N(1,0.5^2)$ limited to within $\pm2$ standard deviations from the mean.
The exponential dependence is considered in \cite{bertsimas2022data} and sine/cosine dependence is considered in \cite{zhang2024optimal}.
We first focus on the risk-averse contextual newsvendor problem \eqref{eqn:newsvendor-expected-cvar} with expected CVaR.

\begin{figure} [tb]
    \centering
    \includegraphics[width=\linewidth]{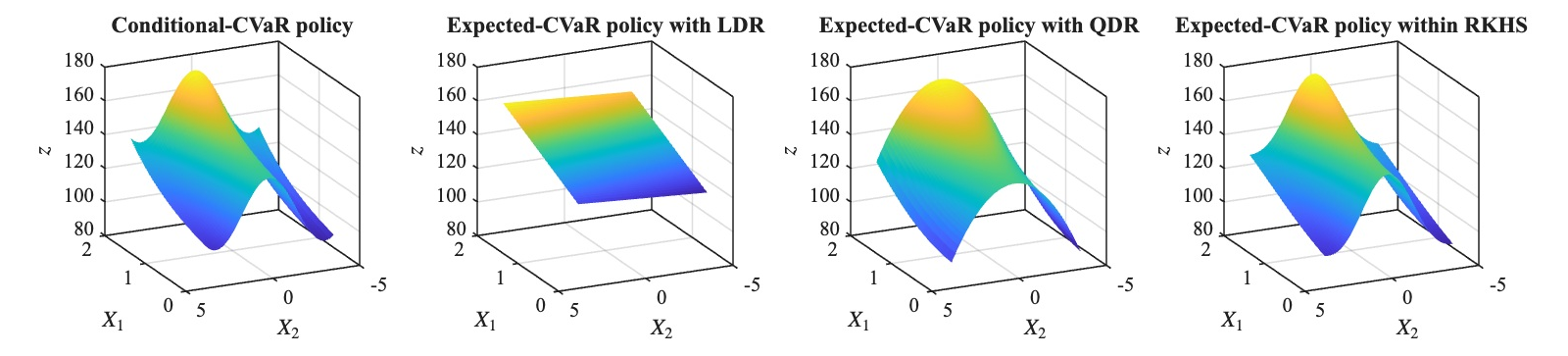}
    \caption{Conditional-CVaR policy and Expected-CVaR policies with hypothesis space being LDR, QDR, and RKHS (training sample size $N=800$).  }
    \label{fig:ex1-nonlinear-policy-cvar}
\end{figure}

\begin{figure} [tb]
    \centering
    \includegraphics[width=0.95\linewidth]{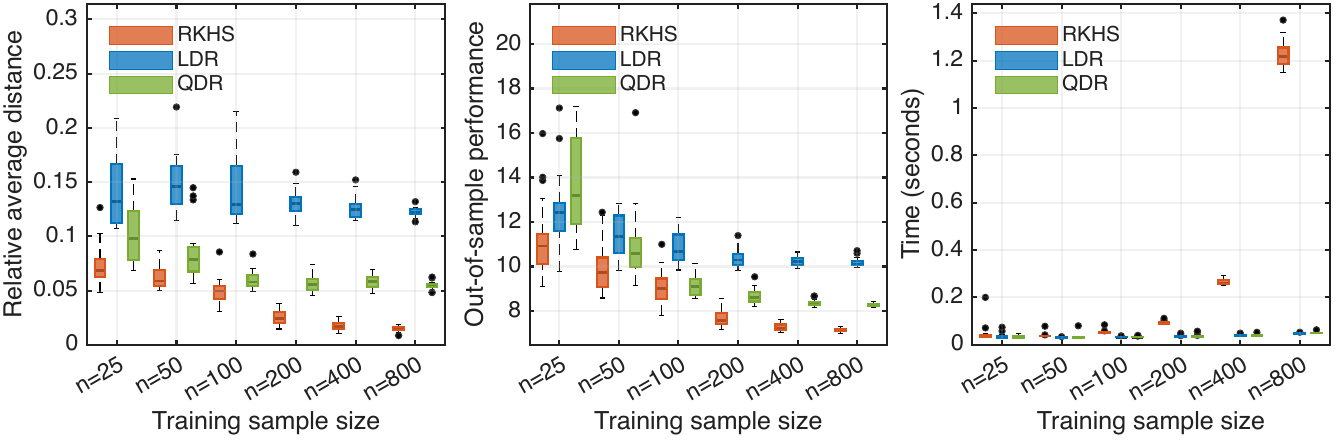}
    \caption{Performance of expected-CVaR policy within various hypothesis spaces: conditional-CVaR policy, RKHS, LDR, and QDR.}
    \label{fig:ex1-nonlinear-boxplot-cvar}
\end{figure}

We vary the number of data points $N\in \{25,50,100,200,400,800\}$ and conduct
20 replicate experiments for each sample size.
{
For expected-CVaR policy within RKHS, the kernel bandwidth is chosen as \(h_{N,d} = h_{\text{ref}} \left(\frac{10}{N} \right)^{1/d+4}\), which corresponds to the classical result found, e.g., in \cite{hall2004cross,stone1982optimal},
the regularization parameter for the decision policy $z(\cdot)$ is set to \(\lambda_N^z = \lambda_{\text{ref}}^z \left(\frac{10}{N} \right)^{1/3}\), and
the regularization parameter for auxiliary function $t(\cdot)$ is set to \(\lambda_N^t = \lambda_{\text{ref}}^t \left(\frac{10}{N} \right)^{1/3}\), which ensures that the regularization strength decreases with the sample size, as discussed in Section~\ref{sec:consist-saa-rkhs}.
For these parameters, the values of $h_{\text{ref}}, \lambda_{\text{ref}}^z,$ and $\lambda_{\text{ref}}^t$ are selected from a set of candidate values; in the end, we choose $h_{\text{ref}}=5, \lambda_{\text{ref}}^z=10^{-6}, \lambda_{\text{ref}}^t=0.8$.
}
We also assess the performance of the expected-CVaR policy with linear decision rule (LDR) and quadratic decision rule (QDR).
We plot the conditional-CVaR policy (for comparative purposes only) and the expected-CVaR policies within these hypothesis spaces in Fig.~\ref{fig:ex1-nonlinear-policy-cvar}.
Fig.~\ref{fig:ex1-nonlinear-boxplot-cvar} displays boxplots illustrating the performance of the different hypothesis spaces across different training sample sizes, including the relative average distance to the conditional-CVaR policy (relative average distance), the out-of-sample performance on the testing data (out-of-sample-performance), and the computation times (times (seconds)).

As shown in Fig.~\ref{fig:ex1-nonlinear-policy-cvar},
since the conditional-CVaR policy exhibits significant nonlinearity, the performance of LDR is poor in this setting,
while both the QDR and policies in RKHS
exhibit
much better performances.
Fig.~\ref{fig:ex1-nonlinear-boxplot-cvar}
shows that
the relative average distance for the LDR policies consistently remains above $10\%$.
In contrast, as the sample size increases, the relative average distance for the QDR policies drops to approximately $5\%$,
whereas the relative average distance for the RKHS policies quickly converges to roughly $1-2\%$.
Similarly, the out-of-sample performance of LDR policies remains above $10$, that of QDR policies decreases to $9$, and that of RKHS policies converges to $8$.
These findings show that RKHS policies achieve strong performance and generalize effectively in this nonlinear setting.
Nevertheless, the RKHS-based method tends to require a longer computation time.

\begin{figure} [tb]
    \centering
    \includegraphics[width=\linewidth]{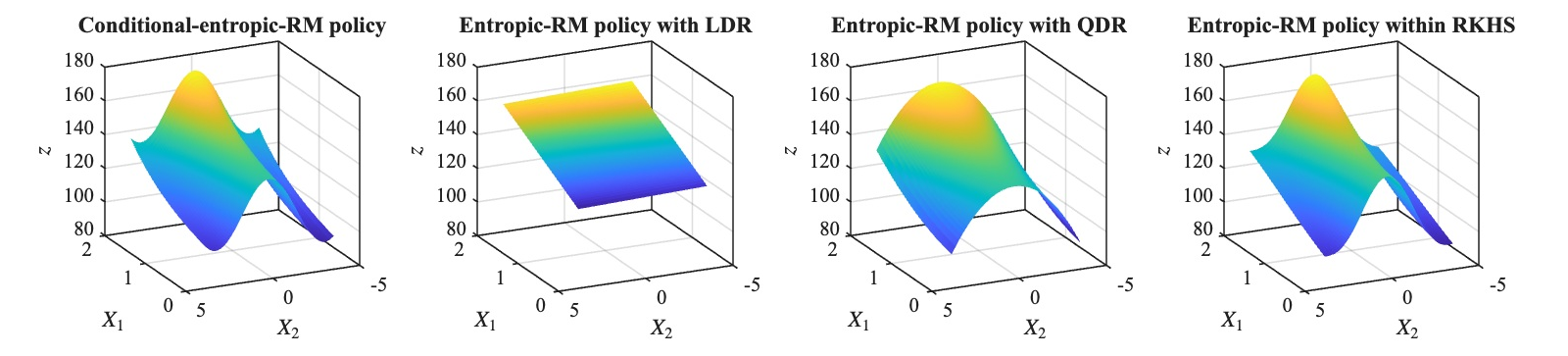}
    \caption{Optimal policy and SAA policies (sample size $N=800$) for risk-averse contextual optimization problem with entropic risk measure.}
    \label{fig:ex1-nonlinear-policy-ent}
\end{figure}

\begin{figure} [ht]
    \centering
    \includegraphics[width=0.95\linewidth]{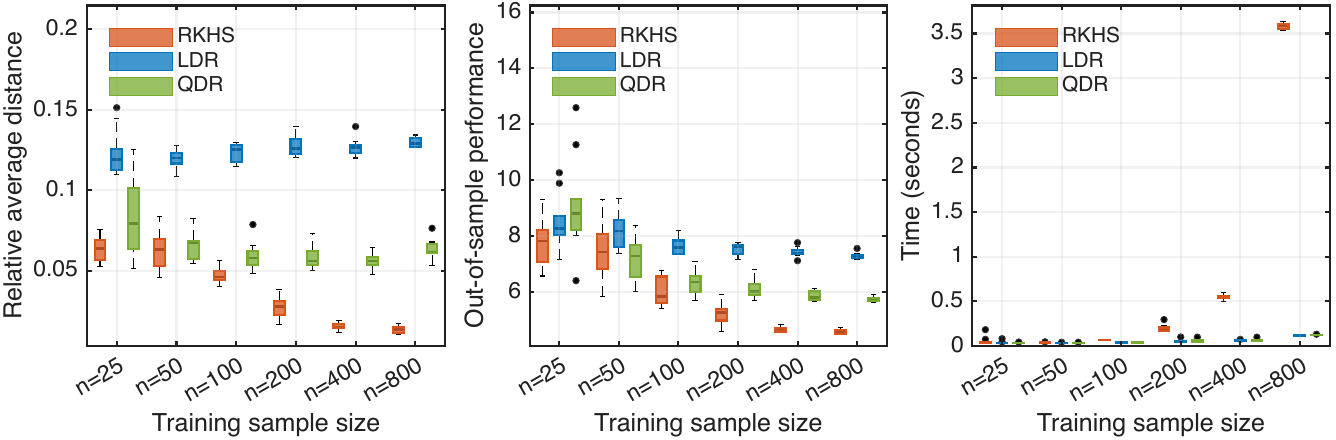}
    \caption{Performance of risk-averse contextual optimization problem with entropic risk measure within various hypothesis spaces.}
    \label{fig:ex1-nonlinear-boxplot-ent}
\end{figure}

For entropic risk measure,
the ex-ante-entropic-RM policy and the expected-entropic-RM policy are identical as discussed in Example \ref{example:expected-ent} and demonstrated in Section \ref{sec:experiment-newsvend-linear}.
Therefore,
we focus only on the ex-ante-entropic-RM policy and refer to it simply as the entropic-RM policy.
We vary the number of data points over $N\in \{25,50,100,200,400,800\}$ and conduct
20 replicate experiments for each sample size,
with the same parameter settings as in the above experiment for the expected-CVaR policy.
Fig.~\ref{fig:ex1-nonlinear-policy-ent} shows
the conditional-entropic-RM policy and the entropic-RM policies under three different hypothesis spaces: RKHS, LDR, and QDR.
Fig.~\ref{fig:ex1-nonlinear-boxplot-ent}
presents boxplots of their performance.
Similar to the results for the expected-CVaR policy, the entropic-RM policy within RKHS rapidly converges to the conditional-entropic-RM policy as
the sample size increases, although this comes with longer running time.

{
\subsubsection{Scalability of the SAA approach in RKHS.}

In the preceding experiments, we demonstrate that the expected-CVaR policy can achieve conditional optimality in a linear setting,
and we further present evidence that the RKHS methods are highly effective for learning a policy in a simple nonlinear setting.
In this part, we investigate the scalability of the proposed expected-CVaR policy within RKHS by varying the dimension of the context $X$.

Specifically, for each \(d\), we let \(X_1, X_2\) be \(d\)-dimensional random vectors whose components are i.i.d. draws from a normal distribution \(N(1, 0.5^2)\), truncated to \(\pm 2\) standard deviations from the mean.
Let \(B \in \mathbb{R}^d\) denote the weight vector defined by
\(
    B = (1,-1,2,-2,3,-3,\dots)^T,
\)
truncated to dimension \(d\).
Let the demand quantity \(Y\) depend nonlinearly on the contextual variables through
\[
Y = \max \left\{
5 \exp\!\left(2 + \frac{1}{2} B^T X_1\right)
+ 20 \cos\!\left(2 B^T X_2\right)
+ 100
+ \left| B^T X_1 + B^T X_2 + 12 \right| \epsilon,\; 0
\right\},
\]
where \(\epsilon\) follows a standard normal distribution \(N(0,1)\) truncated to \(\pm 2\) standard deviations from the mean.
We vary the dimension \(d\in \{2,4,6,8,10,12\}\) and
the training sample size $N\in \{100,200,400,800\}$, and conduct
20 replications for each \(d,N\) combination.

We compare three approaches: the ex-ante-CVaR policy derived from problem \eqref{eqn:newsvendor-ex-ante-cvar} in RKHS,
the expected-CVaR policy from problem \eqref{eqn:newsvendor-expected-cvar} in RKHS,
and a sequential learning and optimization (SLO) approach, in which the DM first uses the kernel regression to estimate the conditional distribution $P_{Y \mid X = x}$ for any new $x \in \mathcal{X}$, and then makes a decision based on this estimated conditional distribution.
For these methods, the kernel bandwidth is chosen as \(h_{N,d} = h_{\text{ref}} \left(\frac{10}{N} \right)^{1/d+4}\), which corresponds to the classical result found, e.g., in \cite{hall2004cross,stone1982optimal}.
For the ex-ante-CVaR policy and expected-CVaR policy, the regularization parameter for the decision policy $z(\cdot)$ is set to \(\lambda_N^z = \lambda_{\text{ref}}^z \left(\frac{10}{N} \right)^{1/3}\);
and for expected-CVaR policy,
the regularization parameter for auxiliary function $t(\cdot)$ is set to \(\lambda_N^t = \lambda_{\text{ref}}^t \left(\frac{10}{N} \right)^{1/3}\), which ensures that the regularization strength decreases with the sample size, as discussed in Section~\ref{sec:consist-saa-rkhs}.
For these parameters, the values of $h_{\text{ref}}, \lambda_{\text{ref}}^z,$ and $\lambda_{\text{ref}}^t$ are selected from a set of candidate values for each $d$.
We summarize the parameters in Table~\ref{tab:exp1-hyperparameters}.

\begin{table}[tb]
\centering
\small
\caption{Bandwidth parameters $h_{\text{ref}}$ and regularization parameters $\lambda_{\text{ref}}^{z},\lambda_{\text{ref}}^{t}$ selected by grid search for the scalability experiment. }
\label{tab:exp1-hyperparameters}
\renewcommand{\arraystretch}{1.15}
\begin{tabular}{c|cc|ccc|c}
\hline
& \multicolumn{2}{c|}{Ex ante CVaR over RKHS}
& \multicolumn{3}{c|}{Expected CVaR over RKHS}
& \multicolumn{1}{c}{SLO (kernel regression)} \\
\cline{2-3}\cline{4-6}\cline{7-7}
$d$
& $h_{\text{ref}}$ & $\lambda_{\text{ref}}^{z}$
& $h_{\text{ref}}$ & $\lambda_{\text{ref}}^{z}$ & $\lambda_{\text{ref}}^{t}$
& $h_{\text{ref}}$ \\
\hline
2  & 10  & $5\times 10^{-7}$ & 10  & $5\times 10^{-7}$ & $0.05$ & 0.3 \\
4  & 10 & $5\times 10^{-7}$ & 10 & $5\times 10^{-7}$ & $0.01$ & 0.3 \\
6  & 7  & $5\times 10^{-7}$ & 7  & $5\times 10^{-7}$ & $0.1$ & 0.3 \\
8  & 4  & $5\times 10^{-7}$ & 7  & $5\times 10^{-7}$ & $0.5$ & 0.3 \\
10 & 3  & $5\times 10^{-7}$ & 7  & $1\times 10^{-7}$ & $0.5$ & 0.3 \\
12 & 2  & $1\times 10^{-7}$ & 4  & $1\times 10^{-7}$ & $0.5$ & 0.3 \\
\hline
\end{tabular}
\end{table}

We evaluate the methods using the following metrics: the relative average distance to the conditional-CVaR policy, the out-of-sample performance, and wall clock time.
The testing sample is generated by a nested Monte Carlo simulation.
We first generate $1,000$ contexts, and then, for each context, we generate $400$ realizations of the problem data.
This design allows us to calculate the conditionally optimal policy on the same context sample and to use it as a benchmark for comparison.
The running time includes both the training time and the decision time on the testing dataset.
This is because the decision rule approaches (ex-ante-CVaR policy and expected-CVaR policy) spend most of their computational effort in the training stage, whereas, after the decision rule $z(\cdot)$ has been learned, the decision for a new context $x$ can be obtained directly from the decision rule as $z(x)$ and is therefore very fast.
By contrast, the SLO method does not require a training stage, and its computational cost is mainly spent on the decision stage, where the conditional distribution of $Y$ must be estimated and the corresponding decision-making problem must be solved for each newly observed context.
The numerical results are shown in Figures~\ref{fig:experiment1_boxplots_d2}-\ref{fig:experiment1_boxplots_d12}.

\begin{figure}[tbp]
    \centering
    \includegraphics[width=0.95\linewidth]{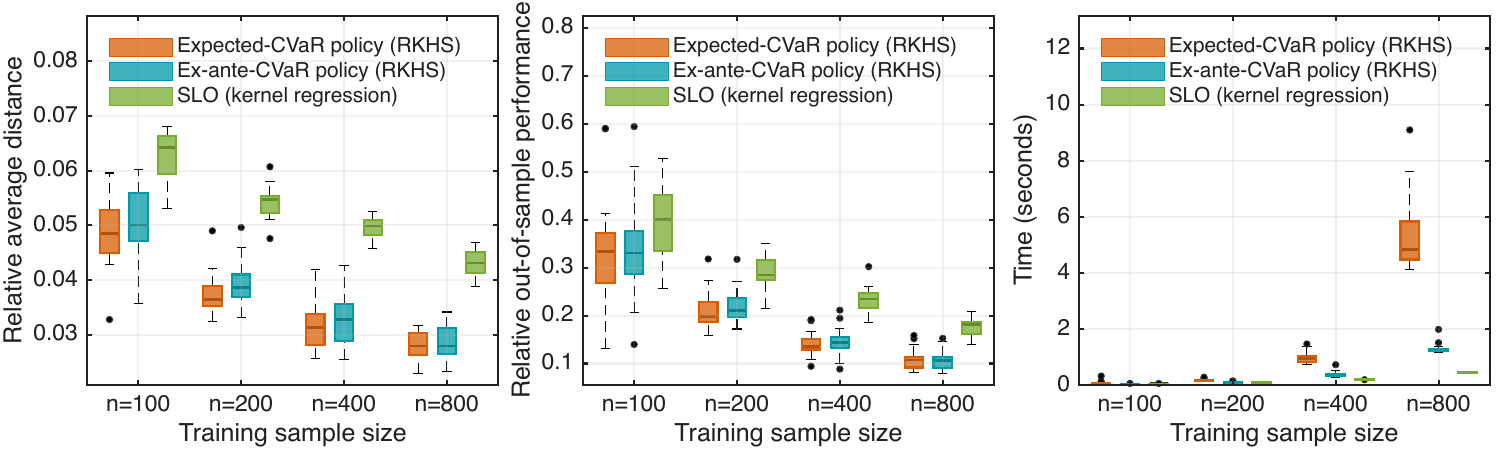}
    \caption{Performance across different training sample sizes when $d = 2$.}
    \label{fig:experiment1_boxplots_d2}
\end{figure}

\begin{figure}[tbp]
    \centering
    \includegraphics[width=0.95\linewidth]{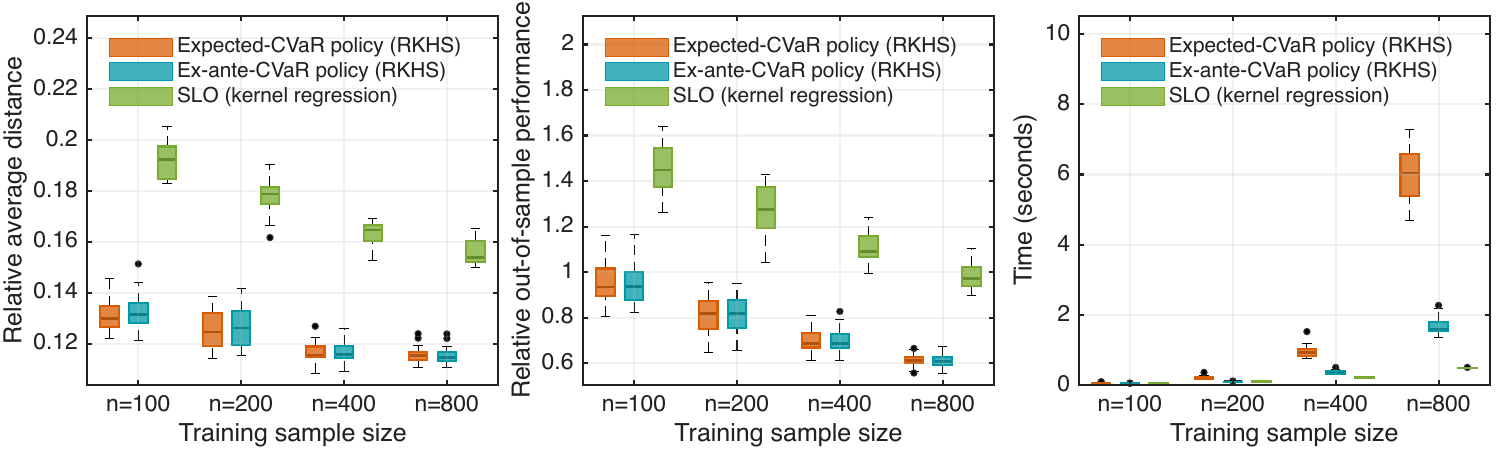}
    \caption{Performance across different training sample sizes when $d = 4$.}
    \label{fig:experiment1_boxplots_d4}
\end{figure}

\begin{figure}[tbp]
    \centering
    \includegraphics[width=0.95\linewidth]{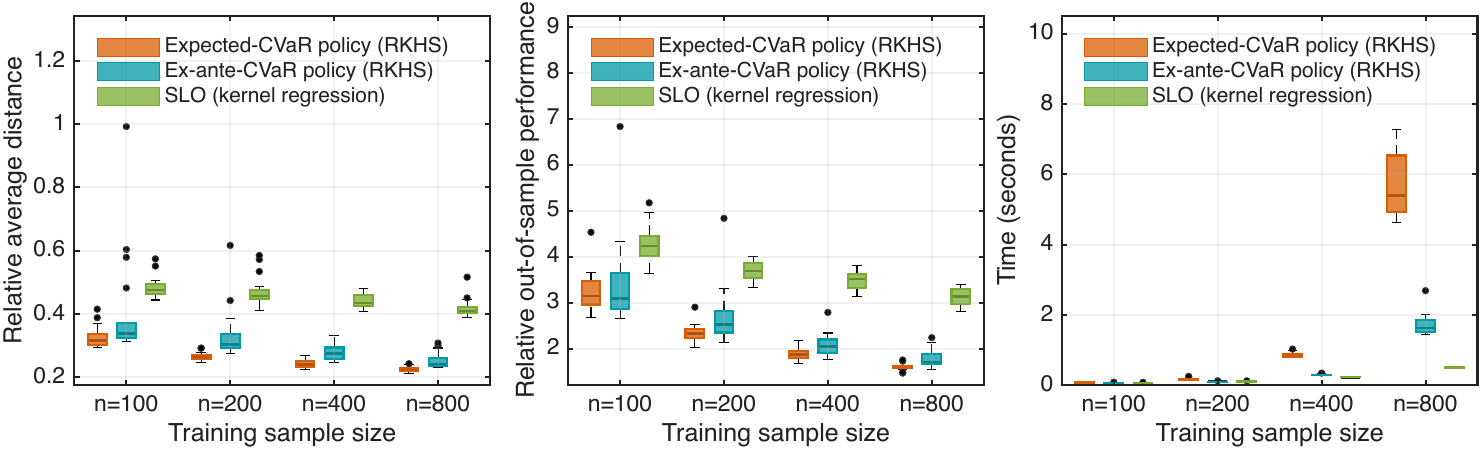}
    \caption{Performance across different training sample sizes when $d = 6$.}
    \label{fig:experiment1_boxplots_d6}
\end{figure}

\begin{figure}[tbp]
    \centering
    \includegraphics[width=0.95\linewidth]{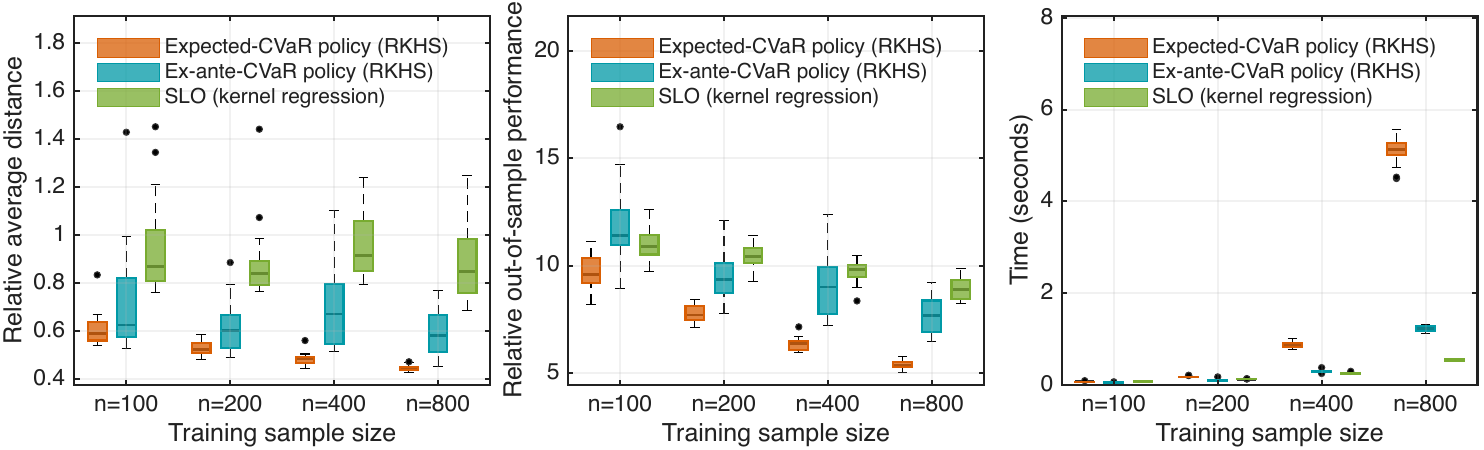}
    \caption{Performance across different training sample sizes when $d = 8$.}
    \label{fig:experiment1_boxplots_d8}
\end{figure}

\begin{figure}[tbp]
    \centering
    \includegraphics[width=0.95\linewidth]{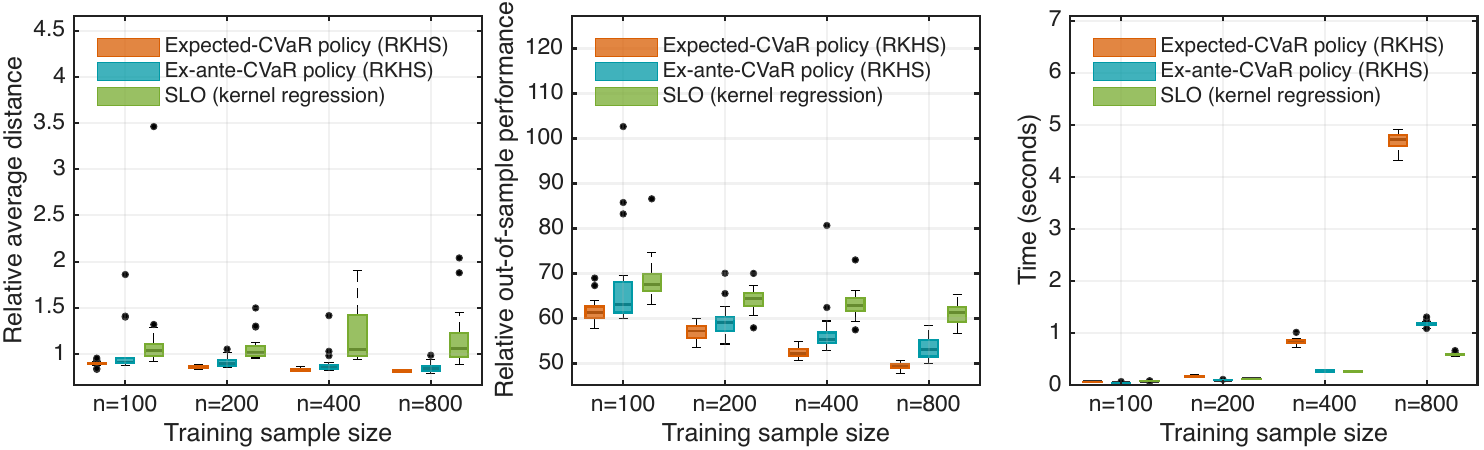}
    \caption{Performance across different training sample sizes when $d = 10$.}
    \label{fig:experiment1_boxplots_d10}
\end{figure}

\begin{figure}[tbp]
    \centering
    \includegraphics[width=0.95\linewidth]{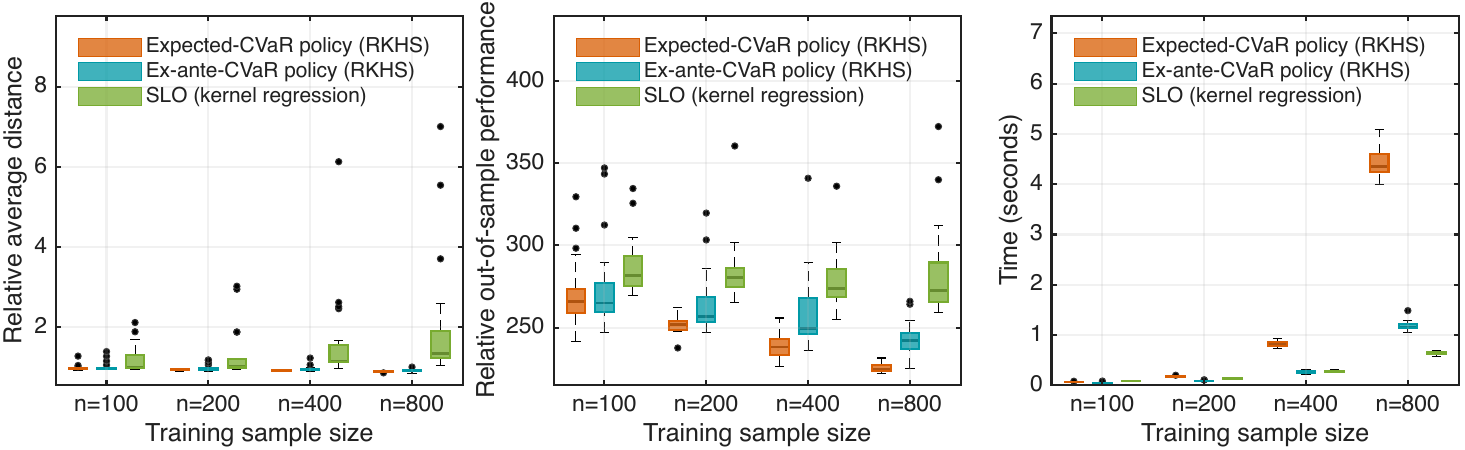}
    \caption{Performance across different training sample sizes when $d = 12$.}
    \label{fig:experiment1_boxplots_d12}
\end{figure}

The experiment results indicate that, for any fixed contextual dimension $d$, all methods achieve better decision quality as the sample size increases, at the cost of increased running time.
However, their performance declines as the dimension of the contextual information grows.
This phenomenon is consistent with the classical curse of dimensionality in nonparametric estimation:
in high-dimensional feature spaces, data become increasingly sparse, which reduces the accuracy of kernel smoothing and conditional distribution estimation \cite{stone1982optimal}.
Actually, improving RKHS-based methods in high-dimensional and large-scale dataset is an active research area in itself, with representative approaches including dimension reduction \cite{fukumizu2004dimensionality,fukumizu2009kernel}, Nystr\"om approximations \cite{gittens2013revisiting,williams2000using}, random features \cite{bach2017equivalence,rahimi2007random}, and sparse kernel methods \cite{smola2000sparse,tipping2001sparse}.
Since this topic is not the main focus of this work, we do not pursue it further and instead refer interested readers to the above literature.

Comparing the three methods, we find that the expected-CVaR policy within RKHS delivers the best overall decision performance, followed by the ex-ante-CVaR policy.
This finding is consistent with the theoretical analysis of decision rule approaches that they incorporate the prediction and decision-making into a unified optimization problem for policy learning, and the expected-CVaR policy can attain the conditional optimality, as discussed in Section \ref{sec:contextual-ra-opt}.
Nevertheless, both the expected-CVaR policy and the ex-ante-CVaR policy become computationally more demanding as the sample size increases, while the SLO approach is significantly faster.

}

\subsection{Portfolio selection}
\label{sec:portfolio}

In this experiment, we consider a portfolio selection problem,
which has been revisited in the setting of contextual optimization \cite{elmachtoub2022smart,nguyen2024robustifying}.
To maximize the expected return and reduce the risk,
Rockafellar and Uryasev \cite{rockafellar2000optimization} propose the mean-CVaR methodology where
the objective function is
\begin{eqnarray}
\label{eqn:mean-cvar}
    \min_{z\in \R^{d_z}_+}  \inmat{CVaR}_{\beta,P_Y}(-Y^Tz) - \eta \mathbb{E}_{P_Y}[Y^Tz],\qquad
    \inmat{s.t.}  \sum_{j=1}^{d_z} z^{(j)}\leq 1,
\end{eqnarray}
where $z\in \R^{d_z}$ denotes the portfolio decision, $Y$ denotes the assets’ future returns, $\eta\geq 0$ is a positive parameter that weights the preference between the portfolio return and the associated risk, and $\beta=0.9$ is fixed and used to form the objective function in subsequent methods.
As discussed in Section \ref{sec:intro}, there are contextual covariates $X$ that can be used when making the portfolio decision.
In this experiment, we set the dimension of covariates $d_x=5$, the number of assets $d_z = 50$ with the risk-free rate $y^{(51)}=0$.
We compare the following methods {for} the portfolio selection problem:
\begin{enumerate}
    \item[(i)] {The equal-weighted model (EW)} where the portfolio decision coincides with the {$1/d_z$}-portfolio \cite{demiguel2009optimal}.

    \item[(ii)] The mean-CVaR model without contextual information (MC):  the portfolio decision is obtained by solving problem (\ref{eqn:mean-cvar}).

    \item[(iii)] The contextual mean-ex-ante-CVaR model (CMEAC):
    the portfolio decision is obtained by solving
        \begin{eqnarray*}
        \min_{{z(\cdot)\in\mathcal{H}^{50}}, t\in\R} &&
        \overbrace{t + \frac{1}{(1-\beta)} \mathbb{E}_{P_{X,Y}} [(-Y^T z(X)-t)_+ ]}^{\textrm{CVaR}_{\beta,{P_{X,Y}}}(-Y^Tz(X))}- {\eta} \mathbb{E}_{P_{X,Y}} [Y^T z(X)] \\
        &&+ \lambda_z \|z(\cdot)\|_{\mathcal{H}^{50}}^2\\
        \inmat{s.t.} && \sum_{j=1}^{50} (z(x))^{(j)} \leq 1,\;
        z(x)\geq0,\; \forall x\in\mathcal{X}
    \end{eqnarray*}
    where $\mathcal{H}^{50}$ is the Cartesian product of $50$ copies of RKHS $\mathcal{H}$ generated by Gaussian kernel.

    \item[(iv)] The contextual mean-expected-CVaR model (CMEC):
    the portfolio decision is obtained by solving
    \begin{eqnarray}
        \min_{z(\cdot)\in\mathcal{H}^{50}, t\in\mathcal{H}} &&
        \overbrace{\mathbb{E}_{P_{X,Y}}[t(X) + \frac{1}{1-\beta}  (-Y^T z(X)-t(X))_+]}^{\mathbb{E}_{P_X}[\textrm{CVaR}_{\beta,{P_{Y|X}}}(-Y^Tz(X))]}  - \eta\mathbb{E}_{P_{X,Y}}[Y^T z(X)]\nonumber\\
        &&+ \lambda_z \|z(\cdot)\|_{\mathcal{H}^{50}}^2 + \lambda_t \|t(\cdot)\|^2, \label{eqm:CMEC-portfolio}\\
        \inmat{s.t.} && \sum_{j=1}^{50} (z(x))^{(j)} \leq 1,\;
        z(x)\geq0,\; \forall x\in\mathcal{X}.\nonumber
    \end{eqnarray}

\end{enumerate}

{
These problems satisfy
the Assumption~\ref{assumption:well-defined_rho},
Assumption~\ref{assumption:convex_finite-minimum},
and Assumption~\ref{assump:saa-convergence-unbounded}
required for the theoretical developments across Sections \ref{sec:contextual-risk-measure}-\ref{sec:DR-rkhs}.
We verify some of them in
Section \ref{sec:veri-assum-portfolio}.
}
In (iii) and (iv), we restrict the policy $z(\cdot)$ to the RKHS generated by Gaussian kernel.
For a new contextual observation $x_0$, we can find the portfolio decision as $z(x_0)$. If $z(x_0)$ is not feasible, then project it onto the feasible set by solving the following problem as discussed in \cite{bertsimas2022data}:
\begin{eqnarray}
    \min_{z\in \R^{50}} && \|z-z(x_0)\|^2 \nonumber\\
    \inmat{s.t.} && \sum_{j=1}^{50} (z(x_i))^{(j)} \leq 1,\;
        z(x_i)\geq0,\; \forall i=1,\dots,N.
\end{eqnarray}

The data are generated according to a similar process as in  \cite{elmachtoub2022smart}. We first generate a matrix $B\in \mathcal{R}^{d_y\times d_x}$ that encodes the problem data of the true model, whereby each entry of $B$ is a Bernoulli random variable that is equal to $1$ with probability 0.5.
We generate a factor loading matrix $L\in \R^{d_y\times 4}$ where each entry of $L$ is uniformly distributed on $[-0.0025\tau,0.0025\tau]$ and $\tau$ is the noise-level parameter.
Then, the training/testing data are generated according to the following process:
\begin{enumerate}
    \item[(i)] The covariates $\{x_i\}_{i=1}^N$ are i.i.d. generated from a multivariate Gaussian distribution with i.i.d. standard normal entries, i.e., $X\sim N(0_{d_x\times 1}, I_{d_x\times d_x})$.
    \item[(ii)] Given the realization of covariates $x_i$, the return is set to
    $$y_i = \left( \frac{0.05}{\sqrt{d_x}}Bx_i + 0.1^{1/p} \right)^p + L\epsilon_{1,i} + 0.01\tau \epsilon_{2,i},$$
    where $\epsilon_{1,i}$ follows $N(0,I_{4\times 4})$, $\epsilon_{2,i}$ follows $N(0,I_{d_y\times d_y})$
    and $p$ is a positive integer parameter to encode the linear/nonlinear relationship between $X$ and $Y$.
    In this way,  we can generate the return such that the mean of the assets' returns is $\bar{y}_i:= \left( \frac{0.05}{\sqrt{d_x}}Bx_i + 0.1^{1/p} \right)^p$ and the covariance of the returns is $\Sigma = LL^T + (0.01\tau)^2 I_{d_y\times d_y}$.
    As discussed in \cite{elmachtoub2022smart}, $y_i\in [0,1]$ with high probability.
\end{enumerate}

We vary the training set size $N\in \{100,500\}$, the parameter of risk preference $\eta \in \{1/3, 1,3\}$,
the noise level $\tau\in \{1,2\}$,
the linear/nonlinear parameter $p\in\{1,4\}$,
and fix the testing set size $N_{\text{test}} = 50000$.
For each combination of parameter $(\eta,\tau,p)$, we run 50 independent trials and report the performance metrics on average value.
Denote the empirical distribution constructed from testing set as $P_{\text{test}}$.
The performance metrics are chosen as
the expected return ($\mathbb{E}$), i.e., $\mathbb{E}_{P_{\text{test}}}[Y^Tz(X)]$,
the CVaR value of negated return (CVaR), i.e., $\inmat{CVaR}_{\beta,P_{\text{test}}}(-Y^T z(X))$,
the tradeoff between expected return and CVaR
($\eta\mathbb{E}-\inmat{CVaR}$), i.e.,
\begin{eqnarray*}
    \eta \mathbb{E}_{P_{\text{test}}}[Y^Tz(X)] - \inmat{CVaR}_{\beta,P_{\text{test}}}(-Y^T z(X)),
\end{eqnarray*}
and the relative regret with respect to the largest return (rela. regret), i.e.,
\begin{eqnarray}
    \frac{1}{N_{\text{test}}}\sum_{i=1}^{N_{\text{test}}} \frac{\max\{0,\max\limits_{j=1,\dots,50}y_i^{(j)}\}-y_i^Tz(x_i)}{\max\{0,\max\limits_{j=1,\dots,50}y_i^{(j)}\}}.
\end{eqnarray}

Tables \ref{tab:experiment2-cvar-eta1/3}-\ref{tab:experiment2-cvar-eta3-N500}
presents the average out-of-sample performance.
Overall, the contextual mean-expected-CVaR model (CMEC) demonstrates a superior performance compared to the other methods in most cases.
It achieves the highest expected return, the smallest relative regret, and the highest $\eta\mathbb{E}-\inmat{CVaR}$ value.
When the decision maker is strongly risk-averse ($\eta=1/3$) or the data noise level is high ($\tau=2$),
the CMEC model yields the lowest CVaR, indicating its robustness in risk-averse environments.
These observations suggest that the CMEC model has great out-of-sample performance and
effectively minimizes the risk while achieving a reasonable balance between revenue and risk.

\begin{table} [hbt]
    \centering
\begin{tabular}{llccccc}
\hline$(\eta,\tau,p)$ & Model & $\mathbb{E}$ $\uparrow$ & CVaR $\downarrow$ & $\eta\mathbb{E} - \inmat{CVaR}$  $\uparrow$ & rela. regret $\downarrow$ \\
\hline
$(1/3,1,1)$
& EW & 0.1000 & -0.0604 & 0.0938 & 0.2924  \\
& MC & 0.1000 & -0.0839 & 0.1172 & 0.2708 \\
& CMEAC & 0.1256 & -0.0872 & 0.1291 & 0.1093  \\
& CMEC & \textbf{0.1276} & \textbf{-0.0874} & \textbf{0.1299} & \textbf{0.0971} \\
$(1/3,1,4)$
& EW & 0.1018 & -0.0728 & 0.1067 & 0.2380  \\
& MC & 0.1002 & -0.0838 & 0.1173 & 0.2354 \\
& CMEAC & 0.1210 & -0.0885 & 0.1288 & 0.0966  \\
& CMEC & \textbf{0.1226} & \textbf{-0.0888} & \textbf{0.1297} & \textbf{0.0855} \\
$(1/3,2,1)$
& EW & 0.1000 & -0.0523 & 0.0857 & 0.3065  \\
& MC & 0.1000 & -0.0678 & 0.1011 & 0.2864 \\
& CMEAC & 0.1248 & -0.0767 & 0.1183 & 0.1292  \\
& CMEC & \textbf{0.1271} & \textbf{-0.0781} & \textbf{0.1204} & \textbf{0.1136} \\
$(1/3,2,4)$
& EW & 0.1017 & -0.0623 & 0.0962 & 0.2538  \\
& MC & 0.1005 & -0.0678 & 0.1013 & 0.2516 \\
& CMEAC & 0.1203 & -0.0769 & 0.1170 & 0.1184  \\
& CMEC & \textbf{0.1223} & \textbf{-0.0782} & \textbf{0.1189} & \textbf{0.1040} \\
\hline
\end{tabular}
    \caption{Average out-of-sample performance with  $N=100,\eta=1/3$.}
    \label{tab:experiment2-cvar-eta1/3}
\end{table}
\vspace{-0.3cm}

\begin{table}[hbt]
    \centering
\begin{tabular}{llccccc}
\hline$(\eta,\tau,p)$ & Model & $\mathbb{E}$ $\uparrow$ & CVaR $\downarrow$ & $\eta\mathbb{E} - \inmat{CVaR}$  $\uparrow$ & rela. regret $\downarrow$ \\
\hline
$(1,1,1)$
& EW & 0.1000 & -0.0604 & 0.1604 & 0.2924  \\
& MC & 0.1000 & -0.0837 & 0.1837 & 0.2711 \\
& CMEAC & 0.1268 & \textbf{-0.0873} & 0.2141 & 0.1016  \\
& CMEC & \textbf{0.1275} & -0.0873 & \textbf{0.2148} & \textbf{0.0968} \\
$(1,1,4)$
& EW & 0.1018 & -0.0728 & 0.1746 & 0.2380  \\
& MC & 0.1003 & -0.0834 & 0.1837 & 0.2356 \\
& CMEAC & 0.1214 & \textbf{-0.0886} & 0.2100 & 0.0942  \\
& CMEC & \textbf{0.1220} & -0.0883 & \textbf{0.2103} & \textbf{0.0901} \\
$(1,2,1)$
& EW & 0.1000 & -0.0523 & 0.1524 & 0.3065  \\
& MC & 0.1000 & -0.0674 & 0.1674 & 0.2869 \\
& CMEAC & 0.1262 & -0.0773 & 0.2035 & 0.1202  \\
& CMEC & \textbf{0.1271} & \textbf{-0.0779} & \textbf{0.2050} & \textbf{0.1138} \\
$(1,2,4)$
& EW & 0.1017 & -0.0623 & 0.1640 & 0.2538  \\
& MC & 0.1006 & -0.0670 & 0.1676 & 0.2519 \\
& CMEAC & 0.1216 & -0.0776 & 0.1992 & 0.1093  \\
& CMEC & \textbf{0.1224} & \textbf{-0.0782} & \textbf{0.2006} & \textbf{0.1035} \\
\hline
\end{tabular}
    \caption{Average out-of-sample performance  with  $N=100,\eta=1$.}
    \label{tab:experiment2-cvar-eta1}
\end{table}
\vspace{-0.3cm}

\begin{table}[hbt]
    \centering
\begin{tabular}{llccccc}
\hline$(\eta,\tau,p)$ & Model & $\mathbb{E}$ $\uparrow$ & CVaR $\downarrow$ & $\eta\mathbb{E} - \inmat{CVaR}$  $\uparrow$ & rela. regret $\downarrow$ \\
\hline
$(3,1,1)$
& EW & 0.1000 & -0.0604 & 0.3602 & 0.2924  \\
& MC & 0.1000 & -0.0819 & 0.3819 & 0.2727 \\
& CMEAC & 0.1270 & \textbf{-0.0869} & 0.4678 & 0.1007  \\
& CMEC & \textbf{0.1272} & -0.0868 & \textbf{0.4685} & \textbf{0.0989} \\
$(3,1,4)$
& EW & 0.1018 & -0.0728 & 0.3781 & 0.2380  \\
& MC & 0.1008 & -0.0791 & 0.3814 & 0.2364 \\
& CMEAC & 0.1224 & \textbf{-0.0889} & 0.4562 & 0.0869  \\
& CMEC & \textbf{0.1226} & -0.0888 & \textbf{0.4567} & \textbf{0.0855} \\
$(3,2,1)$
& EW & 0.1000 & -0.0523 & 0.3524 & 0.3065  \\
& MC & 0.1000 & -0.0647 & 0.3646 & 0.2898 \\
& CMEAC & 0.1264 & -0.0772 & 0.4563 & 0.1191  \\
& CMEC & \textbf{0.1267} & \textbf{ -0.0774} & \textbf{0.4576} & \textbf{0.1166} \\
$(3,2,4)$
& EW & 0.1017 & -0.0623 & 0.3675 & 0.2538  \\
& MC & 0.1006 & -0.0625 & 0.3662 & 0.2529 \\
& CMEAC & 0.1219 & -0.0778 & 0.4434 & 0.1073  \\
& CMEC & \textbf{0.1221} & \textbf{-0.0780} & \textbf{0.4444} & \textbf{0.1053} \\
\hline
\end{tabular}
    \caption{Average out-of-sample performance  with $N=100,\eta=3$.}
    \label{tab:experiment2-cvar2-eta3}
\end{table}

\begin{table}[hbt]
    \centering
\begin{tabular}{llccccc}
\hline$(\eta,\tau,p)$ & Model & $\mathbb{E}$ $\uparrow$ & CVaR $\downarrow$ & $\eta\mathbb{E} - \inmat{CVaR}$  $\uparrow$ & rela. regret $\downarrow$ \\
\hline
$(1/3,1,1)$
& EW & 0.1000 & -0.0604 & 0.0937 & 0.2924  \\
& MC & 0.1000 & -0.0845 & 0.1178 & 0.2700 \\
& CMEAC & 0.1286 & \textbf{-0.0910} & \textbf{0.1339} & 0.0878  \\
& CMEC & \textbf{0.1307} & -0.0900 & 0.1336 & \textbf{0.0745} \\
$(1/3,1,4)$
& EW & 0.1017 & -0.0728 & 0.1067 & 0.2380  \\
& MC & 0.1002 & -0.0845 & 0.1179 & 0.2352 \\
& CMEAC & 0.1237 & \textbf{-0.0912} & 0.1324 & 0.0771  \\
& CMEC & \textbf{0.1254} & -0.0907 & \textbf{0.1325} & \textbf{0.0659} \\
$(1/3,2,1)$
& EW & 0.1000 & -0.0524 & 0.0857 & 0.3065  \\
& MC & 0.1000 & -0.0690 & 0.1023 & 0.2850 \\
& CMEAC & 0.1281 & -0.0805 & 0.1232 & 0.1055  \\
& CMEC & \textbf{0.1305} & \textbf{-0.0808} & \textbf{0.1243} & \textbf{0.0901} \\
$(1/3,2,4)$
& EW & 0.1017 & -0.0623 & 0.0962 & 0.2538  \\
& MC & 0.1004 & -0.0691 & 0.1025 & 0.2512 \\
& CMEAC & 0.1231 & -0.0796 & 0.1207 & 0.0973  \\
& CMEC & \textbf{0.1252} & \textbf{-0.0802} & \textbf{0.1219} & \textbf{0.0835} \\
\hline
\end{tabular}
    \caption{Average out-of-sample performance with  $N=500,\eta=1/3$.}
    \label{tab:experiment2-cvar-eta1/3-N500}
\end{table}

\begin{table}[hbt]
    \centering
\begin{tabular}{llccccc}
\hline$(\eta,\tau,p)$ & Model & $\mathbb{E}$ $\uparrow$ & CVaR $\downarrow$ & $\eta\mathbb{E} - \inmat{CVaR}$  $\uparrow$ & rela. regret $\downarrow$ \\
\hline
$(1,1,1)$
& EW & 0.1000 & -0.0604 & 0.1604 & 0.2924  \\
& MC & 0.1000 & -0.0841 & 0.1844 & 0.2701 \\
& CMEAC & 0.1299 & \textbf{-0.0911} & \textbf{0.2209} & 0.0794  \\
& CMEC & \textbf{0.1307} & -0.0900 & 0.2207 & \textbf{0.0747} \\
$(1,1,4)$
& EW & 0.1017 & -0.0728 & 0.1746 & 0.2380  \\
& MC & 0.1002 & -0.0844 & 0.1846 & 0.2352 \\
& CMEAC & 0.1247 & \textbf{-0.0913} & 0.2160 & 0.0699  \\
& CMEC & \textbf{0.1254} & -0.0907 & \textbf{0.2160} & \textbf{0.0660} \\
$(1,2,1)$
& EW & 0.1000 & -0.0524 & 0.1524 & 0.3065  \\
& MC & 0.1000 & -0.0689 & 0.1689 & 0.2851 \\
& CMEAC & 0.1296 & \textbf{-0.0810} & 0.2106 & 0.0957  \\
& CMEC & \textbf{0.1305} & -0.0807 & \textbf{0.2112} & \textbf{0.0903} \\
$(1,2,4)$
& EW & 0.1017 & -0.0623 & 0.1640 & 0.2538  \\
& MC & 0.1004 & -0.0688 & 0.1693 & 0.2514 \\
& CMEAC & 0.1244 & -0.0801 & 0.2046 & 0.0885  \\
& CMEC & \textbf{0.1252} & \textbf{-0.0802} & \textbf{0.2053} & \textbf{0.0836} \\
\hline
\end{tabular}
    \caption{Average out-of-sample performance with  $N=500,\eta=1$.}
    \label{tab:experiment2-cvar-eta1-N500}
\end{table}

\begin{table}[hbt]
    \centering
\begin{tabular}{llccccc}
\hline$(\eta,\tau,p)$ & Model & $\mathbb{E}$ $\uparrow$ & CVaR $\downarrow$ & $\eta\mathbb{E} - \inmat{CVaR}$  $\uparrow$ & rela. regret $\downarrow$ \\
\hline
$(3,1,1)$
& EW & 0.1000 & -0.0604 & 0.3604 & 0.2924  \\
& MC & 0.1000 & -0.0842 & 0.3842 & 0.2704 \\
& CMEAC & 0.1304 & \textbf{-0.0907} & \textbf{0.4819} & 0.0763  \\
& CMEC & \textbf{0.1306} & -0.0899 & 0.4818 & \textbf{0.0751} \\
$(3,1,4)$
& EW & 0.1017 & -0.0728 & 0.3780 & 0.2380  \\
& MC & 0.1004 & -0.0835 & 0.3847 & 0.2355 \\
& CMEAC & 0.1251 & \textbf{-0.0912} & 0.4665 & 0.0672  \\
& CMEC & \textbf{0.1253} & -0.0907 & \textbf{0.4666} & \textbf{0.0664} \\
$(3,2,1)$
& EW & 0.1000 & -0.0524 & 0.3524 & 0.3065  \\
& MC & 0.1000 & -0.0685 & 0.3685 & 0.2858 \\
& CMEAC & 0.1302 & \textbf{-0.0810} & 0.4715 & 0.0921  \\
& CMEC & \textbf{0.1304} & -0.0807 & \textbf{0.4719} & \textbf{0.0910} \\
$(3,2,4)$
& EW & 0.1017 & -0.0623 & 0.3675 & 0.2538  \\
& MC & 0.1008 & -0.0671 & 0.3695 & 0.2520 \\
& CMEAC & 0.1249 & \textbf{-0.0803} & 0.4549 & 0.0853  \\
& CMEC & \textbf{0.1251} & -0.0801 & \textbf{0.4553} & \textbf{0.0843} \\
\hline
\end{tabular}
    \caption{Average out-of-sample performance with  $N=500,\eta=3$.}
    \label{tab:experiment2-cvar-eta3-N500}
\end{table}

\clearpage

\section{Potential extensions}
\label{sec:extension}

We present some potential extensions of the basic methodology presented in earlier sections.

\subsection{Risk-averse regret-based contextual optimization}
\label{sec:regret-based}

Elmachtoub and Grigas \cite{elmachtoub2022smart} propose
to tackle the contextual optimization from a perspective of regret.
In the risk-averse contextual optimization setting, we provide new insight into the regret-based contextual optimization.
Regret
represents
the additional loss of the decision $z$ compared to the optimal decision $z^*(y)$
under the true problem data $y$, i.e.,
$
    R(z,y):= c(z,y) -  c(z^*(y),y),
$
where
$z^*(y) :=\arg\limits\min_{z\in Z} c(z,y)$.
The risk-averse regret-based contextual optimization problem can be defined as
\begin{eqnarray} \label{eqn:regret-ep}
    &&\min_{g\in\widetilde{\mathcal{G}}} \rho_{P_X}^{(1)}\left(\rho_{P_{Y|X}}^{(2)}\big(R(g(X),Y)\big)\right)\nonumber\\
    &=& \min_{g\in\widetilde{\mathcal{G}}}\rho_{P_X}^{(1)}\left(\rho_{P_{Y|X}}^{(2)}\big(c(g(X),Y) -  c(z^*(Y),Y)\big)\right).
\end{eqnarray}
When $\rho^{(1)},\rho^{(2)}$ are chosen as entropic risk measures $\rho^{\text{ent}}_\gamma$,
model (\ref{eqn:regret-ep})
is
\begin{eqnarray} \label{eqn:regret-based-context-risk-measure}
    \min_{g\in\widetilde{\mathcal{G}}} \rho^{\text{ent}}_{\gamma,{P_X}} \left( \rho^{\text{ent}}_{\gamma,P_{Y|X}} ( R(g(x),y) ) \right)
    =\min_{g\in\widetilde{\mathcal{G}}} \frac{1}{\gamma} \ln  \mathbb{E}_{P} \left[\frac{1}{e^{\gamma c(z^*(Y),Y)}}e^{\gamma c(g(X),Y)}\right].
\end{eqnarray}
Compared
to model (\ref{eqn:example4.4ent}),
the term $\frac{1}{e^{\gamma c(z^*(y),y)}}$ acts as a weight for the negative exponential utility $e^{\gamma c(g(x),y)}$.
An intuitive explanation
for problem (\ref{eqn:regret-based-context-risk-measure})
is that
if
the optimal cost $c(z^*(y),y)$
obtained under the true problem data $y$ is large,
then it is acceptable for a proportionally large cost $c(g(x),y)$ incurred by policy $g$.
In this case, the negative exponential utility $e^{\gamma c(g(x),y)}$ of policy $g$ under $y$ should be assigned less weight when evaluating performance.
In other words,
the negative exponential utility of the policy is reweighed by the optimal negative exponential utility
derived by an informed decision maker.

\subsection{Risk-averse integrated learning and optimization}
\label{sec:regret-based-ILO}

The proposed risk-averse framework
can be used in integrated learning and optimization (ILO) model:
\begin{eqnarray} \label{eqn:ilo-ra}
    \min_{f\in\mathcal{F}} &&
        \rho_{P_X}^{(1)}\left(\rho_{P_{Y|X}}^{(2)}(L(z^*(X,f),X,Y))\right)
        \nonumber\\
        \inmat{s.t.}&& z^*(x,f) \in \arg\min_{z\in \mathcal{Z}} \rho_{f(x)}^{(2)}(c(z,Y)),\; \forall x \in \mathcal{X}.
\end{eqnarray}
As discussed in Section \ref{sec:3.1}, if there is an ex ante risk mapping $\rho^{\text{EA}}_P = \rho_{P_X}^{(1)}\circ\rho_{P_{Y|X}}^{(2)}$, then
model (\ref{eqn:ilo-ra}) is equivalent to
\begin{eqnarray} \label{eqn:ilo-ra-single}
    \min_{f\in\mathcal{F}} &&
        \rho_{P}^{\text{EA}}  \left(L(z^*(X,f),X,Y)\right)\quad
        \nonumber\\
        \inmat{s.t.} &&  z^*(x,f) \in \arg\min_{z\in \mathcal{Z}} \rho_{f(x)}^{(2)}(c(z,Y)),\; \forall x \in \mathcal{X}.
\end{eqnarray}
When the joint distribution $P$ is unknown and approximated by empirical distribution $P^N(\cdot) = \frac{1}{N} \sum_{i=1}^N \mathds{1}_{(x_i,y_i)}(\cdot)$, then the model (\ref{eqn:ilo-ra-single}) can be approximated
by the following SAA problem:
\begin{eqnarray*}
    \min_{f\in\mathcal{F}} &&
        \rho_{P^N}^{\text{EA}}\left(L(z^*(X,f),X,Y)\right)
        \nonumber\\
        \inmat{s.t.} &&  z^*(x_i,f) \in \arg\min_{z\in \mathcal{Z}} \rho_{f(x_i)}^{(2)}(c(z,y_i)),\; \forall i=1,\dots,N.
\end{eqnarray*}
Specifically,
the risk-averse ILO model can be extended to regret-based ILO \cite{elmachtoub2022smart} by letting $L(z,x,y)=R(z,y)=c(z,y)-c(z^*(y),y)$ as discussed in Section \ref{sec:regret-based}.

\subsection{Relation with expected utility models via certainty equivalents}

Bazier-Matte and Delage \cite{bazier2020generalization} investigate the contextual optimization from the perspective of expected utility theory.
As discussed in Example \ref{example:expected-ent}, when $\rho_{P_X}^{(1)},\rho_{P_{Y|X}}^{(2)}$ are chosen as entropic risk measures, then the risk-averse contextual optimization problem (\ref{eqn:example4.4ent}) is equivalent to the expected utility maximization problem with exponential utility $u(t) = 1-e^{-\gamma t}$.
This equivalence can be easily obtained using the tower property of certainty equivalent:
$
    \rho(\xi) := - u^{-1}(\mathbb{E}_P[u(-\xi)]),
$
where
$u: [\underline{\xi},\bar{\xi}] \rightarrow \R$ is a strictly increasing and continuous function, see Theorem 1.4 in \cite{kupper2009representation}.
Let $r(t) := - \frac{u''(t)}{u'(t)}$ be the Arrow–Pratt risk aversion measure, see \cite{pratt1978risk}.
To ensure the convexity of $\rho(\xi)$ and thus facilitate the optimization,
$1/r(t)$ is required to be concave.
As shown in \cite{ben2007old}, a class of risk-averse utility functions satisfying the requirement
is the class where $1/r(t)$ is a linear function, i.e., $\frac{1}{r(t)} = a t+ b $ with $at+b>0$.
These utility functions are in one of the following three forms:
{\begin{eqnarray*}
    u(t)=\left\{\begin{array}{lll}
b\left(1-e^{-t / b}\right) & \text { if } & a=0, b \neq 0, \\
\log (b+t) & \text { if } & a=1, \\
(a t+b)^{\frac{a-1}{a}} & \text { if } & a \neq 0, a \neq 1.
\end{array}\right.
\end{eqnarray*}}
The corresponding risk metrics are respectively defined as entropic risk measure, $\rho(\xi) = -e^{\mathbb{E}\left[\log(b-\xi)\right]} + b$ and $\rho(\xi) = -\frac{1}{a}\left(
\mathbb{E}\left[(-a\xi+b)^{\frac{a-1}{a}}\right] \right)^{\frac{a}{a-1}} + \frac{b}{a}$.
The last two risk metrics are defined for the case where the random loss is upper bounded by $\bar{\xi}\leq b/a$.
This kind of upper boundedness can be satisfied in some practical applications and is assumed as in \cite{bertsimas2022data}.

\section{Concluding remarks}
\label{sec:conclusion}

In this paper, we explore the risk management
in
contextual optimization problems
and do so
by evaluating the risk arising from PDU and CU with appropriate risk measures.
We establish conditions for choosing risk
measures which ensure both the computational tractability
and the optimality in the conditional risk minimization problem based on the conditional distribution of PDU.
We propose a set of risk measures/metrics for PDU
that satisfy the conditions,
including entropic risk measure, optimized certainty equivalents, mean-upper-semideviation, spectral risk measures, and distortion risk metrics.
Specifically, our analysis shows that the policies derived from contextual optimization problems with the ex ante CVaR or nested-CVaR objective do not achieve optimality in conditional risk minimization.
Alternatively, we propose to consider the contextual optimization with
the risk measures/metrics mentioned above, including expected-CVaR objective.

To solve
the proposed risk-averse contextual optimization problem
particularly in a data-driven environment,
we use the
sample average approximation method and demonstrate
under some moderate conditions
the exponential rate of convergence of the optimal value obtained from the SAA problem to its true counterpart as the sample size increases.
Particular attention is paid to the case when the hypothesis space is RKHS.
Finally, we
examine the performance of
the proposed models and computational schemes
by applying them to
a newsvendor problem and a
portfolio selection problem.

A crucial assumption in the convergence/consistency analysis is that the samples are i.i.d.. This assumption
may not be satisfied in real data-driven problems.
Moreover,
the training data
may be
potentially perturbed or even corrupted.
This raises a question as to
whether the optimal values and optimal
policies obtained from SAA models are statistically
stable.
This kind of analysis is conducted in \cite{guo2023statistical,zhang2024statistical} for
general utility preference robust
optimization problems and machine learning problems.
We may follow this stream of research to discuss statistical robustness for the contextual optimization problems and leave this for future research.

\section{Auxiliary results and proofs}

\subsection{Interchangeability principle}

We first
restate an interchangeability principle in
Theorem 2.2 in \cite{hiai1977integrals}.

\begin{lemma} [Theorem 2.2 in \cite{hiai1977integrals}]
\label{lemma-interchange-expectation}
 Consider a real separable Banach space $\mathbb{Z}$ with Borel field $\B(\mathbb{Z})$ and
a complete sample space $\Omega$ associated with
sigma algebra $\F$ and measure
$\mathbb{P}$.
 Let $\mathcal{Z}:\Omega \rightrightarrows \mathbb{Z}$ be a
 measurable set-valued mapping with images in the family of nonempty closed subsets of $\mathbb{Z}$.  Let $\mathcal{L}$ be
 a Banach
 space of measurable functions $ g: \Omega \rightarrow \mathbb{Z}$ with finite $p$-th order moments ($1\leq p\leq \infty$)
 and  $\mathcal{L}_\mathcal{Z} := \{ g \in \mathcal{L}: g(\omega) \in \mathcal{Z}(\omega)
 \subseteq \mathbb{Z},\; \text{for any}\ \omega\in \Omega \}$.
Let $f: \mathbb{Z} \times \Omega \rightarrow \bar{{\R}}$
 be a $\B(\mathbb{Z})\otimes\F$-measurable function.
 Assume that $f(z,\omega)$ is lower semicontinuous in $z$ for every $\omega\in \Omega$
 and
$\mathbb{E}\left[ f(g_0(\omega), \omega) \right]< \infty$ holds for some $g_0\in \mathcal{L}_Z$.
 Then
\bgeqn
\inf_{g \in \mathcal{L}_\mathcal{Z}} \mathbb{E}\left[F_{g}\right]=
\mathbb{E}\left[\inf _{z \in \mathcal{Z}(\omega)
} f(z, \omega)\right],
\label{eq:interchange}
\edeqn
where $F_{g}(\omega):=f(g(\omega), \omega)$.
\end{lemma}

Liu et al. \cite{liu2021multistage} extend the interchangeability principle to a Polish space.

\begin{lemma}[Lemma 1 in \cite{liu2021multistage}]
\label{lemma:interchange-polish-space}
 Consider a Polish space $\mathbb{Z}$ with Borel field $\B(\mathbb{Z})$ and
 a probability space $(\Omega,\F,\mathbb{P})$.
 Let $\mathcal{Z}:\Omega \rightrightarrows \mathbb{Z}$ be a $\F$-measurable set-valued mapping with closed values.
 Let $\mathcal{L}$ be a linear space of measurable functions $ g: \Omega \rightarrow \mathbb{Z}$
 and  $\mathcal{L}_\mathcal{Z} := \{ g \in \mathcal{L}: g(\omega) \in \mathcal{Z}(\omega) \subseteq \mathbb{Z},\; \text{for } a.e.\ \omega\in \Omega \}$.
Let $f: \mathbb{Z} \times \Omega \rightarrow \bar{{\R}}$
 be a Carath\'edory function.
 Suppose that either $\mathbb{E}\left[ \left( \inf_{g\in \mathcal{Z}(\omega)} f(z,\omega) \right)_+ \right]< \infty$ or $\mathbb{E}\left[ \left(- \inf_{g\in \mathcal{Z}(\omega)} f(z,\omega) \right)_+ \right]< \infty$.
 Then
\bgeqn
\inf_{g \in \mathcal{L}_\mathcal{Z}} \mathbb{E}\left[F_{g}\right]=
\mathbb{E}\left[\inf _{z \in \mathcal{Z}(\omega)} f(z, \omega)\right],
\label{eqn:interchange-polish}
\edeqn
where $F_{g}(\omega):=f(g(\omega), \omega)$.
\end{lemma}

\begin{lemma} \label{lemma:polish}
Let $\mathcal{M}([0,1))$ denote the space of probability measures supported on $[0,1)$ endowed with the weak topology. Then $\mathcal{M}([0,1))$ is a Polish space.
\end{lemma}

\noindent
\textbf{Proof.}
\label{proof:Lemma7}
Since any open subset of a Polish space is Polish and $\R$ is Polish, then {$(-1/k,1)$} is Polish for $k=1,2,\dots$.
Observe that
$$
[0,1) = \bigcap_{k=1}^\infty {(-1/k,1)}
$$
is a $G_\delta$ subset of $\R$, i.e., a countable intersection of open subsets of $\R$,
and therefore is Polish, see e.g., Theorem 1 on page 197 in \cite{bourbaki1966elements}.
By Theorem 2.15 in \cite{huber2009robust}, $\mathcal{M}([0,1))$ is Polish.
\hfill $\Box$

{
\subsection{Supplementary results on Wasserstein distance}

\begin{lemma} \label{lemma:bound_wass-dist}
    For probability distributions $P_1,\dots,P_K$ and $\lambda,\lambda' \in \R_+^K$ with $\sum_{k=1}^K \lambda_k=1$ and $\sum_{k=1}^K \lambda'_k=1$,
    \begin{eqnarray} \label{eqn:appendix_wass_result}
        W_r\left(\sum_{k=1}^K \lambda_k P_k, \sum_{k=1}^K \lambda_k' P_k \right)\leq \left(\frac{1}{2} \sum_{k=1}^K  |\lambda_k-\lambda_k'| \right)^{1/r} \max_{0\leq i,j\leq K} W_r(P_i,P_j).
    \end{eqnarray}
\end{lemma}

\noindent
\textbf{Proof.}
For $k=1,\dots,K$, let $\delta_k = (\lambda_k-\lambda_k')_+$ and $\delta'_k = (\lambda_k' - \lambda_k)_+$.
We have $\lambda_k = \min\{\lambda_k,\lambda_k'\}+ \delta_k$, $\lambda_k' = \min\{\lambda_k,\lambda_k'\}+ \delta_k'$, and $\sum_{k=1}^K \delta_k = \sum_{k=1}^K \delta_k'$.

Note that when $\lambda_k = \lambda_k'$ for all $k = 1,\dots,K$,
inequality \eqref{eqn:appendix_wass_result} holds trivially. In what follows, we concentrate on the case that there exists at least one $k$ such that
 $\lambda_k \neq \lambda_k'$, in which case $\omega := \sum_{k=1}^K \delta_k>0$.

For $i,j=1,\dots,K$, let $\rho_{ij} = \frac{\delta_i\delta_j'}{\omega}$.
We have $\sum_{i=1}^K\rho_{ij}= \delta_j'$ and $\sum_{j=1}^K\rho_{ij}= \delta_i$.
By the definition of the Wasserstein distance, for $i,j=1,\dots,K$ and any $\epsilon>0$,
there exists a coupling distribution $\pi_{ij}^*$ such that
\begin{eqnarray*}
    \left(\int \|\xi-\zeta\|^r \pi_{ij}^*(d\xi,d\zeta)\right)^{1/r} < W_r(P_i,P_j) + \epsilon.
\end{eqnarray*}
Let
 $   \hat{\pi} := \sum_{k=1}^K \min\{\lambda_k,\lambda_k'\} \pi_{ii}^* + \sum_{i,j=1}^K  \rho_{ij} \pi_{ij}^*. $
We verify that $\hat{\pi}$
is a measure with marginals
$\sum_{k=1}^K \lambda_k P_i$ and $\sum_{k=1}^K \lambda_k' P_i$.
Observe that
\begin{eqnarray*}
    \int\hat{\pi}(d\xi,d\zeta)&=& \sum_{k=1}^K \min\{\lambda_k,\lambda_k'\} + \sum_{i,j=1}^K \rho_{ij} = \sum_{k=1}^K \min\{\lambda_k,\lambda_k'\} + \sum_{k=1}^K \delta_k \\
    &=& \sum_{k=1}^K \lambda_k =1.
\end{eqnarray*}
Moreover
\begin{eqnarray*}
    \int\hat{\pi}(\xi,d\zeta) &=& \sum_{k=1}^K \min\{\lambda_k,\lambda_k'\}P_i(\xi) + \sum_{i,j=1}^K \rho_{ij} P_i(\xi) \\
&=& \sum_{k=1}^K (\min\{\lambda_k,\lambda_k'\}+ \delta_i)P_i(\xi) = \sum_{k=1}^K \lambda_k P_i(\xi).
\end{eqnarray*}
Likewise, $\int\hat{\pi}(d\xi,\zeta) = \sum_{k=1}^K \lambda_k' P_i(\zeta)$.
By the definition of Wasserstein distance, we have
\begin{eqnarray*}
     &&W_r^r\left(\sum_{k=1}^K \lambda_k P_k, \sum_{k=1}^K \lambda_k' P_k \right) \leq \int \|\xi-\zeta\|^r \hat{\pi}(d\xi,d\zeta) \\
     &&\;
    = \sum_{k=1}^K \min\{\lambda_k,\lambda_k'\} \int \|\xi-\zeta\|^r \pi_{ii}^*(d\xi,d\zeta) + \sum_{i,j=1}^K  \rho_{ij} \int \|\xi-\zeta\|^r \pi_{ij}^*(d\xi,d\zeta)\\
     &&\; < \sum_{k=1}^K \min\{\lambda_k,\lambda_k'\} (W_r(P_i,P_i) + \epsilon)^{r} + \sum_{i,j=1}^K  \rho_{ij} (W_r(P_i,P_j) + \epsilon)^{r}\\
     &&\; \leq \sum_{k=1}^K \min\{\lambda_k,\lambda_k'\} \epsilon^{r} + \sum_{i,j=1}^K  \rho_{ij} (\max_{0\leq i,j \leq K}W_r(P_i,P_j) + \epsilon)^{r}.
\end{eqnarray*}
The first term at rhs of the inequality above goes to $0$ as $\epsilon\to 0$. The second term converges to $\frac{1}{2} \sum_{k=1}^K  |\lambda_k-\lambda_k'| \left(\max_{0\leq i,j\leq K} W_r(P_i,P_j)\right)^r$. This is
 due to the fact that
by the definition of $\delta_{k},\delta_{k}'$ and
$\rho_{ij}$, we have $\sum_{k=1}^K |\lambda_k - \lambda_k'| = \sum_{k=1}^K |\delta_k-\delta_k'| = \sum_{k=1}^K (\delta_k+\delta_k') = 2 \sum_{i,j=1}^K \rho_{ij}$.
Summarizing the discussions, we obtain \eqref{eqn:appendix_wass_result} as desired.
\hfill $\Box$

}

\subsection{Sample average approximation for general hypothesis spaces}
\label{sec:SAA}

In data-driven problems, the true distribution of $(X,Y)$ is unknown
whereas the sample data are
available to decision makers.
This raises a question as to whether
one can use the available sample data to obtain an approximate optimal value and optimal solution of problems
(\ref{eqn:ex_ante}) and (\ref{eqn:ex-post-interchange}).
Let $\{(x_i,y_i)\}_{i=1}^N$ be
i.i.d. samples of $(X,Y)$ and
\begin{eqnarray*}
    P^N(\cdot) = \frac{1}{N} \sum_{i=1}^N \mathds{1}_{(x_i,y_i)}(\cdot),
\end{eqnarray*}
where $\mathds{1}_{(x_i,y_i)}(\cdot)$ denotes Dirac probability measure at $(x_i,y_i)$.
Consider sample average approximation (SAA) of (\ref{eqn:ex_ante})
\bgeqn
\label{eqn:RA-exante-SAA}
\min_{g\in \widetilde{\mathcal{G}}}\rho_{P^N}^{\text{EA}}\big(c(g(X),Y)\big).
\edeqn
In this section, we investigate convergence of the optimal value
obtained from solving
(\ref{eqn:RA-exante-SAA}) to
the optimal value of
problem (\ref{eqn:ex_ante}).
To facilitate the analysis, we concentrate on the case that (\ref{eqn:ex_ante}) takes a particular form
\begin{eqnarray} \label{eqn:primal-form-cont}
    \vt(P) := \min_{(g,s)\in \widetilde{\mathcal{G}}\times \widetilde{\mathcal{S}}}V(\mathbb{E}_{P}[h(g(X),s(X),Y)]),
\end{eqnarray}
where $V:\R \rightarrow \R$ is a non-decreasing continuous function, and $h:\R^{d_z}\times \R \times \R^{d_y} \rightarrow\R$ is a continuous function.
This is because (\ref{eqn:primal-form-cont}) subsumes
two most important risk-averse contextual optimization problems:
risk-averse contextual optimization with entropic risk measure
and risk-averse contextual optimization with expected OCE.
The corresponding SAA problem
can be written as
\begin{eqnarray} \label{eqn:saa-form-cont}
    \vt_N(P_N) :=\min_{(g,s)\in \widetilde{\mathcal{G}}\times \widetilde{\mathcal{S}}} V(\mathbb{E}_{P^N}[h(g(X),s(X),Y)]).
\end{eqnarray}
In particular,
if we let  $V(\mathbb{E}_{P^N}[h(g(X),s(X),Y)]) = \frac{1}{\gamma}\ln(\mathbb{E}_{P^N}[h(g(X),s(X),Y)])$ and $h(z,s,y) = e^{\gamma c(z,y)}$,
then
problem (\ref{eqn:saa-form-cont}) reduces to
risk-averse contextual optimization with entropic risk measure:
\begin{eqnarray} \label{eqn:ent-contextual-opt-emp}
    \min_{g\in \widetilde{\mathcal{G}}} \frac{1}{\gamma} \ln\left(\frac{1}{N} \sum_{i=1}^N e^{\gamma c(g(x_i),y_i)} \right).
\end{eqnarray}
On the other hand, if  $V(\mathbb{E}_{P^N}[h(g(X),s(X),Y)]) = \mathbb{E}_{P^N}[h(g(X),s(X),Y)]$ and $h(z,s,y) = - s - u(-c(z,y)-s)$, then problem (\ref{eqn:saa-form-cont}) reduces to risk-averse contextual optimization with expected OCE:
\begin{eqnarray}
\label{eqn:expected-oce-contextual-opt-emp}
    \min_{g\in\widetilde{\mathcal{G}},s\in \widetilde{\mathcal{S}}} \frac{1}{N} \sum_{i=1}^N\bigg[ - s(x_i) - u\big(-c(g(x_i),y_i)-s(x_i)\big)\bigg].
\end{eqnarray}
Moreover, if $V(\mathbb{E}_{P^N}[h(g(X),s(X),Y)]) = \mathbb{E}_{P^N}[h(g(X),s(X),Y)]$, $h(g(x),s(x),y) = c(g(x),y)$ and $\widetilde{\mathcal{G}}$ is a reproducing kernel Hilbert space, then problem (\ref{eqn:saa-form-cont}) reduces to the risk-neutral problem considered in \cite[Section 5]{bertsimas2022data}.
We begin by establishing the
uniform convergence of
\begin{eqnarray} \label{eqn:saa-inner}
    H_N(g,s) := \mathbb{E}_{P^N}[h(g(X),s(X),Y)]
\end{eqnarray}
to
\begin{eqnarray} \label{eqn:saa-inner-primal}
    H(g,s) := \mathbb{E}_{P}[h(g(X), s(X),Y)].
\end{eqnarray}
To this end,
we need to impose some conditions on $h$ and the hypothesis
spaces $\widetilde{\mathcal{G}},\widetilde{\mathcal{S}}$.

\begin{assumption} \label{assump:SAA_convergence-bounded}
    \begin{enumerate}
        \item[(a)]\label{assump:SAA_convergence-bounded-Zg} For any $\epsilon > 0$, there exists a compact set $\mathcal{X}_\epsilon \times \mathcal{Y}_\epsilon \subseteq \mathcal{X} \times \mathcal{Y}$ such that
        \begin{eqnarray} \label{eqn:assump-sec5-a-growth}
            \sup_{g\in \widetilde{\mathcal{G}},s\in \widetilde{\mathcal{S}}} \mathbb{E}_{P^N}[|h(g(X),s(X),Y)|\mathds{1}_{(\mathcal{X}\times \mathcal{Y})\backslash (\mathcal{X}_\epsilon \times \mathcal{Y}_\epsilon)}(X,Y)] \leq \epsilon, \; \text{w.p.1},
        \end{eqnarray}
        and
        there exist positive numbers $\eta_\epsilon, \Delta_\epsilon$ such that
        \begin{subequations}
        \begin{eqnarray}\label{eqn:assump-equi-continuou}
            \sup_{g\in \widetilde{\mathcal{G}}} \|g(x_1)-g(x_2)\| \leq \epsilon,\; \forall x_1,x_2\in \mathcal{X}_\epsilon \text{ with } \|x_1 - x_2\| \leq \eta_\epsilon,\\
            \sup_{s\in \widetilde{\mathcal{S}}} |s(x_1)-s(x_2)| \leq \epsilon,\; \forall x_1,x_2\in \mathcal{X}_\epsilon \text{ with } \|x_1 - x_2\| \leq \eta_\epsilon,
        \end{eqnarray}
        \end{subequations}
        and $\sup\limits_{g\in \widetilde{\mathcal{G}}_\epsilon} \|g(x)\| \leq \Delta_\epsilon, \sup\limits_{s\in \widetilde{\mathcal{S}}_\epsilon} |s(x)| \leq \Delta_\epsilon$ for all $x\in \mathcal{X}_\epsilon$,
        where $\widetilde{\mathcal{G}}_{\mathcal{X}_\epsilon}$ and $\widetilde{\mathcal{S}}_{\mathcal{X}_\epsilon}$
        are restrictions of $\widetilde{\mathcal{G}}$ and $\widetilde{\mathcal{S}}$ to set $\mathcal{X}_\epsilon$, i.e.,
        \begin{eqnarray*}
            \widetilde{\mathcal{G}}_{\mathcal{X}_\epsilon} := \{g|_{\mathcal{X}_\epsilon}: g\in\widetilde{\mathcal{G}}\},\;
            \widetilde{\mathcal{S}}_{\mathcal{X}_\epsilon} := \{s|_{\mathcal{X}_\epsilon}: s\in\widetilde{\mathcal{S}}\}.
        \end{eqnarray*}

        \item[(b)]\label{assump:SAA_convergence-bounded-h}
        There exist
        a continuous function $r:\mathcal{Y}\rightarrow \R_+$
        and a positive constant $\nu>0$  such that for all $y\in Y$,
        \begin{eqnarray*}
            |h(z_1,s_1,y)-h(z_2,s_2,y)| \leq r(y)\left( \|z_1-z_2\| + |s_1-s_2|\right)^\nu, \forall z_1,z_2\in \mathcal{Z}, s_1,s_2\in \R.
        \end{eqnarray*}

        \item[(c)] \label{assump:SAA-convergence-bounded-mgf} The moment generating function $m(t,g,s) = \mathbb{E}_{P}[e^{t[h(g(X),s(X),Y) -\mathbb{E}_{P}[h(g(X),s(X),Y)]]}]$
        is finite valued for all $g\in\widetilde{\mathcal{G}}$, $s\in \widetilde{\mathcal{S}}$, $t$ in a neighborhood of zero and $\epsilon>0$ close to zero.
    \end{enumerate}
\end{assumption}

Condition (\ref{eqn:assump-sec5-a-growth}) is a kind of uniform integrability condition, see Chapter 3 in \cite{billingsley2013convergence}.
Condition (\ref{eqn:assump-sec5-a-growth}) holds when $\widetilde{\mathcal{G}},\widetilde{\mathcal{S}}$ are bounded and $\mathcal{X}\times \mathcal{Y}$ is a compact set.
In condition (\ref{eqn:assump-equi-continuou}),
we assume the equi-continuity of $\widetilde{\mathcal{G}}$ and $\widetilde{\mathcal{S}}$ over $\mathcal{X}_\epsilon$.
This condition may be fulfilled
under some moderate conditions on $\widetilde{\mathcal{G}},\widetilde{\mathcal{S}}$.
We will return to this
in Section \ref{sec:DR-rkhs}
where $\widetilde{\mathcal{G}},\widetilde{\mathcal{S}}$ are defined within reproducing kernel Hilbert spaces.
Assumption \ref{assump:SAA_convergence-bounded} (b) requires the growth condition of function $h$,
which can be verified in specific problems and is important for the use of Cram\'{e}r’s Large Deviation Theorem.
Assumption \ref{assump:SAA_convergence-bounded} (c) requires that the moment generating function is finite valued, which is a common assumption in SAA literature.
It holds if the probability distribution of the random variable has exponentially decaying tails. In particular, it holds if the random variable is supported on a bounded subset, see \cite{shapiro2008stochastic}.
In Section \ref{sec:consist-saa-rkhs}, we will
show how these assumptions can be verified with
specific hypothesis
spaces, such as the reproducing kernel Hilbert spaces.

\begin{lemma} \label{lemma:saa_convergence_general}
    Let Assumptions \ref{assumption:convex_finite-minimum} (a) and \ref{assump:SAA_convergence-bounded} hold.  Then for any $\delta > 0$, there exist positive constants $\epsilon < \frac{\delta}{6}$, $C(\epsilon,\delta)$ and $\beta(\epsilon,\delta)$, independent of $N$ such that
    \begin{eqnarray} \label{eqn:lemma_saa_convergence_general}
        \inmat{Prob}\left(\sup_{g\in \widetilde{\mathcal{G}},s\in \widetilde{\mathcal{S}}}|H_N(g,s)-H(g,s)|\geq \delta\right) \leq C(\epsilon,\delta)e^{-N{\beta(\epsilon,\delta)}},
    \end{eqnarray}
    for all $N$ sufficiently large,
    where the probability measure `$\inmat{Prob}$' is understood as the product probability measure of $P$ over measurable space $\Xi\times\Xi\times \cdots$ with its product Borel sigma algebra where $\Xi := (\mathcal{X}\times \mathcal{Y})$.
\end{lemma}

With the uniform convergence of the inner term in Lemma \ref{lemma:saa_convergence_general}, we are ready to state
consistency of the objective function $\vt_N(P^N)$ in problem~(\ref{eqn:saa-form-cont}) to $\vt(P)$ in problem~(\ref{eqn:primal-form-cont}).

\begin{theorem} \label{thm: SAA-convergence_V}
    Let Assumptions \ref{assumption:convex_finite-minimum} (a)
    and \ref{assump:SAA_convergence-bounded} hold.
    Then the following assertions hold.
\begin{enumerate}
    \item[(i)] If $V(x)$ is monotone and differentiable and $|V'(t)|$ is upper bounded by $L_V>0$,
    then for any $\delta > 0$, there exist positive constants $\epsilon < \frac{\delta}{6L_V}$, $C(\epsilon,\delta/L_V)$ and $\beta(\epsilon,\delta/L_V)$, independent of $N$ and a positive number $N_0$ such that for all $N\geq N_0$,
    \begin{eqnarray}  \label{eqn:SAA_convergence_V1}
        \inmat{Prob}\left(|\vt_N(P^N)-\vt(P)|\geq \delta\right) \leq C(\epsilon,\delta/L_V)e^{-N{\beta(\epsilon,\delta/L_V)}}.
    \end{eqnarray}
    \item[(ii)] Let $H^*:= \inf\limits_{g\in \widetilde{\mathcal{G}},s \in \widetilde{\mathcal{S}}} \mathbb{E}_{P}[h(g(X), s(X),Y)]$. If the function $V:\R\rightarrow \R$ takes a special form with $V(t) = \ln t$ and $H_N(g,s)>0, H(g,s)>0$, then for any $\delta > 0$, there exist positive constants $\epsilon < \frac{\delta H^*}{6(1+\delta)}$,
    $C(\epsilon,\frac{\delta }{1+ \delta}H^*)$ and $\beta(\epsilon,\frac{\delta }{1+ \delta}H^*)$, independent of $N$, and a positive number $N_0$ such that for all $N\geq N_0$,
    \begin{eqnarray} \label{eqn:SAA_convergence_V}
        \inmat{Prob}\left(|\vt_N(P^N)-\vt(P)|\geq \delta\right) \leq C\left(\epsilon,\frac{\delta }{2+ \delta}H^*\right)e^{-N{\beta(\epsilon,\frac{\delta }{2+ \delta}H^*)}}.
    \end{eqnarray}
\end{enumerate}
\end{theorem}

\noindent
\textbf{Proof.}
See Section \ref{proof:Thm4}.
$\Box$

{
Theorems 2 and 4 both establish the consistency of the SAA method. The distinction is that Theorem 2 restricts the feasible set to the RKHS, $\mathcal{H}^{d_z}\times \mathcal{H}$, whereas Theorem 4 allows for a more general feasible set, $\tilde{\mathcal{G}} \times \tilde{\mathcal{S}}$.
The consistency for the risk-averse contextual optimization problem with expected negative OCE (\ref{eqn:expected-oce-contextual-opt-emp}) and problem with entropic risk measure (\ref{eqn:ent-contextual-opt-emp}) can be obtained by Theorem \ref{thm: SAA-convergence_V} (i) and (ii) respectively.
}

\subsection{Proof of Section \ref{sec:contextual-risk-measure}}

\subsubsection{Proof of Proposition \ref{prop:random_variable_inner_risk}}
\label{proof:Prop1}
Consider $x\in \mathcal{X}$ and $g\in \widetilde{\mathcal{G}}$. Let $z=g(x)$, $Q=P_{Y|X=x}$, and $\theta(z,Q) = (\delta_z,Q)\circ c^{-1}$ be an image probability measure of $\delta_z\otimes Q$ under $c$.
Then, $\rho^{(2)}_{P_{Y|X=x}} ( c(g(x),Y)) = \varrho(\theta(g(x),P_{Y|X=x}))$ as discussed in Section \ref{sec:risk-measure}.
By Assumption \ref{assumption:well-defined_rho} (c) and Lemma~2.74 in \cite{claus2016advancing},
\footnote{
Note that Lemma 2.74 in \cite{claus2016advancing}
is
presented
in the setting
where the underlying decision variables and random variables
are defined in Euclidean spaces $\R^{d_x}, \R^{d_y}, \R^{d_z}$. However, these results also remain valid when the spaces
are replaced
with metric spaces such as
Polish subspaces of $\R^{d_x}, \R^{d_y}, \R^{d_z}$.
This note is applicable to
Lemma~2.81 in \cite{claus2016advancing}, which is to be cited later.}
$\theta(z,Q)$ is continuous in $z$ and $Q$.
Since $z = g(x)$ is continuous in $x$ under
Assumption~\ref{assumption:well-defined_rho}~(b) and $Q=P_{Y|X=x}$ is continuous in $x$ under Assumption \ref{assumption:well-defined_rho} (d),
then $\theta(g(x),P_{Y|X=x})$ is continuous in $x$.
On the other hand,
under Assumption \ref{assumption:well-defined_rho} (a),
$\varrho:\mathcal{M}_1^p\rightarrow\R$ is continuous in $\theta$ by Lemma 2.81 in \cite{claus2016advancing}.
Consequently,
$\rho^{(2)}_{P_{Y|X=x}} ( c(g(x),Y)) = \varrho(\theta(g(x),P_{Y|X=x}))$ is continuous in $x\in\mathcal{X}$.
\hfill $\Box$

\subsubsection{Proof of Proposition \ref{prop:well-definedness-problem}}
\label{proof:Prop2}
By
continuity of functions
$g(\cdot),c(\cdot)$ and the boundedness of
sets $\mathcal{X},\mathcal{Y}$,
$c(g(X),Y)$ is a bounded random function of $(X,Y)$.
{
By Assumption \ref{assumption:well-defined_rho} (a), $\rho_P^{\textbf{EA}}$ is a convex risk measure and thus, by definition, a monetary risk measure.
}
For any $g,g_0\in \widetilde{\mathcal{G}}$,
it follows by Lemma 4.3 in \cite{follmer2011stochastic} that
\begin{eqnarray*}
     |\rho_P^{\text{EA}}\big(c(g(X),Y)\big) - \rho_P^{\text{EA}}\big(c(g_0(X),Y)\big)|
    &\leq& \left\|c(g(X),Y) - c(g_0(X),Y) \right\|_{\infty,\mathcal{X},\mathcal{Y}}\\
    &\leq& \|r\|_{\infty, \mathcal{Y}} \left\| g - g_0 \right\|_{\infty, \mathcal{X}}^\nu,
\end{eqnarray*}
where the second inequality is due to
Assumption~\ref{assumption:well-defined_obj} (b).
This implies the continuity of $\rho_P^{\text{EA}}\big(c(g(X),Y)\big)$ in $g$.
Likewise,  for any $g,g_0\in \widetilde{\mathcal{G}}$,
\begin{eqnarray*}
    &&\left|\rho_{P_X}^{(1)} \left(\rho_{P_{Y|X}}^{(2)} \big(c(g(X),Y)\big)\right) - \rho_{P_X}^{(1)} \left(\rho_{P_{Y|X}}^{(2)} \big(c(g_0(X),Y)\big)\right) \right|\\
    &\leq& \left\| \rho_{P_{Y|X=x}}^{(2)} \big(c(g(x),Y)\big) - \rho_{P_{Y|X=x}}^{(2)} \big(c(g_0(x),Y)\big) \right\|_{\infty,\mathcal{X}}\\
    &\leq&
    \left\| c(g(x),y) - c(g_0(x),y) \right\|_{\infty,\mathcal{X},\mathcal{Y}}\\
    &\leq& \|r\|_{\infty, \mathcal{Y}} \left\| g - g_0 \right\|_{\infty, \mathcal{X}}^\nu,
\end{eqnarray*}
which implies
continuity of $\rho_{P_X}^{(1)}\left(\rho_{P_{Y|X}}^{(2)}\big(c(g(X),Y)\big)\right)$ in $g$.
\hfill $\Box$

\subsection{Proof of Section \ref{sec:contextual-ra-opt}}
\label{ecsec:proof-contextual}

\subsubsection{Proof of Theorem \ref{thm: optimality-of-ent}.}
\label{proof:Thm1}
Part (i).
Similar to the proof of Proposition \ref{prop:random_variable_inner_risk},
we let $Q=P_{Y|X=x}$ and {$\theta(z,Q) = (\delta_z,Q)\circ c^{-1}$ be an image probability measure of $\delta_z\otimes Q$ under $c$.}
Then, $\rho^{(2)}_{P_{Y|X=x}} ( c(z,Y)) = \varrho(\theta(z,P_{Y|X=x}))$.
By Assumption \ref{assumption:well-defined_rho} (c) and Lemma 2.74 in \cite{claus2016advancing},
$\theta(z,Q)$ is continuous in $z$ and $Q$.
Since $Q=P_{Y|X=x}$ is continuous in $x$ under Assumption \ref{assumption:well-defined_rho} (d),
then $\theta(z,P_{Y|X=x})$ is continuous in $(z,x)$.
Under Assumption \ref{assumption:well-defined_rho} (a),
$\varrho:\mathcal{M}_1^p\rightarrow\R$ is continuous in $\theta$ by Lemma 2.81 in \cite{claus2016advancing}.
Therefore,
$\rho^{(2)}_{P_{Y|X=x}} ( c(z,Y))$ is continuous over $\mathcal{Z}\times\mathcal{X}$.
{By Assumption \ref{assumption:convex_finite-minimum}, for any fixed point $\hat{x}$ and all $x$ in a neighborhood of $\hat{x}$, the feasible set of problem \eqref{eqn:conditional_risk_min} can be restricted to a compact subset $\hat{\mathcal{Z}}$, i.e.
\bgeq
\min_{z \in \hat{{\cal Z}}} \ \rho^{(2)}_{P_{Y \mid X=x}}
\bigl(c(z, Y )\bigr).
\edeq
Then by
Berge's maximum theorem \cite{berge1963topological},
the contextual optimal policy mapping $z^*(x)$ is single-valued
and continuous in $x$ over $\mathcal{X}$.}

Part (ii).
    By part (i), $z^*(x)$ is continuous in $x$.
    Therefore, under Assumption \ref{assumption:well-defined_rho} (b), $z^*$ is a feasible solution to {problem (\ref{eqn:ex_ante})}.
    Since $\rho_{P}^{\text{EA}}$ is contextually consistent, $z^*$ is an optimal solution to problem (\ref{eqn:ex_ante}), {which implies the existence of an optimal solution to problem \eqref{eqn:ex_ante}.}
    On the other hand, if $g^*$ is not optimal in $\inf\limits_{z\in Z} \rho_{ P_{Y|X=x}}^{(2)}(c(z,Y))$ on a measurable set of $\mathcal{X}$, then by Assumption \ref{assumption:convex_finite-minimum} (b) and the strict contextual consistency, $\rho_{P}^{\text{EA}}\left(c(g^*(X),Y)\right)$ is strictly larger than $\rho_{P}^{\text{EA}}\left(c(z^*(X),Y)\right)$.
    Therefore, the unique optimal solution to {problem (\ref{eqn:ex_ante})} is $z^*$, which is optimal in all contexts.

Part (iii).
By Part (i), $z^*(x)$ is continuous in $x$.
Therefore, under Assumption \ref{assumption:well-defined_rho} (b), $z^*$ is a feasible solution to problem (\ref{eqn:ex-post-interchange})
with
\begin{eqnarray*}
    \rho_{P_X}^{(1)} \left(\rho_{P_{Y|X}}^{(2)} \big(c(z^*(X),Y)\big)\right) \geq \min\limits_{g\in \widetilde{\mathcal{G}}} \rho_{P_X}^{(1)} \left(\rho_{P_{Y|X}}^{(2)} \big(c(g(X),Y)\big)\right).
\end{eqnarray*}
Next, we
prove $z^*$ is an optimal solution to this problem.
Assume for the sake of
a contradiction that
\begin{eqnarray*}
    z^* \notin \arg\min\limits_{g\in \widetilde{\mathcal{G}}} \rho_{P_X}^{(1)} \left(\rho_{P_{Y|X}}^{(2)} \big(c(g(X),Y)\big)\right).
\end{eqnarray*}
Then
there exists
 $   g^*\in \arg\min\limits_{g\in \widetilde{\mathcal{G}}} \rho_{P_X}^{(1)} \left(\rho_{P_{Y|X}}^{(2)} \big(c(g(X),Y)\big)\right)$
such that
\begin{eqnarray*}
    \rho_{P_X}^{(1)} \left(\rho_{P_{Y|X}}^{(2)} \big(c(g^*(X),Y)\big)\right) < \rho_{P_X}^{(1)} \left(\rho_{P_{Y|X}}^{(2)} \big(c(z^*(X),Y)\big)\right).
\end{eqnarray*}
Since $ \rho_{P_X}^{(1)}$ is a convex risk measure and is thus monotonic, the inequality above implies that
\begin{eqnarray*}
    \rho_{P_{Y|X=x}}^{(2)} \big(c(g^*(x),Y\big) < \rho_{P_{Y|X=x}}^{(2)} \big(c(z^*(x),Y)\big),
\end{eqnarray*}
holds in a subset of $\mathcal{X}$ with non-zero measure,
which contradicts to
the fact that $ z^*(x)\in \arg\min\limits_{z\in \mathcal{Z}} \rho_{ P_{Y|X=x}}(c(z,Y))$
for all $x\in {\cal X}$.
This shows that
$z^*$ is the optimal solution to problem (\ref{eqn:ex-post-interchange}) as desired,
{which implies the existence of an optimal solution to problem \eqref{eqn:ex-post-interchange}.}
Moreover,
by the strict monotonicity of $\rho_{P_X}^{(1)}$,
for
each  $g^*\in \min\limits_{g\in \widetilde{\mathcal{G}}} \rho_{P_X}^{(1)} \left(\rho_{P_{Y|X}}^{(2)} \big(c(g(X),Y)\big)\right)$
\begin{eqnarray*}
    \rho_{P_{Y|X=x}}^{(2)} \big(c(g^*(x),Y)\big) = \rho_{P_{Y|X=x}}^{(2)} \big(c(z^*(x),Y)\big) = \min\limits_{z\in Z} \rho_{P_{Y|X=x}}^{(2)} \big(c(z,Y)\big)
\end{eqnarray*}
for almost every $x\in {\cal X}$.
The proof is complete.
\hfill $\Box$

\subsection{Proof of Section \ref{sec:DR-rkhs}}

\subsubsection{Proof of Lemma \ref{lemma:multidimensional-representer}.}
\label{proof:Lemma4}
{For a given
sample $\left\{ x_i \right\}_{i=1}^N$ of size $N$, }
define
$$\hat{\mathcal{H}}^{d} =  \left\{ \hat{f}\in \mathcal{H}^{d} | \hat{f}^{(t)} = \sum_{i=1}^N k(x_i,\cdot) \alpha_i^{(t)},\; \alpha_1^{(t)},\dots,\alpha_N^{(t)}\in \R,
t=1,\dots, d \right\}.$$
$\hat{\mathcal{H}}^{d}$ is a closed subspace of $\mathcal{H}^{d}$.
Therefore, any $f\in \mathcal{H}^{d}$ can be decomposed into a part $\hat{f}$ that lies in $\hat{\mathcal{H}}^{d}$ and a part $\hat{f}^\perp \in \mathcal{H}^{d}$ which is orthogonal to $\hat{\mathcal{H}}^{d}$, i.e.,
\begin{eqnarray*}
    f^{(t)}(\cdot) = \hat{f}^{(t)}(\cdot)+\hat{f}^{{(t)}^\perp}(\cdot) = \sum_{i=1}^N k(x_i,\cdot) \alpha_i^{(t)} + \hat{f}^{{(t)}^\perp}(\cdot),\; t=1,\dots, d
\end{eqnarray*}
with
\begin{eqnarray*}
    \langle  \hat{f}^{{(t)}^\perp}, k(x_i,\cdot) \rangle = 0, \forall i = 1,\dots,N,\; \inmat{for}\; t=1,\dots, d.
\end{eqnarray*}
By Definition \ref{def: RKHS},
\begin{eqnarray*}
    f^{(t)}(x_i) = \hat{f}^{(t)}(x_i) + \hat{f}^{{(t)}^\perp}(x_i) = \hat{f}^{(t)}(x_i) + \langle \hat{f}^{{(t)}^\perp},k(x_i,\cdot) \rangle = \hat{f}^{(t)}(x_i),\; \inmat{for}\; t=1,\dots, d.
\end{eqnarray*}
Thus, for any function $f\in \mathcal{H}^{d}$,
\begin{eqnarray} \label{eqn:proof-multidim-representer-first-part}
    &&W((x_1,y_1,f(x_1)),\dots,(x_N,y_N,f(x_N))) \nonumber\\
    &&\qquad = W((x_1,y_1,\hat{f}(x_1)),\dots,(x_N,y_N,\hat{f}(x_N))),
\end{eqnarray}
and, by the monotonicity of $R$,
\begin{eqnarray*}
    R(\|f\|_{\mathcal{H}^{d}}) &=& R\left( \sqrt{ \sum_{t=1}^{d} \|f^{(t)}\|^2_{\mathcal{H}}} \right) = R\left( \sqrt{ \sum_{t=1}^{d} \left\|\sum_{i=1}^N k(x_i,\cdot) \alpha_i^{(t)} + \hat{f}^{{(t)}^\perp} \right\|^2_{\mathcal{H}}} \right)\\
    &=& R\left( \sqrt{ \sum_{t=1}^{d} \left(\left\|\sum_{i=1}^N k(x_i,\cdot) \alpha_i^{(t)}\right\|_\mathcal{H}^2 + \left\|\hat{f}^{{(t)}^\perp}\right\|^2_{\mathcal{H}}\right)} \right)\\
    &\geq& R\left( \sqrt{ \sum_{t=1}^{d} \left\|\sum_{i=1}^N k(x_i,\cdot) \alpha_i^{(t)}\right\|_\mathcal{H}^2} \right),
\end{eqnarray*}
where
equality
holds
iff $\hat{f}^{{(t)}^\perp} = 0,
\inmat{for}\; t = 1,\dots, d$.
Since the first term of (\ref{eqn:proof-multidim-representer-ec}) is independent of $\hat{f}^{{(t)}^\perp}$ as shown in (\ref{eqn:proof-multidim-representer-first-part}) and the second term is minimized by $\hat{f}^{{(t)}^\perp} = 0, \forall  t=1,\dots,d$, then any minimizer of (\ref{eqn:proof-multidim-representer-ec}) must have $\hat{f}^{{(t)}^\perp} = 0, \inmat{for}\; t=1,\dots, d$ and we arrive at (\ref{eqn:optimal-solution-form-ec}).
By plugging the optimal form (\ref{eqn:optimal-solution-form-ec}) into (\ref{eqn:proof-multidim-representer-ec}), we immediately arrive at (\ref{eqn:optimal-alpha_form-ec}).
\hfill $\Box$

\subsubsection{Proof of Theorem \ref{theorem:saa_convergence_rkhs}.}
\label{proof:Thm2}
We only prove Part (i).
Parts (ii) and (iii) can be proved based on
Part (i) (analogous to
 the proof of Theorem \ref{thm: SAA-convergence_V}).
Under Assumption \ref{assump:saa-convergence-unbounded} (a), for any $\epsilon >0 $,
there exists a constant $r>0$ such that
\begin{eqnarray*}
    \int_{\mathcal{X} \times \mathcal{Y}} \phi(x,y) \mathds{1}_{(r,\infty)}(\phi(x,y))P(dxdy) \leq \epsilon.
\end{eqnarray*}
Since $\phi$ is coercive, then
there exists a compact set $(\mathcal{X}_\epsilon \times \mathcal{Y}_\epsilon) \in \mathcal{X}\times \mathcal{Y}$ such that
$\{(x,y)\in \mathcal{X}\times \mathcal{Y}:\phi(x,y) \leq r \} \subseteq (\mathcal{X}_\epsilon \times \mathcal{Y}_\epsilon)$ and thus
\begin{eqnarray*}
    \int_{(\mathcal{X}\times \mathcal{Y}) \setminus (\mathcal{X}_\epsilon\times \mathcal{Y}_\epsilon)} \phi(X,Y) P(dxdy) \leq \int_{\mathcal{X}\times \mathcal{Y}} \phi(X,Y) \mathds{1}_{(r,\infty)}(\phi(X,Y))P(dxdy) \leq \epsilon.
\end{eqnarray*}
By Cram\' er's large deviation theory, there exist positive numbers $C_0$ and $\gamma_0$ such that
\begin{eqnarray*}
    \text{Prob}\left( \int_{(\mathcal{X} \times \mathcal{Y}) \setminus (\mathcal{X}_\epsilon \times \mathcal{Y}_\epsilon)} \phi(x,y) P^N(dxdy) \geq 2\epsilon \right) \leq C_0 e^{-\gamma_0 N},
\end{eqnarray*}
and thus
\begin{eqnarray} \label{eqn:saa-unbounded-restricted}
    && \text{Prob}\left( \sup_{g\in \widetilde{\mathcal{G}}}\int_{(\mathcal{X} \times \mathcal{Y}) \setminus (\mathcal{X}_\epsilon \times \mathcal{Y}_\epsilon)} h(z(x),s(x),y) P^N(dxdy) \geq 2\epsilon \right) \nonumber\\
    &\leq&
    \text{Prob}\left( \int_{(\mathcal{X} \times \mathcal{Y}) \setminus (\mathcal{X}_\epsilon \times \mathcal{Y}_\epsilon)} \phi(x,y) P^N(dxdy) \geq 2\epsilon \right) \leq C_0 e^{-\gamma_0 N}.
\end{eqnarray}

Next, we prove the relative compactness of
the hypothesis set
restricted to $\mathcal{X}_\epsilon \times \mathcal{Y}_\epsilon$.
For this
purpose,
we first prove the equi-continuity of the set.
Under Assumption \ref{assump:saa-convergence-unbounded} (b) and (c), for any $\epsilon>0$, there exists $\eta_\epsilon>0$ such that for any $z(\cdot) \in \bar{\mathcal{H}}_z^{d_z}$,
    {
\begin{eqnarray*}
    &&\| z(x) - z(x') \|^2 \\
    &&\qquad = \sum_{t=1}^{d_z}| z_{t}(x) - z_{t}(x') |^2 = \sum_{t=1}^{d_z}| \langle z_{t}, k(\cdot, x) \rangle - \langle z_{t}, k(\cdot, x') \rangle |^2\\
    &&\qquad \leq  \sum_{t=1}^{d_z} \|z_t\|_{\mathcal{H}^{d_z}}^2 \| k(\cdot, x) - k(\cdot, x') \|_\mathcal{H}^2
    \leq \beta_1^2 \epsilon^2,
    \; \forall x,x'\in \mathcal{X}_\epsilon \text{ with } \|x-x'\|\leq \eta_\epsilon.
\end{eqnarray*}
}
Likewise, we can show that for any $s(\cdot)\in \bar{\mathcal{H}}_s$,
\begin{eqnarray*}
    |s(x)-s(x')| \leq \beta_2 \epsilon,\; \forall x,x'\in \mathcal{X}_\epsilon \text{ with } \|x-x'\|\leq \eta_\epsilon.
\end{eqnarray*}
This implies the equi-continuity of the functions in
the hypothesis  set restricted to
$\mathcal{X}_\epsilon$.
On the other hand, for any $z \in \bar{\mathcal{H}}_z^{d_z}$,
{
\begin{eqnarray*}
    \sup_{x\in \mathcal{X}_\epsilon} \| z(x) \|^2 = \sup_{x\in \mathcal{X}_\epsilon} \sum_{t = 1} ^ {d_z} | \langle z_t, k(\cdot,x) \rangle |^2 \leq \| z \|_{\mathcal{H}^{d_z}}^2 \sup_{x\in \mathcal{X}_\epsilon} \| k(\cdot, x) \|_\mathcal{H}^2 \leq \beta_1^2 \sup_{x\in \mathcal{X}_\epsilon} \| k(\cdot, x) \|_\mathcal{H}^2.
\end{eqnarray*}
}
Similarly, for any $s \in \bar{\mathcal{H}}_s$, we can show that
 $   \sup_{x\in \mathcal{X}_\epsilon} | s(x) | \leq \beta_2 \sup_{x\in \mathcal{X}_\epsilon} \| k(\cdot, x) \|_\mathcal{H}.$
Observe that $\| k(\cdot, x) \|_\mathcal{H} = \sqrt{k(x,x)}$ is bounded on a compact set $\mathcal{X}_\epsilon$. Therefore, we arrive at the uniform boundedness of the hypothesis sets
restricted to $\mathcal{X}_\epsilon$, i.e., $\sup\limits_{z(\cdot)\in \bar{\mathcal{H}}_z^{d_z}}\|z(\cdot)\|_\infty$ and $\sup\limits_{s(\cdot)\in \bar{\mathcal{H}}_s}\|s(\cdot)\|_\infty$ are bounded on $\mathcal{X}_\epsilon$.
By Ascoli-Arzela Theorem \cite[Theorem A5]{rudin1973functional}, the hypothesis sets
restricted to
$\mathcal{X}_\epsilon$ are  relatively compact.
Thus, we can choose an $\bar\epsilon$-net
of $\bar{\mathcal{H}}_z^{d_z}\times\bar{\mathcal{H}}_s$, i.e., for any $\bar{\epsilon}>0$, there exists a set of finite number of functions $\{(z_k(\cdot),s_k(\cdot))\}_{k=1}^K\subset \bar{\mathcal{H}}_z^{d_z}\times\bar{\mathcal{H}}_s$
such that $\bar{\mathcal{H}}_z^{d_z}\times\bar{\mathcal{H}}_s = \cup_{k=1}^K(\bar{\mathcal{H}}_z^{d_z}\times\bar{\mathcal{H}}_s)_k^{\bar{\epsilon}}$,
where $$(\bar{\mathcal{H}}_z^{d_z}\times\bar{\mathcal{H}}_s)_k^{\bar{\epsilon}} := \{(z(\cdot),s(\cdot))\in \bar{\mathcal{H}}_z^{d_z}\times\bar{\mathcal{H}}_s| \sup_{x\in \mathcal{X}_\epsilon}(||z(x)-z_k(x)||+ |s(x)-s_k(x)|) \leq\bar{\epsilon}\}.
$$
To facilitate the exposition, let
\begin{eqnarray*}
    R_1(z_k,s_k,\epsilon,P^N) &:=& \mathbb{E}_{P^N}[h(z_k(X),s_k(X),Y)\mathds{1}_{\mathcal{X}_\epsilon \times \mathcal{Y}_\epsilon}(X,Y)] \\
    && - \mathbb{E}_{P}[h(z_k(X),s_k(X),Y)\mathds{1}_{\mathcal{X}_\epsilon \times \mathcal{Y}_\epsilon}(X,Y)].
\end{eqnarray*}
Then
\begin{eqnarray} \label{eqn:proof-sec6-2}
    &&\sup_{(z(\cdot),s(\cdot))\in \bar{\mathcal{H}}_z^{d_z}\times\bar{\mathcal{H}}_s} \left|\mathbb{E}_{P^N}[h(z(X),s(X),Y)\mathds{1}_{\mathcal{X}_\epsilon \times \mathcal{Y}_\epsilon}(X,Y)]- \mathbb{E}_{P}[h(z(X),s(X),Y)\mathds{1}_{\mathcal{X}_\epsilon \times \mathcal{Y}_\epsilon}(X,Y)]\right|\nonumber\\
    &\leq& \sup_{k\in\{1,\dots,K\}}\sup_{(z(\cdot),s(\cdot))\in (\bar{\mathcal{H}}_z^{d_z}\times\bar{\mathcal{H}}_s)_k^{\bar{\epsilon}}} \big\{ \mathbb{E}_{P^N}\left[|h(z(X),s(X),Y)-h(z_k(X),s_k(X),Y)|\mathds{1}_{\mathcal{X}_\epsilon \times \mathcal{Y}_\epsilon}(X,Y)\right]\nonumber\\
    && +  \left|\mathbb{E}_{P^N}[h(z_k(X),s_k(X),Y)\mathds{1}_{\mathcal{X}_\epsilon \times \mathcal{Y}_\epsilon}(X,Y)]-\mathbb{E}_{P}[h(z_k(X),s_k(X),Y)\mathds{1}_{\mathcal{X}_\epsilon \times \mathcal{Y}_\epsilon}(X,Y)]\right| \nonumber\\
    && + \mathbb{E}_{P}\left[|h(z_k(X),s_k(X),Y) - h(z(X),s(X),Y)|\mathds{1}_{\mathcal{X}_\epsilon \times \mathcal{Y}_\epsilon}(X,Y)\right]\big\} \nonumber\\
    &\leq & \mathbb{E}_{P^N}[r(X,Y){\bar{\epsilon}}^\nu] + \sup_{k\in\{1,\dots,K\}} \left| R_1(z_k,s_k,\epsilon,P^N)\right| + \mathbb{E}_{P}[r(X,Y){\bar{\epsilon}}^\nu]\nonumber\\
    &\leq &  \sup_{k\in\{1,\dots,K\}} \left| R_1(z_k,s_k,\epsilon,P^N)\right|+2 \epsilon,
\end{eqnarray}
where the last inequality holds by letting $\bar{\epsilon}\leq (\epsilon/\bar{r})^{1/\nu}$ and $\bar{r} := \max\{|r(x,y)|:(x,y)\in \mathcal{X}_\epsilon \times \mathcal{Y}_\epsilon\}$.
Consequently, we have
\begin{eqnarray*}
    && \sup_{(z(\cdot),s(\cdot))\in \bar{\mathcal{H}}_z^{d_z}\times\bar{\mathcal{H}}_s} \left|\mathbb{E}_{P^N}[h(z(X),s(X),Y)] +\lambda_N R(z,s) - \mathbb{E}_{P}[h(z(X),s(X),Y)]\right|\\
    &\leq &  \sup_{(z(\cdot),s(\cdot))\in \bar{\mathcal{H}}_z^{d_z}\times\bar{\mathcal{H}}_s} \left|\mathbb{E}_{P^N}[h(z(X),s(X),Y)\mathds{1}_{\mathcal{X}_\epsilon \times \mathcal{Y}_\epsilon}(X,Y)]- \mathbb{E}_{P}[h(z(X),s(X),Y)\mathds{1}_{\mathcal{X}_\epsilon \times \mathcal{Y}_\epsilon}(X,Y)]\right|\\
    && + \sup_{(z(\cdot),s(\cdot))\in \bar{\mathcal{H}}_z^{d_z}\times\bar{\mathcal{H}}_s}  \int_{(\mathcal{X} \times \mathcal{Y}) \setminus (\mathcal{X}_\epsilon \times \mathcal{Y}_\epsilon)} |h(g(X),s(X),Y)| P^N(dxdy) \\
    && + \sup_{(z(\cdot),s(\cdot))\in \bar{\mathcal{H}}_z^{d_z}\times\bar{\mathcal{H}}_s}  \int_{(\mathcal{X} \times \mathcal{Y}) \setminus (\mathcal{X}_\epsilon \times \mathcal{Y}_\epsilon)} |h(g(X),s(X),Y)| P(dxdy) + \epsilon
    \\
    &\leq & \sup_{k\in\{1,\dots,K\}} \left| R_1(z_k,s_k,\epsilon,P^N)\right| \\
    && + \sup_{(z(\cdot),s(\cdot))\in \bar{\mathcal{H}}_z^{d_z}\times\bar{\mathcal{H}}_s}  \int_{(\mathcal{X} \times \mathcal{Y}) \setminus (\mathcal{X}_\epsilon \times \mathcal{Y}_\epsilon)} |h(z(X),s(X),Y)| P^N(dxdy) + 4\epsilon,
\end{eqnarray*}
where the first inequality holds due to $\lambda_N\leq \epsilon/(\beta_1^2+\beta_2^2)$ for $N\geq N_0$ and $R(z,s) = \sum_{t=1}^{d_z} \|z^{(t)}\|_\mathcal{H}^2 + \|s\|_\mathcal{H}^2 \leq \beta_1^2+\beta_2^2$, and the second inequality comes from (\ref{eqn:proof-sec6-2}).
On the other hand, by setting $\epsilon < \frac{\delta}{6}$, we have
\begin{eqnarray} \label{eqn:proof-theorem6.1-1}
    &&\inmat{Prob} \left( \sup_{k\in\{1,\dots,K\}} \left| R_1(z_k,s_k,\epsilon,P^N)\right|  \geq \delta - 6\epsilon \right)\nonumber\\
    &\leq& \sum_{k=1}^K \inmat{Prob} \bigg( \left|  R_1(z_k,s_k,\epsilon,P^N)\right|  \geq \delta  - 6\epsilon \bigg)\nonumber\\
    &\leq& \sum_{k=1}^K \bigg(\inmat{Prob} \bigg(   R_1(z_k,s_k,\epsilon,P^N)  \geq \delta - 6\epsilon \bigg) + \inmat{Prob} \bigg(   R_1(z_k,s_k,\epsilon,P^N)  \leq -( \delta - 6\epsilon) \bigg)\bigg)\nonumber\\
    &\leq& 2\sum_{k=1}^K \left(e^{-N [I(\delta - 6\epsilon,z_k,s_k)\wedge I(-(\delta - 6\epsilon),z_k,s_k)]}\right),
\end{eqnarray}
where the last inequality comes from Cram\'{e}r’s Large Deviation Theorem with $I(u,z_k,s_k): = \sup_{t\in \R} \{ut-\log m(t,z_k,s_k)\}$ and
\begin{eqnarray*}
    m(t,z_k,s_k) &:=& \mathbb{E}_{P}[e^{t[h(z_k(X),s_k(X),Y) -\mathbb{E}_{P}[h(z_k(X),s_k(X),Y)]]}].
\end{eqnarray*}
Combining (\ref{eqn:saa-unbounded-restricted}) and (\ref{eqn:proof-theorem6.1-1}),
we obtain
\begin{eqnarray*}
    &&\inmat{Prob}\left( \sup_{(z(\cdot),s(\cdot))\in \bar{\mathcal{H}}_z^{d_z}\times\bar{\mathcal{H}}_s} \left| \mathbb{E}_{P^N}[h(z(X),s(X),Y)] + \lambda_N R(z,s)- \mathbb{E}_{P}[h(z(X),s(X),Y)]\right| \geq \delta \right)\\
    &\leq& \inmat{Prob} \left( \sup_{k\in\{1,\dots,K\}} \left| R_1(z_k,s_k,\epsilon,P^N)\right|  \geq \delta - 6\epsilon \right)\\
    && + \text{Prob}\left( \sup_{(z(\cdot),s(\cdot))\in \bar{\mathcal{H}}_z^{d_z}\times\bar{\mathcal{H}}_s} \int_{(\mathcal{X} \times \mathcal{Y}) \setminus (\mathcal{X}_\epsilon \times \mathcal{Y}_\epsilon)} |{h(z(X),s(X),Y)}| P(dxdy) \geq 2\epsilon \right) \\
    &\leq& 2\sum_{k=1}^K \left(e^{-N [I(\delta - 6\epsilon,z_k,s_k)\wedge I(-(\delta - 6\epsilon),z_k,s_k)]}\right) + C_0 e^{-\gamma_0 N},
\end{eqnarray*}
Note that $K$ depends on $\epsilon$ and thus on $\delta$.
By letting $C(\epsilon, \delta) = 2K + C_0 $ and $\beta(\epsilon,\delta) = \min\{\gamma_0,\min\limits_k\{I(\delta - 6\epsilon,z_k,s_k)\wedge I(-(\delta - 6\epsilon),z_k,s_k)\}\}$, we arrive at (\ref{eqn:saa_unbounded-rkhs-uniform-boundedness-rkhs}).
\hfill $\Box$

\subsubsection{Proof of Theorem \ref{theorem:optimal-rkhs}.}
\label{proof:Thm3}
Part (i).
By the definition of the universal kernel, for any $\bar{\epsilon}>0$, there exists $(\hat{z}(\cdot),\hat{s}(\cdot))\in \mathcal{H}^{d_z}\times \mathcal{H}$ { approximating $(z^*(\cdot),s^*(\cdot))$ componentwise} such that
{
\begin{eqnarray*}
    \sup_{x\in \mathcal{X}} \left|(\hat{z}(x))^{(t)} - (z^*(x))^{(t)}\right| &\leq& \frac{\bar\epsilon}{2\sqrt{d_z}},\\
    \sup_{x\in \mathcal{X}} |\hat{s}(x) - s^*(x)| &\leq& \frac{\bar\epsilon}{2},
\end{eqnarray*}
and thus
}
$\sup\limits_{x\in \mathcal{X}} (\|\hat{z}(x) - z^*(x)\| + |\hat{s}(x) - s^*(x)|) \leq \bar{\epsilon}$.
Under Assumption \ref{assump:saa-convergence-unbounded} (d),
\begin{eqnarray*}
    && \mathbb{E}_P[h(\hat{z}(X),\hat{s}(X),Y)] - \mathbb{E}_P[h(z^*(X),s^*(X),Y)]\\
    & = & \left|\mathbb{E}_P[h(\hat{z}(X),\hat{s}(X),Y)] - \mathbb{E}_P[h(z^*(X),s^*(X),Y)]\right|\\
    &\leq & \int_{\mathcal{X}\times \mathcal{Y}} |h(\hat{z}(X),\hat{s}(X),Y) - h(z^*(X),s^*(X),Y)| P(dxdy) \\
    &\leq& \int_{\mathcal{X}\times \mathcal{Y}} {r(Y)} (\|\hat{z}(X) - z^*(X)\| + |\hat{s}(X) - s^*(X)|)^\nu P(dxdy)\\
    &\leq& \int_{\mathcal{X}\times \mathcal{Y}} {r(Y)}{\bar{\epsilon}}^\nu P(dxdy) \leq \delta_1,
\end{eqnarray*}
where the first inequality comes from the fact that $(\hat{z}(\cdot),\hat{s}(\cdot))\in \bar{\mathcal{H}}_z^{d_z}\times\bar{\mathcal{H}}_s \subseteq {\mathcal{C}(\mathcal{X})}^{d_z}\times \mathcal{C}(\mathcal{X})$, and the last inequality holds by letting $\bar{\epsilon}\leq (\delta_1/\bar{r})^{1/\nu}$ and {$\bar{r} := \mathbb{E}_P[r(Y)].$}

Part (ii).
{
Before proceeding to the proof, we recall the model and associated notations to facilitate reading.
Since we focus on the setting where $V(t)=t$, the optimal value under the true distribution specified in (\ref{eqn:saa-form-rkhs-primal}) is
\begin{eqnarray*}
    \vt(P) &:=& \min_{(z(\cdot),s(\cdot))\in \mathcal{H}^{d_z}\times\mathcal{H}} \mathbb{E}_{P}[h(z(X),s(X),Y)],
\end{eqnarray*}
and
the optimal value of the SAA problem given in (\ref{eqn:saa-form-rkhs}) is
\begin{eqnarray*}
    \vt_N(P_N,\lambda_N) := \min_{(z(\cdot),s(\cdot))\in \mathcal{H}^{d_z}\times\mathcal{H}} \mathbb{E}_{P^N}[h(z(X),s(X),Y)]+ \lambda_N R(z,s).
\end{eqnarray*}
}
Next, we consider the hypothesis set $\mathcal{H}^{d_z}$ for $z(\cdot)$ with $\|z(\cdot)\|_{\mathcal{H}^{d_z}}\leq \beta_1$
and $\mathcal{H}$ for $s(\cdot)$ with $\|s(\cdot)\|_{\mathcal{H}}\leq \beta_2$.
Since $(\hat{z}(\cdot),\hat{s}(\cdot))\in \bar{\mathcal{H}}_z^{d_z}\times\bar{\mathcal{H}}_s$ and $\bar{\mathcal{H}}_z^{d_z}\times\bar{\mathcal{H}}_s \subseteq {\mathcal{C}(\mathcal{X})}^{d_z}\times \mathcal{C}(\mathcal{X})$, then
\begin{eqnarray*}
     \mathbb{E}_P[h(z^*(X),s^*(X),Y)]
    \leq  \vt(P)
    \leq  \mathbb{E}_P[h(\hat{z}(X),\hat{s}(X),Y)]  ,
\end{eqnarray*}
and thus
\begin{eqnarray} \label{eqn:rkhs-optimal-proof-1}
\vt(P) - \mathbb{E}_P[h(z^*(X),s^*(X),Y)] \leq \delta_1.
\end{eqnarray}
By Theorem \ref{theorem:saa_convergence_rkhs} (i),
for any $\delta_2 > 0$, there exist positive constants $\epsilon < \frac{\delta_2}{6}$, $\hat{C}(\epsilon,\delta_2,\delta_1)$ and $\beta(\epsilon,\delta_2,\delta_1)$, independent of $N$ such that
    \begin{eqnarray}  \label{eqn:rkhs-optimal-proof-2}
        && \inmat{Prob}\left( \left| \vt_N(P_N,\lambda_N) - \vt(P) \right|\geq \delta_2\right)\nonumber\\
        &\leq& \inmat{Prob}\left(\sup_{(z(\cdot),s(\cdot))\in \bar{\mathcal{H}}_z^{d_z}\times\bar{\mathcal{H}}_s} \bigg| \mathbb{E}_{P^N}[h(z(X),s(X),Y)] + \lambda_N R(z,s) - \mathbb{E}_{P}[h(z(X),s(X),Y)] \bigg|\geq \delta_2\right)\nonumber\\
        &\leq& \hat{C}(\epsilon,\delta_2,\delta_1)e^{-N{\hat\beta(\epsilon,\delta_2,\delta_1)}},
    \end{eqnarray}
where $\hat{C}(\epsilon,\delta_2,\delta_1)$ and $\hat\beta(\epsilon,\delta_2,\delta_1)$
depend on the upper bounds of RKHS norms $\beta_1$ and $\beta_2$ which further depend on the choice of $(\hat{z},\hat{s})$ and thus on $\delta_1$.
Finally, we have
\begin{eqnarray*}
    && \inmat{Prob}\bigg( \bigg|\vt_N(P_N,\lambda_N) - \mathbb{E}_{P}[h(z^*(X),s^*(X),Y)] \bigg| \geq \delta_1 + \delta_2 \bigg)\\
    &\leq & \inmat{Prob}\bigg( \bigg|\vt_N(P_N,\lambda_N)  - \vt(P) \bigg| + \bigg| \vt(P)  -\mathbb{E}_{P}[h(z^*(X),s^*(X),Y)] \bigg| \geq \delta_1+\delta_2 \bigg)\\
    &\leq& \inmat{Prob}\bigg( \bigg|\vt_N(P_N,\lambda_N)  - \vt(P) \bigg| + \delta_1 \geq \delta_1+\delta_2 \bigg)\\
    & = & \inmat{Prob}\bigg( \bigg| \vt_N(P_N,\lambda_N)  - \vt(P) \bigg|\geq \delta_2 \bigg) \\
    &\leq& C(\epsilon,\delta_2,\delta_1)e^{-N{\beta(\epsilon,\delta_2,\delta_1)}},
\end{eqnarray*}
where the second inequality comes from (\ref{eqn:rkhs-optimal-proof-1}) and the last inequality comes from (\ref{eqn:rkhs-optimal-proof-2}).
\hfill $\Box$

\subsection{Proof of Section \ref{sec:SAA}}

\subsubsection{Proof of Lemma \ref{lemma:saa_convergence_general}.}
\label{proof:Lemma8}
Under Assumption \ref{assump:SAA_convergence-bounded} (d),
for any $\epsilon>0$, there exists a compact set $\mathcal{X}_\epsilon \times \mathcal{Y}_\epsilon \subseteq \mathcal{X} \times \mathcal{Y}$ such that
$\widetilde{\mathcal{G}}_{\mathcal{X}_\epsilon}$ and $\widetilde{\mathcal{S}}_{\mathcal{X}_\epsilon}$ are uniformly bounded and  equi-continuous on $\mathcal{X}_\epsilon$.
By Ascoli-Arzela Theorem \cite[Theorem A5]{rudin1973functional}, $\widetilde{\mathcal{G}}_{\mathcal{X}_\epsilon}\times \widetilde{\mathcal{S}}_{\mathcal{X}_\epsilon}$ is relatively compact, which means that there exists an $\bar{\epsilon}$-net of $\widetilde{\mathcal{G}}_{\mathcal{X}_\epsilon}\times \widetilde{\mathcal{S}}_{\mathcal{X}_\epsilon}$, i.e., for any $\bar{\epsilon}>0$, there exists a finite set of functions $\{(g_k,s_k)\}_{k=1}^K\subset \widetilde{\mathcal{G}}_{\mathcal{X}_\epsilon} \times \widetilde{\mathcal{S}}_{\mathcal{X}_\epsilon}$ such that $\widetilde{\mathcal{G}}_{\mathcal{X}_\epsilon} = \cup_{k=1}^K(\widetilde{\mathcal{G}})_k^{\bar{\epsilon}}$ and $\widetilde{\mathcal{S}}_{\mathcal{X}_\epsilon} = \cup_{k=1}^K(\widetilde{\mathcal{S}})_k^{\bar{\epsilon}}$, where $(\widetilde{\mathcal{G}}_{\mathcal{X}_\epsilon})_k^{\bar{\epsilon}} := \{g\in \widetilde{\mathcal{G}}_{\mathcal{X}_\epsilon}| \|g-g_k\|_\infty\leq {2^{\frac{1-\nu}{\nu}}}\bar{\epsilon}\}$ and $(\widetilde{\mathcal{S}}_{\mathcal{X}_\epsilon})_k^{\bar{\epsilon}} := \{s\in \widetilde{\mathcal{S}}_{\mathcal{X}_\epsilon}| \|s-s_k\|_\infty\leq {2^{\frac{1-\nu}{\nu}}}\bar{\epsilon}\}$.
To facilitate the exposition, let
\begin{eqnarray*}
    R_1(g_k,s_k,\epsilon,P^N) &:=& \mathbb{E}_{P^N}[h(g_k(X),s_k(X),Y)\mathds{1}_{\mathcal{X}_\epsilon \times \mathcal{Y}_\epsilon}(X,Y)] \\
    && - \mathbb{E}_{P}[h(g_k(X),s_k(X),Y)\mathds{1}_{\mathcal{X}_\epsilon \times \mathcal{Y}_\epsilon}(X,Y)].
\end{eqnarray*}
By the definition of $H_N(g,s)$ and $H(g,s)$,
\begin{eqnarray*}
    && \sup_{g\in\widetilde{\mathcal{G}}, s\in\widetilde{\mathcal{S}}} |H_N(g,s)-H(g,s)|\\
    &=& \sup_{g\in \widetilde{\mathcal{G}}, s\in\widetilde{\mathcal{S}}}\left|\mathbb{E}_{P^N}[h(g(X),s(X),Y)]- \mathbb{E}_{P}[h(g(X),s(X),Y)]\right|\\
    & \leq & \sup_{g\in \widetilde{\mathcal{G}}_{\mathcal{X}_\epsilon},s\in \widetilde{\mathcal{S}}_{\mathcal{X}_\epsilon}}\left|\mathbb{E}_{P^N}[h(g(X),s(X),Y) \mathds{1}_{\mathcal{X}_\epsilon \times \mathcal{Y}_\epsilon}(X,Y)]\right. \\
    && \left. - \mathbb{E}_{P}[h(g(X),s(X),Y)\mathds{1}_{\mathcal{X}_\epsilon \times \mathcal{Y}_\epsilon}(X,Y)]\right|  + 2\epsilon \\
    &\leq & \sup_{k\in\{1,\dots,K\}}\sup_{g\in (\widetilde{\mathcal{G}}_{\mathcal{X}_\epsilon})_k^{\bar{\epsilon}}, s\in (\widetilde{\mathcal{S}}_{\mathcal{X}_\epsilon})_k^{\bar{\epsilon}}} \big\{ \mathbb{E}_{P^N}\left[ |h(g(X),s(X),Y)-h(g_k(X),s_k(X),Y)|\mathds{1}_{\mathcal{X}_\epsilon \times \mathcal{Y}_\epsilon}(X,Y)\right]\\
    &&+  \left|\mathbb{E}_{P^N}[h(g_k(X),s_k(X),Y)\mathds{1}_{\mathcal{X}_\epsilon \times \mathcal{Y}_\epsilon}(X,Y)]-\mathbb{E}_{P}[h(g_k(X),s_k(X),Y)\mathds{1}_{\mathcal{X}_\epsilon \times \mathcal{Y}_\epsilon}(X,Y)]\right|\\
    &&+ \mathbb{E}_{P}\left[|h(g_k(X),s_k(X),Y) - h(g(X),s(X),Y)|\mathds{1}_{\mathcal{X}_\epsilon \times \mathcal{Y}_\epsilon}(X,Y)\right]\big\} + 2\epsilon\\
    \\
    &\leq & 2\mathbb{E}_{P^N}[r(Y)\bar{\epsilon}^\nu] + \sup_{k\in\{1,\dots,K\}}
    \left|R_1(g_k,s_k,\epsilon,P^N)\right|
    + 2\mathbb{E}_{P}[r(Y)\bar{\epsilon}^\nu] + 2\epsilon\\
    &\leq & \sup_{k\in\{1,\dots,K\}}  \left|R_1(g_k,s_k,\epsilon,P^N)\right| +6\epsilon,
\end{eqnarray*}
where the last inequality holds by letting $\bar{\epsilon}\leq (\epsilon/\bar{r})^{1/\nu}$ and $\bar{r} := \max\{|r(x)|:x\in \mathcal{X}_\epsilon\}$.

\begin{eqnarray*}
    &&\inmat{Prob}\left( \sup_{g\in\widetilde{\mathcal{G}}} |H_N(g)-H(g)| \geq \delta \right)\\
    &\leq& \inmat{Prob} \left( \sup_{k\in\{1,\dots,K\}}  \left|R_1(g_k,s_k,\epsilon,P^N)\right|  + 6\epsilon  \geq \delta \right)\\
    &\leq& \sum_{k=1}^K \inmat{Prob} \bigg( \left|R_1(g_k,s_k,\epsilon,P^N)\right|  \geq \delta - 6\epsilon \bigg)\\
    &\leq& \sum_{k=1}^K \bigg[\inmat{Prob} \bigg(  R_1(g_k,s_k,\epsilon,P^N)  \geq \delta - 6\epsilon \bigg)  + \inmat{Prob} \bigg(  R_1(g_k,s_k,\epsilon,P^N)  \leq -( \delta  -6\epsilon) \bigg)\bigg]\\
    &\leq& 2\sum_{k=1}^K \left(e^{-N [I(\delta - 6\epsilon,g_k,s_k)\wedge I(-(\delta - 6\epsilon),g_k,s_k)]}\right),
\end{eqnarray*}
where the last inequality comes from Cram\'{e}r’s Large Deviation Theorem with $I(u,g_k,s_k): = \sup_{t\in \R} \{ut-\log m(t,g_k,s_k)\}$ and
\begin{eqnarray*}
    m(t,g_k,s_k) &:=& \mathbb{E}_{P}[e^{t[h(g_k(X),s_k(X),Y) -\mathbb{E}_{P}[h(g_k(X),s_k(X),Y)]]}].
\end{eqnarray*}
Note that $K$ depends on $\epsilon$ and thus depends on $\delta$.
By letting $C(\epsilon, \delta) = 2K $ and $\beta(\epsilon,\delta) = \min\limits_k\{I(\delta - 6\epsilon,g_k,s_k)\wedge I(-(\delta - 6\epsilon),g_k,s_k)\}$, we arrive at (\ref{eqn:lemma_saa_convergence_general}).
\hfill $\Box$

\subsubsection{Proof of Theorem \ref{thm: SAA-convergence_V}}
\label{proof:Thm4}
Let $H_N(g,s)$ and $H(g,s)$ be defined as in (\ref{eqn:saa-inner}) and (\ref{eqn:saa-inner-primal}).
Under Assumption \ref{assumption:convex_finite-minimum} (a),
both
$\vt_N(P^N)$
and $\vt(P)$
are bounded.
Moreover
\begin{eqnarray*}
&&\left|\vt_N(P^N)-\vt(P)\right|\nonumber\\ & = &
\left|\min_{g\in \widetilde{\mathcal{G}},s\in \widetilde{\mathcal{S}}} \left\{V\left( \mathbb{E}_{P^N}[h(g(X),s(X),Y)]\right)\right\}
-\min_{g\in \widetilde{\mathcal{G}},s\in \widetilde{\mathcal{S}}} V(\mathbb{E}_{P}[h(g(X),s(X),Y)])
\right|
\\
& \leq & \sup_{g\in \widetilde{\mathcal{G}},s\in \widetilde{\mathcal{S}}} \left\{|V\left( H_N(g,s)\right)
-V(H(g,s))
|
\right\}\\
&=&
\sup_{g\in \widetilde{\mathcal{G}},s\in \widetilde{\mathcal{S}}} \left[
|H_N(g,s)-H(g,s)| \left|\int_0^1
 V'\left(
 H(g,s)
 + \upsilon d_N \right)d\upsilon \right|
 \right],
\end{eqnarray*}
where $d_N := H_N(g,s)-H(g,s) = \mathbb{E}_{P^N}[h(g(X),s(X),Y)] - \mathbb{E}_{P}[h(g(X),s(X),Y)]$.

{
Part (i). Since $|V'(t)|$ is assumed to be upper bounded by $L_V>0$,
\begin{eqnarray*}
    |\vt_N(P^N)-\vt(P)| \leq \sup_{g\in \widetilde{\mathcal{G}},s\in \widetilde{\mathcal{S}}} \left[
L_V |H_N(g,s)-H(g,s)| \right].
\end{eqnarray*}
Let $\hat{\delta}:= \frac{\delta}{L_V}$.
By Lemma \ref{lemma:saa_convergence_general}, for any $\hat{\delta}$, there exist $\epsilon<\frac{\hat{\delta}}{6}$, $C(\epsilon, \hat{\delta})$, and $\beta(\epsilon,\hat{\delta})$, such that
\begin{eqnarray*}
        \inmat{Prob}\left(|\vt_N(P^N)-\vt(P)|\geq \delta\right) &\leq&
        \inmat{Prob}\left(\sup_{g\in \widetilde{\mathcal{G}},s\in \widetilde{\mathcal{S}}}|H_N(g,s)-H(g,s)|\geq \hat{\delta}\right)\\
        &\leq& C(\epsilon,\hat{\delta})e^{-N{\beta(\epsilon,\hat{\delta})}},
    \end{eqnarray*}
which indicates that (\ref{eqn:SAA_convergence_V1}) is satisfied.
}

Part (ii).
Under the
condition that $V(\cdot)= \ln \cdot$,
we can use the mean value theorem to obtain
\begin{eqnarray} \label{eqn:saa-proof-mvt}
    |\vt_N(P^N)-\vt(P)| \leq \sup_{g\in \widetilde{\mathcal{G}},s\in \widetilde{\mathcal{S}}} \left[
\left|
\frac{H_N(g,s)-H(g,s)}{H(g,s)+\upsilon(g,s) (H_N(g,s)-H(g,s))}
\right|
 \right],
\end{eqnarray}
where $\upsilon(g,s)\in [0,1]$ is a positive constant dependent on $(g,s)$.
If  $H_N(g,s)\leq H(g,s)$, then
\begin{eqnarray} \label{eqn:saa-proof-H(g)-1}
\left|
\frac{H_N(g,s)-H(g,s)}{H(g,s)+\upsilon(g,s) (H_N(g,s)-H(g,s))}
\right|
\geq \delta
&\Rightarrow&
\frac{H(g,s)-H_N(g,s)}{H_N(g,s)} \geq \delta \nonumber
\end{eqnarray}
and hence
\begin{eqnarray}
H(g,s)-H_N(g,s) \geq \frac{\delta}{1+\delta} H(g,s)
\geq \frac{\delta}{1+\delta} H^*,
\end{eqnarray}
where $H^*:= \inf\limits_{g\in \widetilde{\mathcal{G}},s \in \widetilde{\mathcal{S}}} \mathbb{E}_{P}[h(g(X), s(X),Y)]$ as defined above.
On the other hand, if
$H_N(g,s)> H(g,s)$, then
\begin{eqnarray} \label{eqn:saa-proof-H(g)-2}
\left|
\frac{H_N(g,s)-H(g,s)}{H(g,s)+\upsilon' (H_N(g,s)-H(g,s))}
\right|
\geq \delta
&\Rightarrow& \frac{H_N(g,s)-H(g,s)}{H(g,s)} \geq \delta
\nonumber
\end{eqnarray}
and subsequently
\begin{eqnarray}
H_N(g,s)-H(g,s) \geq \delta H(g,s) \geq \frac{\delta}{1+\delta} H^* .
\end{eqnarray}
Let $(g^*_N, s^*_N )\in \arg \sup\limits_{g\in \widetilde{\mathcal{G}},s \in \widetilde{\mathcal{S}}}
\left|
\frac{H_N(g,s)-H(g,s)}{H(g,s)+\upsilon' (H_N(g,s)-H(g,s))}
\right|$, which depends on the realization of $P^N$.
From (\ref{eqn:saa-proof-H(g)-1}) and (\ref{eqn:saa-proof-H(g)-2}), we have
\begin{eqnarray} \label{eqn:saa-proof-equivalent-of-sup}
&&    \sup_{g\in \widetilde{\mathcal{G}},s \in \widetilde{\mathcal{S}}}
\left|
\frac{H_N(g,s)-H(g,s)}{H(g,s)+\upsilon' (H_N(g,s)-H(g,s))}
\right|\geq \delta
\\
&&\quad\Rightarrow |H_N(g^*_N,s^*_N)-H(g^*_N,s^*_N)| \geq \frac{\delta}{1+\delta} H^* \nonumber
\end{eqnarray}
and hence
\begin{eqnarray}
     \sup_{g\in \widetilde{\mathcal{G}},s \in \widetilde{\mathcal{S}}}
\left|
H_N(g,s)-H(g,s) \right|\geq \frac{\delta}{1+\delta} H^*.
\end{eqnarray}
Then
\begin{eqnarray}
    &&\inmat{Prob}\left(|\vt_N(P^N)-\vt(P)|\geq \delta\right) \nonumber\\
    & \leq & \inmat{Prob}\left(
    \sup_{g\in \widetilde{\mathcal{G}},s\in \widetilde{\mathcal{S}}} \left[
\left|
\frac{H_N(g,s)-H(g,s)}{H(g,s)+\upsilon' (H_N(g,s)-H(g,s))}
\right|
 \right]
    \geq \delta\right) \nonumber\\
& \leq & \inmat{Prob}\left(
    \sup_{g\in \widetilde{\mathcal{G}},s\in \widetilde{\mathcal{S}}} \left[
\left|
H_N(g,s)-H(g,s)
\right|
 \right]
    \geq \delta/(1+ \delta) H^* \right) \nonumber\\
& \leq & C(\epsilon,\delta/(1+ \delta) H^*)e^{-N{\beta(\epsilon,\delta/(1+ \delta) H^*)}},
\end{eqnarray}
where
the first inequality
comes from (\ref{eqn:saa-proof-mvt}),
the second inequality
comes from (\ref{eqn:saa-proof-equivalent-of-sup}),
and the last inequality
is based on Lemma \ref{lemma:saa_convergence_general} with $\epsilon < \frac{\delta H^*/(1+\delta)}{6}$.
\hfill $\Box$

\subsection{Verification of key assumptions in applications}
\label{sec:vei-assumption-application}

{
For completeness,
we verify that the assumptions
required for the theoretical developments across Sections \ref{sec:contextual-risk-measure}-\ref{sec:DR-rkhs}
are satisfied
by the examples in
Section~\ref{sec:experiment}.
These assumptions are related to the cost function $c(z,y)$, including Assumption~\ref{assumption:well-defined_rho}(c), Assumption~\ref{assumption:convex_finite-minimum}, and Assumption~\ref{assump:saa-convergence-unbounded}~(a),~(d),~(e).
The other assumptions can be fulfilled by appropriately choosing the risk measures or the feasible set of decision policies in the setting of RKHS.
}

\subsubsection{Newsvendor problem}
\label{sec:veri-assum-newsvendor}

{
In the newsvendor problem, the cost function is
$
    c(z,y) = h(z-y)_+ + b(y-z)_+,
$
where $h,b\in \R_+$.

\begin{itemize}
    \item[(i)] \textbf{Assumption \ref{assumption:well-defined_rho}(c)}.  Observe that
\begin{eqnarray*}
    |c(z,y)| \leq (h+b)|z-y| \leq (h+b)(|z|+1)(|y|+1).
\end{eqnarray*}
By setting $\eta(z):= (h+b)(|z|+1)$ and $\iota=1$, we verify the condition.

\item[(ii)] \textbf{Assumption~\ref{assumption:convex_finite-minimum}(a)}.
For any fixed $\hat x \in \mathcal{X}$ and a closed neighborhood $U$ of $\hat{x}$ relative to $\mathcal{X}$, let
$\kappa:=
\sup_{x\in U}\rho_{ P_{Y|X=x}}^{(2)}(bY) + 1$.
Since  $\rho_{P_{Y|X=x}}^{(2)}(bY)$ is continuous in $x$ as ensured by
Proposition~\ref{prop:random_variable_inner_risk},
then $\kappa< \infty$.
Observe that $\sup_{x\in U}\rho_{ P_{Y|X=x}}^{(2)}(c(0,Y))
=\sup_{x\in U}\rho_{ P_{Y|X=x}}^{(2)}(bY).
$
Then
$0\in L_x(\kappa):=\left\{ z\in \mathcal{Z}: \rho_{ P_{Y|X=x}}^{(2)}(c(z,Y))\leq \kappa \right\}$ for all $x\in U$ and thus $L_x(\kappa)\neq \emptyset$.
On the other hand,
{
observe that the cost function satisfies $c(z,y)\geq h(z-y)$ for every $z,y$,
}
then
$$
\rho_{ P_{Y|X=x}}^{(2)}(c(z,Y)) \geq \rho_{ P_{Y|X=x}}^{(2)}(h(z-Y)) \geq \mathbb{E}_{P_{Y|X=x}}[h(z-Y)] = hz - h\mathbb{E}_{P_{Y|X=x}}[Y],
$$
{
where the first inequality follows from the monotonicity of $\rho_{P_{Y|X}}^{(2)}$, and the second inequality comes from the fact that
$\rho_{ P_{Y|X=x}}^{(2)}(Y) \geq \mathbb{E}[Y]$ for any normalized convex risk measure $\rho_{ P_{Y|X=x}}^{(2)}$; see,
 e.g., \cite[Corollary 4.65]{follmer2011stochastic}.
}
By Assumption \ref{assumption:well-defined_rho}(d), $\mathbb{E}_{P_{Y|X=x}}[Y]$ is continuous in $x$ and thus there exists a $M_U< \infty$ such that $\sup_{x\in U}\mathbb{E}_{P_{Y|X=x}}[Y]\leq M_U$.
Summarizing the discussions above, we conclude that
$$
\kappa \geq \rho_{ P_{Y|X=x}}^{(2)}(c(z,Y)) \geq hz -hM_U, \forall x\in U, z\in  L_x(\kappa),
$$
and thus $z\leq M_U + \frac{\kappa}{h}$. Assumption \ref{assumption:convex_finite-minimum}(a) therefore holds with $\hat{Z}:= [0, M_U + \frac{\kappa}{h}]$.

\textbf{Assumption \ref{assumption:convex_finite-minimum}(b).}
Gotoh and Takano \cite{gotoh2007newsvendor} show that when $P_{Y|X}$ is a continuous distribution and $\rho_{ P_{Y|X=x}}^{(2)}=\text{CVaR}_{1-\alpha}$, the  optimal solution to problem \eqref{eqn:conditional_risk_min} has a closed-form
\begin{eqnarray*}
    z^{*}(x) &=& \frac{h}{h+b}F_{Y|X=x}^{-1}\left(\frac{b(1-\beta)}{h+b}\right) + \frac{b}{b+h}F_{Y|X=x}^{-1}\left(\frac{h\beta+b}{b+h}\right),
\end{eqnarray*}
which ensures the uniqueness of the optimal solution.
Problem \eqref{eqn:conditional_risk_min} with entropic risk measure can be equivalently written as
$$\min_{z\in\mathcal{Z}} \mathbb{E}_{P_{Y|X=x}}[\exp \left( \gamma (h(z-Y)_+ + b(Y- z)_+) \right)].$$
Since the objective function is strictly convex, the optimal solution is unique.

\item[(iii)] \textbf{Assumption \ref{assump:saa-convergence-unbounded}(a)}.
By the setting of risk-averse contextual optimization problem \eqref{eqn:ex-post-interchange} with expected-CVaR structure, the $h(z,s,y)$ function in problem \eqref{eqn:saa-form-rkhs} takes the form of
\begin{eqnarray*}
    h(z,s,y) = s+ \frac{1}{1-\alpha}\left[h(z-y)_+ + b(y - z)_+-s\right]_+.
\end{eqnarray*}
Then
\begin{eqnarray*}
    |h(z(x),s(x),y)| &=& \left|s(x)+ \frac{1}{1-\alpha}\left[h(z(x)-y)_+ + b(y - z(x))_+-s(x)\right]_+\right|\\
    &\leq& |s(x)|+ \frac{h}{1-\alpha}|z(x)-y| + \frac{b}{1-\alpha}|y - z(x)| + \frac{1}{1-\alpha}|s(x)|\\
    &\leq & \frac{2-\alpha}{1-\alpha}|s(x)| + \frac{h+b}{1-\alpha}(|z(x)|+|y|).
\end{eqnarray*}
Note that
\begin{eqnarray*}
    |z(x)| \leq |\langle z,k(\cdot,x) \rangle|\leq \|z\|_\mathcal{H} \|k(\cdot,x)\|_\mathcal{H},\quad
    |s(x)| \leq |\langle s,k(\cdot,x) \rangle| \leq \|s\|_\mathcal{H} \|k(\cdot,x)\|_\mathcal{H}.
\end{eqnarray*}
If there exist positive constants $\beta_1,\beta_2$ such that $\|z\|_\mathcal{H}\leq \beta_1$ and $\|s\|_\mathcal{H} \leq \beta_2$ as
stated in Assumption \ref{assump:saa-convergence-unbounded}(b), then, for certain kernels with bounded $\|k(\cdot, x)\|_\mathcal{H}$, we can obtain an upper bound on $|h(z(x),s(x),y)|$. Consider,
for example,  a radial basis function (RBF) kernel that depends only on the distance between inputs,
i.e., $k(x',x'') = \tilde{k}(\|x' - x''\|)$ for some function $\tilde{k}:\R \to \R$,
such as the Gaussian kernel defined by $k(x',x''):= \exp\left(-\frac{\|x'-x''\|^2}{2\sigma^2} \right)$, for which $\|k(\cdot,x)\|_\mathcal{H}^2 = \tilde{k}(0)= 1$.
In this case,
we can set
\begin{eqnarray} \label{eqn:newsvendor-4a-1}
    \phi(x,y) &:=& \frac{2-\alpha}{1-\alpha}\beta_2 + \frac{h+b}{1-\alpha}(\beta_1 + |y|) + \|x\| \geq \left|h(z(x),s(x),y)\right|,
\end{eqnarray}
where $\int_{\mathcal{X}\times \mathcal{Y}} \phi(x,y) P(dx,dy) < \infty$
for all
$P\in \mathcal{M}_1^1$, and $\phi(x,y)\rightarrow +\infty$ as $\|(x,y)\|\rightarrow \infty$.
Likewise,
in the
risk-averse contextual optimization problem \eqref{eqn:ex_ante} with entropic risk measure where
$
  h(z,s,y) = e^{\gamma(h(z-y)_+ + b(y - z)_+)},
$
we have
\begin{eqnarray*}
    |h(z(x),s(x),y)| &=& \left|e^{\gamma(h(z(x)-y)_+ + b(y - z(x))_+)}\right| \\
    &\leq&  e^{\gamma(h+b)|z(x)-y|}\leq  e^{\gamma(h+b)|z(x)|+|y|}.
\end{eqnarray*}
Note that $
    |z(x)| \leq \|z\|_\mathcal{H}\|k(\cdot,x)\|_\mathcal{H}.
$
If there is a positive number $\beta_1$ such that $\|z\|_\mathcal{H}\leq \beta_1$ as
specified in Assumption \ref{assump:saa-convergence-unbounded}(b) and $k$ is a RBF kernel, such as Gaussian kernel with $\|k(\cdot,x)\|_\mathcal{H}^2 = 1$,
then
we can set
\begin{eqnarray*}
    \phi(x,y) &:=& e^{\gamma(h+b)\beta_1} +e^{\gamma(h+b) |y|} + \|x\|
    \geq h(z(x),s(x),y).
\end{eqnarray*}
To ease the exposition, we
consider a special case
that $(X,Y)$ follows a standard multivariate normal distribution.
Then
\begin{eqnarray} \label{eqn:newsvendor-4a-2}
    \mathbb{E}_P\left[ e^{\gamma(h+b) |Y|}\right] &=& \int_{-\infty}^{+\infty} \exp\left(\gamma(h+b) |y|\right)  \frac{1}{\sqrt{2\pi}} \exp\left(-\frac{y^2}{2}\right) dy \nonumber\\
    &=& 2\int_0^{+\infty} \exp\left(\gamma(h+b) y\right)  \frac{1}{\sqrt{2\pi}}\exp\left(-\frac{y^2}{2}\right) dy\nonumber\\
    &=& 2 \exp\left( \frac{\gamma^2(h+b)^2}{2} \right) \int_0^{+\infty} \frac{1}{\sqrt{2\pi}}\exp\left(-\frac{(y-\gamma(h+b))^2}{2}\right) dy\nonumber\\
    &=& 2 \exp\left( \frac{\gamma^2(h+b)^2}{2} \right) \Phi(\gamma(h+b))<\infty,
\end{eqnarray}
where $\Phi(\cdot)$ is the cumulative distribution function of a standard normal distribution.
Moreover, $\mathbb{E}_P[\|X\|] \leq \mathbb{E}_P\left[ \|X\|^2 \right]^{1/2} < \infty$.
Thus, we can conclude that  $\int_{\mathcal{X}\times \mathcal{Y}} \phi(x,y) P(dxdy) < \infty$.
The same conclusion can be drawn
for general multivariate normal distributions.

\item[(iv)] \textbf{Assumption \ref{assump:saa-convergence-unbounded}(d)}. Observe that
$$
h(z,s,y) = s+ \frac{1}{1-\alpha}\left[h(z-y)_+ + b(y - z)_+-s\right]_+
$$
is Lipschitz continuous in $(z,s)$ and therefore satisfies the assumption.
Likewise, $h(z,s,y) = e^{\gamma(h(z-y)_+ + b(y - z)_+)}$ is Lipschitz continuous on any compact subset of $\mathcal{Z}\times\mathcal{S}\times \mathcal{Y}$ which fulfills the assumption.

\item[(v)] \textbf{Assumption \ref{assump:saa-convergence-unbounded}(e)}.
This assumption requires
\begin{eqnarray}
\mathbb{E}_{P}\big[h(z(X),s(X),Y)\big] < \infty,\quad \mathbb{E}_{P}\big[e^{th(z(X),s(X),Y)}\big]< \infty. \label{eqn:newsvendor-4e}
\end{eqnarray}
Based on the previous analysis for Assumption \ref{assump:saa-convergence-unbounded}(a), we
conclude that, for problem \eqref{eqn:ex-post-interchange} with the expected CVaR and a Gaussian kernel $k$, we will have
\begin{eqnarray*}
    |h(z(x),s(x),y)| \leq \frac{2-\alpha}{1-\alpha}\beta_2 + \frac{h+b}{1-\alpha}(\beta_1+ |y|).
\end{eqnarray*}
Moreover, following
a proof analogous to that of
\eqref{eqn:newsvendor-4a-1} and \eqref{eqn:newsvendor-4a-2}, we can establish \eqref{eqn:newsvendor-4e}
when $Y$ follows a  normal distribution.
In contrast, for problem \eqref{eqn:ex-post-interchange} with the entropic risk measure where $h(z,s,y) = e^{\gamma(h(z-y)_+ + b(y - z)_+)}$,
\begin{eqnarray*}
    \mathbb{E}_{P}\big[e^{th(z(X),s(X),Y)}\big] &=& \int e^{t\exp\left(\gamma(h(z(x)-y)_+ + b(y - z(x))_+)\right)} P(dx,dy)\\
    & \geq & \int e^{t\exp\left(\gamma b(y - z(x))\right)} P(dx,dy)\\
    & \geq & \int e^{t\exp\left(\gamma b(y - \beta_1)\right)} P_Y(dy),
\end{eqnarray*}
where the second inequality comes from the fact that $\|z(x)\|\leq \|z\|_\mathcal{H} \|k(\cdot,x)\|_\mathcal{H}\leq \beta_1$.
Therefore, the exponential moment $\mathbb{E}_{P}\big[e^{th(z(X),s(X),Y)}\big]$ diverges to infinity under a normal distribution for \(Y\), e.g., for a standard normal distribution,
\begin{eqnarray}
    \mathbb{E}_{P}\big[e^{th(z(X),s(X),Y)}\big] \geq \int_{-\infty}^{+\infty} e^{t\exp\left(\gamma b(y - \beta_1)\right)} \frac{1}{\sqrt{2\pi}} e^{-\frac{y^2}{2}} dy = \infty,
\end{eqnarray}
since $\exp(\gamma b y)$ dominates $-\frac{y^2}{2}$ as $y\rightarrow + \infty$.
To avoid the issue in
the numerical experiments,
we adopt a truncated normal distribution for $Y$.

\end{itemize}

\subsubsection{Portfolio selection}
\label{sec:veri-assum-portfolio}

{
In the portfolio selection problem, we focus on the contextual mean-expected-CVaR model \eqref{eqm:CMEC-portfolio}, where the cost function is $c(z,y) = -y^Tz$ and the feasible set is $\mathcal{Z}:=\{z\in \R^{d_z}: \sum_{j=1}^{d_z} z^{(j)}\leq 1, z\geq 0\}$.
}

\begin{itemize}
    \item[(i)] \textbf{Assumption \ref{assumption:well-defined_rho}(c).} Setting $\eta(z):= \|z\|$ and $\iota=1$, we have
    \begin{eqnarray}
        |c(z,y)|:= \left|-y^T z \right| \leq \|z\| \|y\| \leq \eta(z) (\|y\|^\iota +1).
    \end{eqnarray}

    \item[(ii)]  \textbf{Assumption  \ref{assumption:convex_finite-minimum}(a).} Since the feasible set $\mathcal{Z}$ is compact, Assumption \ref{assumption:convex_finite-minimum}(a) is automatically satisfied.

    \textbf{Assumption  \ref{assumption:convex_finite-minimum}(b).}
    This assumption requires additional conditions on the conditional distribution $P_{Y|X}$.
    We concentrate on a special case where for any $x\in \mathcal{X}$, $P_{Y|X=x}$ follows a normal distribution with positive definite covariance matrix $\Sigma_{Y|X=x}$.
    By \cite[Example 2.18]{mcneil2015quantitative}, the objective function {of the mean-CVaR problem \eqref{eqn:mean-cvar} under normal distribution $P_{Y|X=x}$} can be reformulated  as
\begin{eqnarray}
    f(z)&=&\inmat{CVaR}_{\beta,P_{Y|X=x}}(-Y^T z) - \eta \mathbb{E}_{P_{Y|X=x}}[Y^T z]\nonumber\\
    &=& -\mathbb{E}_{P_{Y|X=x}}[Y^Tz] + \frac{\phi(\Phi^{-1}(\beta))}{1-\beta} \sqrt{z^T\Sigma_{Y|X=x}z} - \eta \mathbb{E}_{P_{Y|X=x}}[Y^T z]\nonumber\\
    &=& \frac{\phi(\Phi^{-1}(\beta))}{1-\beta} \sqrt{z^T\Sigma_{Y|X=x}z} - (1+\eta) \mathbb{E}_{P_{Y|X=x}}[Y^T z],\label{eqn:portfolio-obj-equi}
\end{eqnarray}
where $\phi(\cdot)$ and $\Phi(\cdot)$ are probability density function and cumulative distribution function of a standard normal distribution.
Let $a:= \frac{\phi(\Phi^{-1}(\beta))}{1-\beta}$ and $b:=(1+\eta) \mathbb{E}_{P_{Y|X=x}}[Y]$, and $\|z\|_\Sigma := \sqrt{z^T\Sigma_{Y|X=x}z}$ be a norm of $\R^{d_z}$.
For $z,z'\in \mathcal{Z}$, we have
\begin{eqnarray*}
    f(t z + (1-t) z') &=& a \|t z + (1-t) z'\| - b^T (t z + (1-t) z')\\
    &\leq& t a \| z \| - t b^T z + (1- t) a\| z' \| - (1 - t) b^T z'\\
    & = & t f(z) + (1-t) f(z'),\quad
    \forall t\in [0,1],
\end{eqnarray*}
where equation holds only if $z' = \lambda z$ or $z = \lambda z'$ for some $\lambda \geq 0$.
Therefore, the objective function is not strictly convex on $\mathcal{Z}$. However, it is strictly convex on
    $\Delta: = \{z\in \R^{d_z}: \sum_{j=1}^{d_z} z^{(j)}= 1, z\geq 0\}$ as it excludes the case that $z' = \lambda z$ or $z = \lambda z'$ for any $\lambda \geq 0$.
    Note that objective function is positively homogeneous for $\lambda \in [0,1]$, i.e.
    \begin{eqnarray*}
        &&\inmat{CVaR}_{\beta,P_{Y|X=x}}(-Y^T(\lambda z)) - \eta \mathbb{E}_{P_{Y|X=x}}[Y^T(\lambda z)] \\
        && \quad = \lambda \left( \inmat{CVaR}_{\beta,P_{Y|X=x}}(-Y^T z) - \eta \mathbb{E}_{P_{Y|X=x}}[Y^T z] \right).
    \end{eqnarray*}
    Let
    \begin{eqnarray} \label{eqn:portfolio-obj-delta}
    \vartheta:=\min_{z\in \Delta} f(z) :=
    \inmat{CVaR}_{\beta,P_{Y|X=x}}(-Y^T z) - \eta \mathbb{E}_{P_{Y|X=x}}[Y^T z],
\end{eqnarray}
and
$z^*\in \Delta$
be an optimal solution.
    We consider the following three cases.

    \textbf{Case 1.}
$\vartheta > 0$.
Then the minimum value of $f(z)$
over $\mathcal{Z}$ is zero, i.e.,
\begin{eqnarray} \label{eqn:portfolio-obj}
    \min_{z\in \mathcal{Z}} \inmat{CVaR}_{\beta,P_{Y|X=x}}(-Y^T z) - \eta \mathbb{E}_{P_{Y|X=x}}[Y^T z] = 0,
\end{eqnarray}
with unique solution $0 \cdot z^* = 0$.

\textbf{Case 2.}
$\vartheta = 0$.
Then every $\lambda z^*$ with $\lambda \in [0,1]$ is an optimal solution of problem \eqref{eqn:portfolio-obj}.
In this case, the solution is not unique, and we will avoid this case in the numerical experiment.

\textbf{Case 3.}
$\vartheta < 0$.
Due to the positive homogeneity of $f(z)$,
the optimum must be attained on $\Delta$,
and $z^*$ is the optimal solution to \eqref{eqn:portfolio-obj}.

In the numerical experiment, we will adopt a distribution
{and proper parameter values
$\eta$ and $\beta$ such that
\begin{eqnarray*}
    \frac{\phi(\Phi^{-1}(\beta))}{(1-\beta)(1+\eta)}  <  \min_{z\in \mathcal{Z}}\frac{\mathbb{E}_{P_{Y|X=x}}[Y^T z]}{\sqrt{z^T\Sigma_{Y|X=x}z}},
\end{eqnarray*}
}
which leads to \textbf{case 3}.

    \item[(iii)] \textbf{Assumption \ref{assump:saa-convergence-unbounded}(a)}. By the setting of problem \eqref{eqn:ex-post-interchange} with expected-CVaR structure,
\begin{eqnarray} \label{eqn:portfolio-h}
    h(z,t,y) = t+ \frac{1}{1-\beta}\left[ -y^Tz-t\right]_+ - \eta y^Tz.
\end{eqnarray}
Note that, if $k$ is a Gaussian kernel with $\|k(\cdot,x)\|_\mathcal{H}^2=1$ for all $x$, then
\begin{eqnarray}
    |h(z(x),t(x),y)| &\leq&  |t(x)|+ \frac{1}{1-\alpha}\left( \left| y^Tz(x) \right| + |t(x)|\right) + \eta \left| y^Tz(x) \right|\nonumber\\
    &\leq & \frac{2-\alpha}{1-\alpha} \|t\|_\mathcal{H} \|k(\cdot,x)\|_\mathcal{H} + \frac{1+\eta -\eta\alpha}{1-\alpha} \|y\|_\infty \|z\|_1\nonumber\\
    &\leq & \frac{2-\alpha}{1-\alpha} \beta_2+ \frac{1+\eta -\eta\alpha}{1-\alpha} \|y\| + \|x\|, \label{eqn:portfolio-h-bound}
\end{eqnarray}
where the third inequality follows from the fact that $\|t\|_\mathcal{H}\leq \beta_2$ under Assumption \ref{assump:saa-convergence-unbounded}(b),  $\|z\|_1\leq 1$ for any feasible solution $z$, and
$\|y\|\geq \|y\|_\infty$. Consequently, we can set
$
    \phi(x,y) := \frac{2-\alpha}{1-\alpha} \beta_2+ \frac{1+\eta -\eta\alpha}{1-\alpha} \|y\| + \|x\|,
$
where $\int_{\mathcal{X}\times \mathcal{Y}} \phi(x,y) P(dxdy) < \infty$
for
$P\in \mathcal{M}_{d_z}^1$, and $\phi(x,y)\rightarrow +\infty$ as $\|(x,y)\|\rightarrow \infty$.

\textbf{Assumption \ref{assump:saa-convergence-unbounded}(d)}.
The function $h(z,t,y)$ in \eqref{eqn:portfolio-h} is Lipschitz on $\mathcal{Z}\times \mathcal{T}$, and therefore satisfies the assumption.

\textbf{Assumption \ref{assump:saa-convergence-unbounded}(e)}.
This assumption requires \eqref{eqn:newsvendor-4e}.
From \eqref{eqn:portfolio-h-bound}, for the Gaussian kernel we have
\begin{eqnarray*}
    |h(z(x),s(x),y)| \leq \frac{2-\alpha}{1-\alpha}\beta_2 + \frac{1+\eta -\eta\alpha}{1-\alpha}\|y\|.
\end{eqnarray*}
For simplicity, we consider a special case that $Y$ follows a standard normal distribution $\mathcal{N}(0,I^{d_y\times d_y})$. Then $\mathbb{E}[\|Y\|] \leq  \mathbb{E}[\|Y\|^2]^{1/2} = \text{tr}(I^{d_y\times d_y})^{1/2} = \sqrt{d_y}$ and, similar to \eqref{eqn:newsvendor-4a-2}, we can derive
\begin{eqnarray*}
    \mathbb{E}\left[ e^{\|Y\|} \right] &\leq& \mathbb{E}\left[ e^{ \sum_{i=1}^{d_y}|Y_i|} \right] = \prod_{i=1}^{d_y} \mathbb{E}\left[e^{|Y_i|}\right] = 2e^{1/2} \Phi(1),
\end{eqnarray*}
where the first equality follows from the fact that the $Y_i$ are i.i.d.
The assumption is fulfilled when $Y$ follows a standard normal distribution.
The same conclusion can be drawn
for general multivariate normal distributions.

\end{itemize}

}

\noindent
\textbf{Acknowledgments.} The authors would like to thank the three anonymous referees and the Associate Editor for valuable comments which helped them significantly improve the presentation of the paper.

\section*{Declaration of competing interest}
The authors declare that they have no known competing financial interests or personal relationships that could have appeared to influence the work reported in this paper.

\bibliographystyle{plainnat}
\bibliography{literature}

\end{document}